\documentclass[12pt,a4]{amsart}
\usepackage[a4paper, left=28mm, right=28mm, top=28mm, bottom=34mm]{geometry}
\usepackage{amsmath}
\usepackage{amssymb}
\usepackage{amsthm}
\usepackage{comment}
\usepackage{adjustbox}

\usepackage[all]{xy,xypic}
\usepackage{tikz}
\usetikzlibrary{matrix, arrows}
\usepackage{color, pstricks}
\usepackage{enumerate}
\usepackage{mathrsfs}
\usepackage{amsmath,amsthm, amssymb, amsgen,amsxtra,amsfonts, amsbsy} 
\usepackage{verbatim}
\usepackage{tensor}
\newcommand{\mathds}[1]{{\mathbb #1}}

\begin{document}
%
%
%
\theoremstyle{definition}
\newtheorem{Definition}{Definition}[section]
\newtheorem*{Definitionx}{Definition}
\newtheorem{Convention}{Definition}[section]
\newtheorem{Construction}{Construction}[section]
\newtheorem{Example}[Definition]{Example}
\newtheorem{Examples}[Definition]{Examples}
\newtheorem{Remark}[Definition]{Remark}
\newtheorem*{Remarkx}{Remark}
\newtheorem{Remarks}[Definition]{Remarks}
\newtheorem{Caution}[Definition]{Caution}
\newtheorem{Conjecture}[Definition]{Conjecture}
\newtheorem*{Conjecturex}{Conjecture}
\newtheorem{Question}[Definition]{Question}
\newtheorem{Questions}[Definition]{Questions}
\newtheorem*{Acknowledgements}{Acknowledgements}
\newtheorem*{Organization}{Organization}
\newtheorem*{Disclaimer}{Disclaimer}
\theoremstyle{plain}
\newtheorem{Theorem}[Definition]{Theorem}
\newtheorem*{Theoremx}{Theorem}
\newtheorem{Theoremy}{Theorem}
\newtheorem{Proposition}[Definition]{Proposition}
\newtheorem*{Propositionx}{Proposition}
\newtheorem{Lemma}[Definition]{Lemma}
\newtheorem{Corollary}[Definition]{Corollary}
\newtheorem*{Corollaryx}{Corollary}
\newtheorem{Fact}[Definition]{Fact}
\newtheorem{Facts}[Definition]{Facts}
\newtheoremstyle{voiditstyle}{3pt}{3pt}{\itshape}{\parindent}%
{\bfseries}{.}{ }{\thmnote{#3}}%
\theoremstyle{voiditstyle}
\newtheorem*{VoidItalic}{}
\newtheoremstyle{voidromstyle}{3pt}{3pt}{\rm}{\parindent}%
{\bfseries}{.}{ }{\thmnote{#3}}%
\theoremstyle{voidromstyle}
\newtheorem*{VoidRoman}{}

%
\newcommand{\prf}{\par\noindent{\sc Proof.}\quad}
\newcommand{\blowup}{\rule[-3mm]{0mm}{0mm}}
\newcommand{\cal}{\mathcal}
\newcommand{\Aff}{{\mathds{A}}}
\newcommand{\BB}{{\mathds{B}}}
\newcommand{\CC}{{\mathds{C}}}
\newcommand{\DD}{{\mathds{D}}}
\newcommand{\EE}{{\mathds{E}}}
\newcommand{\FF}{{\mathds{F}}}
\newcommand{\GG}{{\mathds{G}}}
\newcommand{\HH}{{\mathds{H}}}
\newcommand{\KK}{{\mathds{K}}}
\newcommand{\NN}{{\mathds{N}}}
\newcommand{\ZZ}{{\mathds{Z}}}
\newcommand{\PP}{{\mathds{P}}}
\newcommand{\QQ}{{\mathds{Q}}}
\newcommand{\RR}{{\mathds{R}}}
\newcommand{\Liea}{{\mathfrak a}}
\newcommand{\Lieb}{{\mathfrak b}}
\newcommand{\Lieg}{{\mathfrak g}}
\newcommand{\Liem}{{\mathfrak m}}
\newcommand{\ideala}{{\mathfrak a}}
\newcommand{\idealb}{{\mathfrak b}}
\newcommand{\idealg}{{\mathfrak g}}
\newcommand{\idealm}{{\mathfrak m}}
\newcommand{\idealp}{{\mathfrak p}}
\newcommand{\idealq}{{\mathfrak q}}
\newcommand{\idealI}{{\cal I}}
\newcommand{\lin}{\sim}
\newcommand{\num}{\equiv}
\newcommand{\dual}{\ast}
 \newcommand{\iso}{\cong}
\newcommand{\homeo}{\approx}
\newcommand{\mm}{{\mathfrak m}}
\newcommand{\pp}{{\mathfrak p}}
\newcommand{\qq}{{\mathfrak q}}
\newcommand{\rr}{{\mathfrak r}}
\newcommand{\pP}{{\mathfrak P}}
\newcommand{\qQ}{{\mathfrak Q}}
\newcommand{\rR}{{\mathfrak R}}
%
%
\newcommand{\OO}{{\cal O}}
\newcommand{\calA}{{\cal A}}
\newcommand{\calO}{{\cal O}}
\newcommand{\calU}{{\cal U}}
\newcommand{\numero}{{n$^{\rm o}\:$}}
\newcommand{\mf}[1]{\mathfrak{#1}}
\newcommand{\mc}[1]{\mathcal{#1}}
\newcommand{\MazurL}{{\cal L}}
\newcommand{\into}{{\hookrightarrow}}
\newcommand{\onto}{{\twoheadrightarrow}}
\newcommand{\Spec}{{\rm Spec}\:}
\newcommand{\BigSpec}{{\rm\bf Spec}\:}
\newcommand{\Spf}{{\rm Spf}\:}
\newcommand{\Proj}{{\rm Proj}\:}
\newcommand{\Pic}{{\rm Pic }}
\newcommand{\Alb}{{\rm Alb}}
\newcommand{\Br}{{\rm Br}}
\newcommand{\NS}{{\rm NS}}
\newcommand{\MF}{{\rm MF}}
\newcommand{\Fil}{{\rm Fil}}
\newcommand{\Sym}{{\mathfrak S}}
\newcommand{\Aut}{{\rm Aut}}
\newcommand{\Autp}{{\rm Aut}^p}
\newcommand{\End}{{\rm End}}
\newcommand{\Hom}{{\rm Hom}}
\newcommand{\Ext}{{\rm Ext}}
\newcommand{\ord}{{\rm ord}}
\newcommand{\coker}{{\rm coker}\,}
\newcommand{\divisor}{{\rm div}}
\newcommand{\Def}{{\rm Def}}
\newcommand{\et}{{\rm \acute{e}t}}
\newcommand{\fppf}{{\rm fppf}}
\newcommand{\loc}{{\rm loc}}
\newcommand{\ab}{{\rm ab}}
\newcommand{\piet}{{\pi_1^{\rm \acute{e}t}}}
\newcommand{\pitop}{{\pi_1^{\rm top}}}
\newcommand{\pietloc}{{\pi_{\rm loc}^{\rm \acute{e}t}}}
\newcommand{\piN}{{\pi^{\rm N}_1}}
\newcommand{\piNloc}{{\pi_{\rm loc}^{\rm N}}}
\newcommand{\Het}[1]{{H_{\rm \acute{e}t}^{{#1}}}}
\newcommand{\Hfl}[1]{{H_{\rm fl}^{{#1}}}}
\newcommand{\Hcris}[1]{{H_{\rm cris}^{{#1}}}}
\newcommand{\Hrig}[1]{{H_{\rm rig}^{{#1}}}}
\newcommand{\Hlogcris}[1]{{H_{\rm log-cris}^{{#1}}}}
\newcommand{\HdR}[1]{{H_{\rm dR}^{{#1}}}}
\newcommand{\hdR}[1]{{h_{\rm dR}^{{#1}}}}
\newcommand{\defin}[1]{{\bf #1}}
\newcommand{\oX}{\cal{X}}
\newcommand{\oA}{\cal{A}}
\newcommand{\oY}{\cal{Y}}
\newcommand{\calC}{{\cal{C}}}
\newcommand{\calL}{{\cal{L}}}
\newcommand{\bmu}{\boldsymbol{\mu}}
\newcommand{\balpha}{\boldsymbol{\alpha}}
\newcommand{\bM}{{\mathbf{M}}}
\newcommand{\BD}{{\mathbf{BD}}}
\newcommand{\BT}{{\mathbf{BT}}}
\newcommand{\BI}{{\mathbf{BI}}}
\newcommand{\BO}{{\mathbf{BO}}}
\newcommand{\SL}{{\mathbf{SL}}}
\newcommand{\MC}{{\mathbf{MC}}}
\newcommand{\GL}{{\mathbf{GL}}}

\title[Tate conjectures and abeloid varieties]{$\mathbf{p}$-adic Tate conjectures and abeloid varieties}

\author{Oliver Gregory}
\address{TU M\"unchen, Zentrum Mathematik - M11, Boltzmannstr. 3, 85748 Garching bei M\"unchen, Germany}
\email{gregory@ma.tum.de}

\author{Christian Liedtke}
\address{TU M\"unchen, Zentrum Mathematik - M11, Boltzmannstr. 3, 85748 Garching bei M\"unchen, Germany}
\email{liedtke@ma.tum.de}

\date{October 30, 2019}

\begin{abstract}
We explore Tate-type conjectures over $p$-adic fields, especially a conjecture
of Raskind \cite{Raskind} that predicts the surjectivity of 
\begin{equation*}  
\left({\rm NS}(X_{\overline{K}}) \otimes_{\ZZ}\QQ_p \right)^{G_K}
  \longrightarrow
   \Het{2}\left(X_{\overline{K}},\QQ_p(1)\right)^{G_K}
\end{equation*}
if $X$ is smooth and projective over a $p$-adic field $K$ and has totally
degenerate reduction.
Sometimes, this is related to $p$-adic uniformisation. For abelian varieties, Raskind's conjecture is equivalent 
to the question whether
$$
 \Hom(A,B)\otimes\QQ_p \,\to\, \Hom_{G_K}(V_p(A),V_p(B))
$$
is surjective if $A$ and $B$ are abeloid varieties over a $p$-adic field.

Using $p$-adic Hodge theory and Fontaine's functors, we reformulate both problems
into questions about the interplay of $\QQ$- versus $\QQ_p$-structures
inside filtered $(\varphi,N)$-modules.
Finally, we disprove all of these conjectures and questions
by showing that they can fail for algebraisable abeloid surfaces, that is,
for abelian surfaces with totally degenerate reduction.
\end{abstract}

\maketitle

section{Introduction}\label{intro}
Let $F$ be a field and let $G_F:={\rm Gal}(F^{\rm sep}/F)$ be the Galois group 
of a separable closure $F^{\rm sep}$ of $F$.
If $X$ is a smooth and proper variety over $F$ and $\ell$ is a prime
different from the characteristic $p$ of $F$, then the first Chern class map gives rise to  
an injective homomorphism of $\QQ_\ell$-vector spaces
\begin{equation*}\tag{$\star$} 
 \label{tatehom}
 c_{1,\ell}:{\rm NS}(X)\otimes_{\ZZ}\QQ_\ell
  \longrightarrow
   \Het{2}\left(X_{F^{\rm sep}},\QQ_\ell(1)\right)^{G_F}.
\end{equation*}
This is far from being an isomorphism in general. For example, if one takes $F$ to be a separably closed field of characteristic $p=0$, then the image
of $c_{1,\ell}$ is a proper subspace for any smooth and proper variety $X$ with $h^2(\OO_X)\neq0$. 

\subsection{The classical Tate conjecture}
However, the {\em Tate conjecture (for divisors)} predicts that  \eqref{tatehom} is surjective if $F$ is finitely generated over its 
prime field (for example, if $F$ is a number field or a finite field), see \cite{Tate1, Tate2}. 
This conjecture
is known to be true if $X$ is an abelian variety \cite{FaltingsFields, TateAbelian, Zar75},
if $X$ is a hyperk\"ahler variety and $F$ is finitely generated over $\QQ$ 
\cite{Andre, Tankeev}, as well as if $X$ is a K3 surface and $F$ is a finite field
\cite{Charles, KM, Madapusi, Maulik, Nygaard, NO}.
In \cite{Moonen}, it has been established for surfaces $X$ with $h^2(\OO_X)=1$
if $F$ is finitely generated over $\QQ$ and under the assumption that 
the Hodge structure on $H^2(X_\CC,\QQ)$ varies sufficiently 
non-trivially in some family.
We refer to \cite{Totaro} for the current state of the Tate conjecture.

\subsection{Raskind's $p$-adic Tate conjecture}
Suppose now that $F$ is a number field, choose a prime ideal $\mathfrak{p}\in\Spec\OO_{F}$ and let 
$F_{\idealp}$ be the $\idealp$-adic completion of $F$. 
Let $X$ be a smooth and proper variety over $F$. 
Then a standard argument (see  Proposition \ref{prop: relation to Tate}) shows that if $\eqref{tatehom}$ 
is surjective for the completion $X_{F_{\idealp}}$ at $\idealp$, then $X$ satisfies the Tate conjecture. 
It is therefore a natural idea to study the homomorphism $\eqref{tatehom}$ for varieties defined 
over $p$-adic fields (by which we shall mean finite extensions of $\QQ_{p}$). 
Such fields are \emph{not} finitely generated over their prime field.

In light of the previous paragraph, let $X$ be a smooth and proper variety over 
a $p$-adic field $K$. 
Then it is well-known that one cannot expect $\eqref{tatehom}$ to be surjective 
without imposing some further conditions on $X$ (see Appendix \ref{appendix} for counter-examples). 
Nevertheless, Raskind \cite{Raskind} has made a series of conjectures of Tate-type over such fields. 
In codimension one, that is, for divisors, they specialise to the following.    

\begin{Conjecture}[Raskind]\label{conj: raskind}
  Let $K$ be a $p$-adic field, let $\ell=p$, and let $X$ be a smooth and proper variety over $K$
  with totally degenerate reduction.
  Then, \eqref{tatehom} is surjective.
\end{Conjecture}

Of course, one has to specify what one means by
{\em totally degenerate reduction}:
roughly speaking, Raskind requires that $X$ has a strictly semi-stable model ${\cal X}\to\Spec\OO_K$, 
he asks that the Chow groups of all intersections of all components of the special fibre $\mathcal{X}_{0}$ 
to be as trivial as possible, and he requires $\mathcal{X}_{0}$ to be ordinary, 
see \cite[\S1]{Raskind}, \cite[Definition 1]{RX}, and Section \ref{subsec: Raskind}  for details.

In some sense, Conjecture \ref{conj: raskind} has a long history. 
We point out that the conjecture is true for varieties that are $p$-adically uniformisable by Drinfeld's 
upper half space (Proposition \ref{prop: drinfeld}).
This is an easy consequence of an observation of Rapoport and incorporated in work of 
Ito \cite[Appendix]{ItoWeight}, which relies on a result of Schneider and Stuhler \cite{SSp}. 
We also point out that it is a result of Tate that the conjecture is true when $X$ is the 
product of two Tate elliptic curves, as is explained in \cite[Appendix A.1.4]{Serre} 
(see also \cite[Corollary 19]{RXJacobians}). 
Moreover, a somewhat related conjecture has been formulated and established in several cases, 
such as abelian varieties,  by Tankeev \cite[\S 1,2 and 3]{Tankeev Degenerate} for varieties 
over function fields over $\CC$. 
Let us also remark that the restriction on $\ell$ being equal to $p$ and bad reduction 
is really necessary (see Appendix \ref{appendix}).

\subsection{Translation into a variational log-Tate conjecture}
Although Raskind's conjecture superficially looks to have a similar form as the Tate conjecture, we first 
show that it is in fact a variational conjecture. 
More precisely, using Fontaine's $\DD_{st}$-functor we identify the conjectural image of \eqref{tatehom} 
with
$$
 \Het{2}(X_{\overline{K}},\QQ_p(1))^{G_K} \,\cong\,\Hlogcris{2}({\cal X}_0/K_0)^{\varphi=p,N=0} \,\cap\, \Fil^1\HdR{2}(X/K),
$$
where $K_{0}$ is the maximal unramified extension of $\QQ_{p}$ contained in $K$, 
$\varphi$ and $N$ denote the Frobenius and monodromy operator
on log-crystalline cohomology, respectively, and $\Fil^\bullet$ denotes the Hodge filtration.
This translates Conjecture \ref{conj: raskind} into a ``variational
log-Tate conjecture'' as follows:
\begin{enumerate}
 \item By an appropriate log-version of Tate's conjecture for ${\cal X}_0$ over $k$, 
  one might expect 
 $\Hlogcris{2}({\cal X}_0/K_0)^{\varphi=p,N=0}$ to be equal to the $\QQ_p$-span
 of classes of invertible sheaves on ${\cal X}_0$, see Section  \ref{sec: translation}.
\item Since ${\cal X}_0$ is totally degenerate, there exists a combinatorial description
of $\Hlogcris{2}({\cal X}_0/K_0)$, $\varphi$ and $N$.
In fact, this cohomology group and its operators
arise naturally from a $\QQ$-vector space; a \emph{rational structure}
in the sense of Definition \ref{def: rational structure}.
In particular, one can be fairly explicit about the $\QQ$-span and the $\QQ_p$-span
of classes of invertible sheaves.
\item The intersection with the Hodge filtration $\Fil^1\HdR{2}(X)$ as a necessary
and sufficient condition to deform invertible sheaves from ${\cal X}_0$ to $X$ 
in the smooth case is a theorem of Berthelot and Ogus \cite[Theorem 3.8]{Berthelot Ogus}
(based on ideas of Deligne and Illusie, see \cite[p.~124 b)]{De81}), which has
been extended to the semi-stable situation by Yamashita \cite[Theorem 3.1]{Yamashita}.
\end{enumerate}

Now, the crucial point is that the just-mentioned theorems of Berthelot, Ogus, and Yamashita deal with the deformation of the {\em $\QQ$-span} of classes of invertible sheaves from the special
to the generic fibre, whereas Raskind's conjecture predicts this to be true even for the {\em $\QQ_p$-span} of classes of invertible sheaves when $\mathcal{X}_{0}$ is totally degenerate, see Remark \ref{rem: key observation}. This translates Conjecture \ref{conj: raskind} into a question about the 
interplay and the intersection of certain $\QQ$-vector spaces, certain $\QQ_p$-vector spaces, 
and the filtration step $\Fil^1$ of a filtered $(\varphi,N)$-module.
This leads to the notion of such a module being {\em Raskind-admissible} 
(Definition \ref{def: raskind admissible}) 
and we obtain the following reformulation
of Raskind's conjecture for divisors:

\begin{Theoremx}[Theorem \ref{thm: raskind translation}]
 Let $X$ be a smooth and proper variety over a $p$-adic field $K$ 
 with totally degenerate reduction.
 Then, the following are equivalent:
 \begin{enumerate}
 \item The homomorphism \eqref{tatehom} is surjective for $\ell=p$, that is,
   Conjecture \ref{conj: raskind} is true for $X$.
 \item The filtered $(\varphi,N)$-module $\DD_{{\rm st}}(\Het{2}(X_{\overline{K}},\QQ_p))$
   is Raskind-admissible with respect to the rational structure arising from
   ${\cal X}_0$.
 \end{enumerate}
\end{Theoremx}
One benefit of the reformulation is that it illuminates the known examples of varieties that 
satisfy Raskind's conjecture, and another is that Raskind-admissibility is a testable property in practice. 
For example, since (2) is a statement about filtered
$(\varphi,N)$-modules, it is tempting to approach
Conjecture \ref{conj: raskind} via this statement in semi-linear
algebra. In fact, in Section \ref{subsec: product tate} we use this approach to give a simple proof for the product of two Tate elliptic curves. 

\subsection{Abelian and abeloid varieties}\label{abelian and abeloid}
After seeing that Conjecture \ref{conj: raskind} is actually a variational conjecture about $\QQ_{p}$-classes 
of line bundles on totally degenerate varieties, we turn our attention to the study of 
Raskind's conjecture for abelian varieties. 
Let $A$ and $B$ be abelian varieties over a $p$-adic field $K$ and let $\ell$ be any prime number (possibly $\ell=p$). Functoriality gives a natural homomorphism of $\ZZ_\ell$-modules
\begin{equation*}\tag{$\star\star$}
\label{tatehom2}
   {\rm Hom}(A,B)\otimes_\ZZ\ZZ_\ell\,\to\,{\rm Hom}_{G_K}\left(T_\ell(A),T_\ell(B)\right),
\end{equation*}
where the subscript $G_K$ on the right indicates homomorphisms that are
$G_K$-equivariant.
This homomorphism is injective and its cokernel is torsion free. 
Then a classical K\"{u}nneth argument of Tate \cite{TateAbelian} relates 
Conjecture \ref{conj: raskind} to the following question. 

\begin{Question}\label{question: tate}
  Let $K$ be a $p$-adic field, let $\ell=p$, and let $A$ and $B$ be abelian varieties over $K$,
  both of which have totally degenerate reduction.
  Is it true that \eqref{tatehom2} is surjective?
\end{Question}

Raskind and Xarles have established Question \ref{question: tate} if both $A$ and $B$ are 
Tate elliptic curves \cite[Theorem 18]{RXJacobians}. 
Moreover, Mumford's results on $p$-adic uniformisation of abelian varieties  with totally degenerate 
reduction \cite{Mumford} make a positive answer plausible. 
Just as with Conjecture \ref{conj: raskind}, all of the assumptions are necessary. 
For example, if $\ell\neq p$ and if $A$ and $B$ are elliptic curves with good reduction, then 
Lubin and Tate \cite[\S3.5]{LT66} have given an example where surjectivity of \eqref{tatehom2} fails. 
For a comprehensive list of counter-examples where all combinations of assumptions are dropped, 
see Appendix \ref{appendix}. 
We also point out in passing that using the Kuga-Satake correspondence \cite{KS}, 
Conjecture \ref{conj: raskind} for abelian varieties implies a version of Conjecture \ref{conj: raskind} 
for projective hyperk\"ahler varieties (by adapting the arguments of Andr\'e \cite{Andre} or Tankeev \cite{Tankeev}).

There are several natural candidates for what it means for an abelian variety over a $p$-adic field $K$ to have 
totally degenerate reduction, but they turn out to be equivalent up to base change by a finite field extension, 
see Proposition \ref{prop: abelian totally degenerate}. 
It turns out that the point of view of admitting $p$-adic uniformisation in the sense of 
Mumford \cite{Mumford} is a convenient framework for our studies and computations. 
In particular, from Section \ref{subsec: abeloid varieties} onward, we will be working with lattices in 
$\GG_{m,K}^g$ and {\em abeloid varieties}, which are rigid analytic varieties over $K$ 
that are not necessarily algebraic schemes. 
We study Question \eqref{question: tate} in the enlarged context of abeloid varieties. 
Along the way prove and make use of the following results and computations, some of which may 
be of independent interest.
\begin{enumerate}
\item We describe the abelian groups
$$
  \Hom\left(A,B\right) \mbox{ \qquad and \qquad }\Hom_{G_K}\left(T_\ell(A),\,T_\ell(B)\right)
$$
for abeloid varieties over a $p$-adic field $K$ in terms of their lattices
(Theorem \ref{thm: abeloid morphisms}, Proposition \ref{prop: tate module}).
We note that the description of $\Hom(A,B)$ is essentially due to 
Gerritzen \cite{Gerritzen, GerritzenRiemann}, see also \cite{Ka07}.
\item We explicitly compute the filtered $(\varphi,N)$-module $\DD_{\rm st}(V_p(A))$
for an abeloid variety $A$ over a $p$-adic field 
in terms of a period matrix associated to a lattice (Theorem \ref{thm: abeloid phi N}),
where $V_p$ denotes the rational Tate module.
This generalises work of Berger \cite{Berger}, Coleman \cite{Co00}, 
Coleman--Iovita \cite{CI99}, and Le~Stum \cite{LeStum}.
As an application, we also describe 
$$
\Hom_{\mathrm{MF}_{K}^{\mathrm{wa},\varphi,N}}\left(\DD_{\mathrm{st}}(V_{p}(A)),\DD_{\mathrm{st}}(V_{p}(B))\right)
$$
in terms of lattices (Proposition \ref{prop: tate module2}).
\item We introduce an $\MazurL$-invariant for abeloid varieties that generalises the $\MazurL$-invariant
of a Tate elliptic curve. 
If the abeloid variety is the Jacobian $J$ of a Mumford curve $C$, then we show 
the Coleman-$\MazurL$-invariant of $C$ introduced by Besser and de~Shalit 
\cite{BdS16} coincides with our $\MazurL$-invariant for $J$ (Proposition \ref{prop: besser}). 
\end{enumerate}

\subsection{Counter-examples}
Crucially, just as with Raskind's conjecture above, these descriptions reformulate 
Question \ref{question: tate} into a question about the interplay of certain 
{\em $\QQ$-vector spaces} versus certain  {\em $\QQ_\ell$-vector spaces}, see 
Section \ref{subsec: reformulation question}. 
Using these computations we are able to give counter-examples to both 
Conjecture \ref{conj: raskind} and Question \ref{question: tate}.

\begin{Theoremx}[Theorem \ref{thm: counterexample}]
 Let $p$ be a prime with $p\geq5$ and $p\equiv 1\mod 3$.
 \begin{enumerate}
  \item There exists a Tate elliptic curve $A$ and an algebraisable abeloid surface (that is,
  an abelian surface with totally degenerate reduction) $B$  over $\QQ_p$
  such that \eqref{tatehom2} is {\em not} surjective for $\ell=p$.
  \item If $X=B$ with $B$ as in (1), then \eqref{tatehom} is {\em not} surjective for $\ell=p$.
 \end{enumerate}
\end{Theoremx}

On the other hand, our reformulation and following computations of Le~Stum \cite{LeStum} 
and Serre \cite{Serre} allows us to confirm Conjecture \ref{conj: raskind} and 
Question \ref{question: tate} for abeloid varieties 
which are isogenous to arbitrary products of Tate elliptic curves.

\begin{Propositionx}[Proposition \ref{prop: isogeny product tate}]
 Let $K$ be a $p$-adic field and let $A$ and $B$ be abelian varieties over $K$,
 both of which are isogenous to products of Tate elliptic curves.
 \begin{enumerate}
 \item Conjecture \ref{conj: raskind} is true for $A$, that is, \eqref{tatehom} is surjective
 for $\ell=p$.
 \item Question \ref{question: tate} is true for $A$ and $B$, that is, \eqref{tatehom2}
  is surjective for $\ell=p$.
 \end{enumerate}
\end{Propositionx}

Finally, let us remark that in this article, we study Conjecture \ref{conj: raskind} and 
Question \ref{question: tate} over $p$-adic fields. 
Of course, they can also be formulated and studied over local fields of equicharacteristic $p>0$. 
However, we expect that after replacing Yamashita's semi-stable Lefschetz theorem on $(1,1)$ 
classes \cite{Yamashita} with results of Morrow, Lazda, and P\'al \cite{Morrow, LazdaPal}, 
one should be able to set up everything in equicharacteristic $p>0$ and then, we expect that 
counter-examples similar to those of Theorem \ref{thm: counterexample} and 
Appendix \ref{appendix} should disprove them.

\subsection{Organisation}
The article is organised as follows:

In Section \ref{sec: generalities}, we establish general reduction
steps for Raskind's conjecture, such as the behaviour under field
extensions, dominant, or birational maps.
These are familiar from the analogous results for the classical Tate conjectures.
We also treat several simple cases and relate Raskind's conjecture to the
classical Tate conjectures over number fields.

In Section \ref{sec: translation} we translate Conjecture \ref{conj: raskind}
into semi-linear algebra and filtered $(\varphi,N)$-modules, we introduce the notion
of a rational structure, and we show that Conjecture \ref{conj: raskind}
is in fact equivalent to a problem in semi-linear algebra.

In Section \ref{sec: abelian varieties}, we first establish reduction
steps for Question \ref{question: tate} similar to those in Section \ref{sec: generalities}.
Then, we reformulate Raskind's notion of total degeneration for abelian varieties.
As a result, we focus on abeloid varieties: we describe their homomorphisms,
their $\ell$-adic Tate modules, and the filtered $(\varphi,N)$-modules associated
to their $p$-adic Tate modules.

In Section \ref{sec: product tate}, we do explicit computations with filtered
$(\varphi,N)$-modules arising from the product of two Tate curves.
This way, we prove Conjecture \ref{conj: raskind} for these varieties,
but we also produce admissible $(\varphi,N)$-modules 
that do not satisfy a more general version of Conjecture \ref{conj: raskind}.

In Section \ref{sec: counterexample}, we construct explicit examples that
disprove Conjecture \ref{conj: raskind}
and show that Question \ref{question: tate} has a negative answer.

In Appendix \ref{appendix}, we collect examples which show that 
Conjecture \ref{conj: raskind} and Question \ref{question: tate} also have 
a negative answer if one allows $\ell\neq p$ or if one does not consider
totally degenerate reduction.
Here, we claim only little originality.

\begin{Acknowledgements}
 We thank Valentina Di~Proietto, Thomas Geisser, Andreas Langer, Christopher Lazda,
 Bernard Le~Stum, Johannes Nicaise, Wies\l awa Nizio\l, Frans Oort, Otto Overkamp,
 Bjorn Poonen, Michael Rapoport, Wayne Raskind, and Burt Totaro for discussions and comments.
 We also thank the referee for comments and clarifications.
 Both authors are supported by the ERC Consolidator Grant 681838 ``K3CRYSTAL''.
\end{Acknowledgements}

\section*{Notations and Conventions}

Throughout the whole article, we fix the following notations
$$
\resizebox{1.0\hsize}{!}{$
\begin{array}{ll} 
  K & \mbox{ a $p$-adic field, that is, a finite field extension of $\QQ_p$} \\
  \OO_K & \mbox{ its ring of integers with maximal ideal } \idealm_K \\
  k & \mbox{ its residue field, that is, }\OO_K/\idealm_K\\
  \pi_K & \mbox{ a uniformiser of }\OO_K\\
  W(k) & \mbox{ the ring of Witt vectors, which we consider as subring of $\OO_K$}\\
  K_0 & \mbox{ the field of fractions of $W(k)$, which we consider as a subfield of $K$}\\
  \sigma & \mbox{ the Frobenius on $W(k)$ and $K_0$ }\\
  \overline{K}, \overline{k} & \mbox{ algebraic closures of $K$ and $k$, respectively}\\
  G_K, G_k & \mbox{ their absolute Galois groups} \\
   \nu_p & \mbox{ the extension of the standard valuation from $\QQ_p$ to $K$, that is, $\nu_p(p)=1$}\\
   \log_p & \mbox{ the Iwasawa logarithm, normalised such that $\log_p(p)=0$}
 \end{array}$}
$$
By a \emph{variety} over a field $F$, we mean a geometrically integral scheme of finite type
over $\Spec F$.
If $F'/F$ is a field extension and $X$ is a scheme over $F$, 
then we define $X_{F'}:=X\times_{\Spec F}\Spec F'$.

\section{Generalities}
\label{sec: generalities}

In this section, we recall some generalities concerning conjectures
of Tate-type for divisors.
These are well-known to the experts and we do not claim much, if any, 
originality.

\subsection{Setup}\label{subsec: Raskind setup}
Let $F$ be a field of characteristic $p\geq0$, 
let $F^{\rm sep}$ be a separable
closure, and let $G_F:={\rm Gal}(F^{\rm sep}/F)$ 
be its absolute Galois group.
If $X$ is a smooth and proper 
variety over $F$ and $\ell$ is a prime different from 
$p$, then the first Chern class map induces a $G_F$-equivariant
and injective homomorphism of finite-dimensional 
$\QQ_\ell$-vector spaces
$$
  c_{1,\ell}\,:\, {\rm NS}(X_{F^{\rm sep}})\otimes_\ZZ\QQ_\ell 
  \,\longrightarrow\, \Het{2}\left(X_{F^{\rm sep}},\QQ_\ell(1)\right).
$$
Taking $G_F$-invariants, we obtain an inclusion
of finite dimensional $\QQ_\ell$-vector spaces 
\eqref{tatehom}.
It is a natural question, whether this inclusion is in fact
an isomorphism, that is, whether it is surjective.

\subsection{Field extensions}\label{subsec: Raskind field extn}
Concerning this question,
we have the following remarks, which are well-known
in the context of the classical Tate conjecture
(see, for example \cite{Tate1, Tate2}), 
but perhaps not in this context:

\begin{Proposition}
  \label{fieldextension}
  Let $X$ be a smooth and proper variety over a field $F$,
  let $F\subset F'$ be a finite and separable field extension
  and let $\ell$ be a prime.
  If \eqref{tatehom} is surjective with respect to $X_{F'}$, $G_{F'}$,
  and $\ell$, then 
   \eqref{tatehom} is surjective with respect to $X$, $G_F$, and $\ell$.
\end{Proposition}

\prf
This is well-known, but we give a proof here since we are not aware of a reference. 
We consider $G_{F'}:={\rm Gal}(F'^{\rm sep}/F')={\rm Gal}(F^{\rm sep}/F')$
as subgroup of $G_F={\rm Gal}(F^{\rm sep}/F)$. 
Let $n$ be the degree of the extension $F'/F$. 
Suppose that \eqref{tatehom} is an isomorphism for $X_{F'}$. 
Let $\alpha\in\Het{2}(X_{F^{\rm sep}},\mathbb{Q}_{\ell}(1))^{G_{F}}$. 
Then, $\alpha$ is fixed by the open subgroup $G_{F'}=\mathrm{Gal}(F^{\rm sep}/F')\subset G_F$ 
and hence, there is a $z\in\NS(X_{F'})\otimes_{\ZZ}\QQ_{\ell}$ with $c_{1,\ell}(z)=\alpha$.
(Technical point: here, we are using that 
$(\NS(X_{F^{\rm sep}})\otimes_{\ZZ}\QQ_{\ell})^{G_{F'}}\iso\NS(X_{F'})\otimes_{\ZZ}\QQ_\ell$. 
This is not always true integrally since the Brauer group of $F$ may be non-trivial.
However, it is always true rationally since the Brauer group of a field is torsion.)
Let $f:X_{F'}\rightarrow X$ be the finite morphism given by the base extension. 
Then, $f_{\ast}(z)\in\NS(X)\otimes_{\ZZ}\QQ_{\ell}$ and $c_{1,\ell}(f_{\ast}(z))=nc_{1,\ell}(z)=n\alpha$. 
Since we are working with rational coefficients, we see that $\alpha$ is the class of a cycle on $X$. 
\qed\medskip

\subsection{Dominant and birational maps}
Next, we study the question whether surjectivity of \eqref{tatehom} is preserved under birational
maps and dominant maps. 
To do so, we adapt Tate's arguments from \cite{Tate2} to the $p$-adic case.

\begin{Proposition}\label{dombirat}
  Let $K$ be a field that is finitely generated over its prime subfield or a $p$-adic field.
   Let $X$ and $Y$ be smooth and proper varieties over $K$ and assume that
   \eqref{tatehom} is surjective for $X$. 
   \begin{enumerate}
    \item If there exists a dominant and rational map $X\dashrightarrow Y$
     of varieties over $K$, or
    \item if $X$ and $Y$ are birationally equivalent varieties over $K$,
   \end{enumerate}
   then \eqref{tatehom} is surjective for $Y$.
\end{Proposition}

\prf
This is \cite[Theorem 5.2(b)]{Tate2} in the case that $K$ is finitely generated over its prime subfield. 
Tate shows that the \eqref{tatehom} is an isomorphism for $X$ if and only if \eqref{tatehom} 
is an isomorphism for an arbitrary dense open $U\subset X$ 
by using the Gysin sequence for $U\hookrightarrow X$. 
A weight argument then reduces showing that ``the Tate conjecture for divisors on $X$ is 
equivalent to the Tate conjecture for divisors on $U$'' to showing that numerical equivalence 
coincides with $\ell$-adic homological equivalence for divisors on $X_{\overline{K}}$, 
where $\overline{K}$ is an algebraic closure of $K$. 
The coincidence of numerical equivalence and homological equivalence 
(defined using any Weil cohomology theory) is known over algebraically closed fields, see
\cite{Ma57} or \cite[Proposition 3.4.6.1]{An04}. 
To prove the proposition for $p$-adic fields, the same proof works when $K$ since there is an 
appropriate theory of weights (see \cite{Ja10} for a summary of both cases $\ell\neq p$ and $\ell=p$).
\qed\medskip

\subsection{A simple case}
As an easy consequence of the Lefschetz theorem on $(1,1)$-classes
and the Lefschetz principle, we obtain the following corollary
of Proposition \ref{fieldextension}.

\begin{Proposition}
 \label{prop: simple case}
 Let $F$ be a field of characteristic zero.
 Let $X$ be a smooth and proper variety over $K$
 with $H^2(X,\OO_X)=0$.
 Then, \eqref{tatehom} is surjective for all primes $\ell$.
\end{Proposition}

\prf
Being of finite type over $K$, there exists a subfield $F'\subseteq F$
that is finitely generated over $\QQ$ such that $X$ can be defined
over $F'$.
Being finitely generated over $\QQ$, we may choose an embedding
$F'\hookrightarrow\CC$.
Let $X_\CC$ be the base change of a model of $X$ over $F'$
to $\CC$.
Since $X_\CC$ also satisfies $H^2(X_\CC,\OO_{X_\CC})=0$,
the Lefschetz theorem on $(1,1)$-cycles shows that the rank
of $\NS(X_\CC)$ is equal to the second Betti number
$b_2(X_\CC)$.
Thus, already the rank of $\NS(X_{\overline{F'}})$ is
equal to the $\QQ_\ell$-dimension of $\Het{2}(X_{\overline{F'}},\QQ_\ell)$
for some algebraic closure $\overline{F'}$ of $F'$ inside $\CC$.
Since the N\'eron-Severi group is finitely generated, 
there exists a finite field extension $F'\subseteq F''$,
such that the rank of $\NS(X_{F''})$ is equal 
to the $\QQ_\ell$-dimension of $\Het{2}(X_{F''},\QQ_\ell)$.
Thus, the $G_{F''}$-actions on $\NS(X_{\overline{F'}})$ and 
$\Het{2}(X_{F''},\QQ_\ell(1))$ are trivial and
\eqref{tatehom} is an isomorphism for $X_{F''}$.
Thus, \eqref{tatehom} is surjective for $X$
by Proposition \ref{fieldextension}.
\qed\medskip

\begin{Remark}
 \label{remark}
 For example, this includes varieties that 
 are birationally equivalent to smooth and proper varieties 
 over $\overline{F}$ that are
  rationally connected (these satisfy $H^2(X,\OO_X)=0$),
  which includes rational and unirational varieties, and
  Fano varieties.
  It also includes geometrically ruled surfaces and
  Calabi-Yau varieties of dimension at least three
  (even in the liberal sense of varieties whose canonical divisor class is numerically trivial and
  that satisfy $H^i(X,\OO_X)=0$ for $0<i<\dim(X)$).
    
 In particular, the conjectures of Tate and Raskind
 for divisor holds for these classes of varieties even without 
 extra assumptions on finite generation of the field over its prime subfield 
 or on total degeneration.
\end{Remark}

\subsection{Drinfeld's upper half space}
A second simple case is the following:
let $X$ be a smooth and proper variety over a number field $F$.
Assume that there exists a finite extension $F\subset F'$ and a finite 
place $w$ of $F'$ such that  $X\times_FF'_w$ is isomorphic to 
$\Gamma\backslash\widehat{\Omega}_{F'_w}^d$,
where $F'_w$ denotes the $w$-adic completion,
where $\widehat{\Omega}^d_{F'_w}$ denotes the Drinfeld upper
half space of dimension $d\geq1$, 
and where $\Gamma\subset{\rm PGL}_{d+1}(F'_w)$ is a cocompact and torsion
free discrete subgroup.
In \cite[Theorem 7.1]{ItoWeight}, Ito established the Tate conjecture for such varieties,
which is based on ideas of Rapoport.
Essentially the same proof shows the following observation, see also 
\cite[Remark 1.3]{ItoWeight}.

\begin{Proposition}\label{prop: drinfeld}
 Let $K$ be a $p$-adic field, let  $\Gamma\subset{\rm PGL}_{d+1}(K)$ be 
 a cocompact and torsion free discrete subgroup, and set
 $X_\Gamma:=\Gamma\backslash\widehat{\Omega}_K^d$,
 Then, \eqref{tatehom} is surjective for $X_\Gamma$ and all primes $\ell$.
\end{Proposition}

\prf
If $\ell\neq p$, then $\Het{2}(X_{\Gamma, \overline{K}},\QQ_\ell)$ is
one-dimensional by \cite[Theorem 4]{SSp}.
After choosing an embedding $K\into \CC$ and using comparison theorems
with singular cohomology, it follows that also $\Het{2}(X_{\Gamma,\overline{K}},\QQ_p)$ is one-dimensional.
It follows from work of Kurihara and Mustafin \cite{Kurihara, Mustafin} 
(see also the discussion in \cite[\S 6]{ItoWeight})
that $X_\Gamma$ is projective and thus,  $\NS(X_\Gamma)\otimes\QQ$ is non-zero.
Thus, \eqref{tatehom} is surjective for all primes $\ell$ for 
dimensional reasons.
\qed\medskip

\begin{Remark}
 Such an $X_\Gamma$ admits a semi-stable model over $\OO_K$, whose
special fibre is totally degenerate in the sense of Raskind, see 
\cite[Remark 1.3]{ItoWeight} and \cite[Example 1.(iii)]{RX}.
\end{Remark}

\subsection{The Tate conjecture over number fields}\label{sec: relation to Tate}
Let us also relate the conjecture of Raskind over $p$-adic fields
to the conjecture of Tate over number fields, which was
already observed by Raskind \cite{Raskind}.
For a number field $F$, we let $\OO_F$ be its ring of integers
and for a prime ideal $\idealp\in\Spec\OO_F$, we denote by
$F_\idealp$ the $\idealp$-adic completion of $F$. 
The following is a slight generalisation of \cite[Proposition 1]{Raskind}.

\begin{Proposition}\label{prop: relation to Tate}
  Let $X$ be a smooth and proper variety over a number field $F$
  and assume that there exists a  prime $\ell$ and
  a prime ideal $\idealp\subset\OO_F$ such that
  \eqref{tatehom} is surjective for $X_{F_\idealp}$
  $G_{F_\idealp}$, and $\ell$.
  Then, \eqref{tatehom} is surjective for $X$,
  $G_F$, and $\ell$.
\end{Proposition}

\prf
Let $\alpha\in\Het{2}(X_{\overline{F}},\QQ_{\ell}(1))^{G_{F}}$. 
Then because $\Het{2}(X_{\overline{F}},\QQ_{\ell}(1))\cong\Het{2}(X_{\overline{F_{\idealp}}},\QQ_{\ell}(1))$ 
and $G_{F_{\idealp}}\subset G_{F}$, we see that 
$\alpha\in\Het{2}(X_{\overline{F_{\idealp}}},\QQ_{\ell}(1))^{G_{F_{\idealp}}}$. 
Therefore there is a $z\in\mathrm{NS}(X_{F_{\idealp}})\otimes_{\ZZ}\QQ_{\ell}$ 
with $c_{1,\ell}(z)=\alpha$ by assumption. 
Since $\mathrm{NS}(X_{\overline{F}})=\mathrm{NS}(X_{\overline{F_{\idealp}}})$ 
(the N\'{e}ron-Severi group is invariant under algebraically closed base extension), 
we see that $z\in\mathrm{NS}(X_{\overline{F}})\otimes_{\ZZ}\QQ_{\ell}$, and in particular 
$z\in\mathrm{NS}(X_{F'})\otimes_{\ZZ}\QQ_{\ell}$ for some finite extension $F'/F$. 
By making a finite extension if necessary, we may assume that $F'/F$ is Galois. 
Summing over the $\mathrm{Gal}(F'/F)$-conjugates of $z$, we obtain a $G_F$-invariant class 
$z'$, such that $c_{1,\ell}(z')$ is a non-zero multiple of $c_{1,\ell}(z)$. 
Thus, $c_{1,\ell}(z)$ lies in the image of \eqref{tatehom}.
\qed\medskip

\section{A translation into semi-linear algebra}
\label{sec: translation}

In this section, we use Fontaine's functor $\DD_{\rm st}$ and Yamashita's
$p$-adic semi-stable Lefschetz $(1,1)$-theorem to translate
Raskind's conjecture for divisors (Conjecture \ref{conj: raskind})
into a question about semi-linear algebra and filtered $(\varphi,N)$-modules.

More precisely, we define the notion of a {\em rational structure} on a $(\varphi,N)$-module
and show how the special fibre of a model ${\cal X}\to\Spec\OO_K$
of some smooth and proper variety $X$ over $K$ with total degeneration gives
rise to such a structure.
Finally, we introduce the notion of {\em Raskind-admissibility}, which 
is the semi-linear algebra version of the Raskind conjecture for divisors
on the level of filtered $(\varphi,N)$-modules with rational structure.

\subsection{Translation into filtered $(\varphi,N)$-modules}
\label{logcrys}
Let $X$ be a smooth and proper variety over a $p$-adic field $K$ that admits
a proper and semi-stable model 
$$
  \pi\,:\,{\cal X}\,\to\,\Spec\OO_K
$$
over the ring of integers $\OO_K$ of $K$, that is, 
${\cal X}$ is a regular scheme, $\pi$ is a proper and flat morphism,
the generic fibre of $\pi$ is isomorphic to $X$, and the special
fibre ${\cal X}_0$ is a semi-stable scheme over the residue field $k$ of $\OO_K$.
Here, semi-stable means that ${\cal X}_0$ is a strict normal crossing divisor.
In particular, the components of ${\cal X}_0$ are smooth and geometrically integral
over $k$.
Let $W(k)$ be the ring of Witt vectors, which we consider as subring of $\OO_K$, and let
$K_0$ be the field of fractions of $W(k)$.
Then, $K_0$ is the maximal unramified extension of $\QQ_p$ inside $K$.
Endow $\mathcal{X}$ with the log structure induced by $\mathcal{X}_{0}$, and 
let $M$ denote the pullback of this log structure on $\mathcal{X}_{0}$. 
Then, $(\mathcal{X}_{0},M)$ is a fine and log-smooth log scheme over 
$(\Spec k, \NN\rightarrow k, 1\mapsto 0)$, see \cite[2.13.2]{Hyodo Kato}. 
We shall write
$$
 \Hlogcris{n}(\mathcal{X}_{0}/K_{0}) \,:=\,
 H^{n}(((\mathcal{X}_{0},M)/(W(k),\NN))_{\rm cris},\OO_{(\mathcal{X}_{0},M)/(W(k),\NN)})\otimes_{W(k)}K_{0},
$$
where $H^{n}(((\mathcal{X}_{0},M)/(W(k),\NN))_{\rm cris},\OO_{(\mathcal{X}_{0},M)/(W(k),\NN)})$ 
is the log-crystalline cohomology of $(\mathcal{X}_{0},M)\rightarrow(\Spec k,\NN)$. 
Then $\Hlogcris{n}(\mathcal{X}_{0}/K_{0})$ is equipped with a semi-linear endomorphism $\varphi$ (Frobenius)
and a linear endomorphism $N$ (monodromy) satisfying the relation $N\varphi=p\varphi N$, 
making the triple 
$$
\left(\Hlogcris{n}(\mathcal{X}_{0}/K_{0}),\varphi, N\right)
$$
a $(\varphi,N)$-module. 
We refer to \cite[\S3]{Hyodo Kato} for the details.  

Since $X$ has semi-stable reduction, the $G_K$-representation 
on $\Het{n}(X_{\overline{K}},\QQ_p)$ is semi-stable in the sense of Fontaine for every $n$ \cite[Theorem 0.2]{Ts99}.
Next, let  $B_{\rm st}$ be Fontaine's period ring and if $V$ is a finite dimensional $\QQ_p$-vector space 
with a continuous $G_K$-action,
that is, a $p$-adic Galois-representation, then we have a filtered $(\varphi,N)$-module
over $K$
$$
   \DD_{\rm st}(V) \,:=\, \left( V\otimes_{\QQ_p} B_{\rm st} \right)^{G_K},
$$
that is, a $K_0$-vector space with a semi-linear operator $\varphi$, a linear operator $N$,
and a filtration $\Fil^\bullet$ on this vector space tensored with $K$.
We recall that $V$ is said to be \emph{semi-stable} if the inequality
$\dim_{K_0}   \DD_{\rm st}(V)\leq\dim_{\QQ_p}(V)$ is an equality. 
We refer to \cite{CF00} for details. Fontaine's functor $\DD_{\rm st}$ establishes an equivalence 
of categories between the category of semi-stable $G_K$-representations and the 
category $MF_{K}^{\text{wa}, \varphi, N}$ of admissible filtered $(\varphi,N)$-modules 
over $K$ \cite[Th\'{e}orem\`{e} A]{CF00} .

By the semi-stable comparison theorem \cite[Theorem 0.2]{Ts99}, the admissible filtered 
$(\varphi,N)$-module $\DD_{\rm st}(\Het{n}(X_{\overline{K}},\QQ_p))$ is equal to 
$$
   D^n \,:=\, \left( \Hlogcris{n}({\cal X}_0/K_0),\, \Fil^{\bullet}\HdR{n}(X/K),\,\varphi,\,N \right).
$$
Using this translation, we have the following.

\begin{Proposition}\label{prop: logcrys}
 Let $X$ be a smooth and proper variety over $K$ and assume that there exists
 a proper and semi-stable model ${\cal X}\to\Spec\OO_K$ of $X$.
 Let ${\cal X}_0$ be the special fibre.
 Then, there exists an isomorphism of $\QQ_p$-vector spaces
 \begin{equation}
  \label{eq: logcrys}
  \Het{n}\left(X_{\overline{K}},\QQ_p(m)\right)^{G_K} \,\iso\, 
    \Hlogcris{n}({\cal X}_0/K_0)^{\varphi=p^m, N=0}\,\cap\, \Fil^m \HdR{n}(X/K)
 \end{equation}
 for all non-negative integers $m,n$.
\end{Proposition}

\prf 
This follows from the equalities and isomorphisms
\begin{eqnarray*}
   \Het{n}\left(X_{\overline{K}},\QQ_p(m)\right)^{G_K}  &=& \Hom_{G_K}\left( \QQ_p,\, \Het{n}(X_{\overline{K}},\QQ_p(m))\right) \\
   &\cong& \Hom_{\MF_K^{\text{wa},\varphi,N}}\left(K,\,D^n(m)\right) \\
   &\cong& \left\{ x\in \Hlogcris{n}({\cal X}_0/K_0)\,:\, \varphi(x)=p^m\cdot x,\, N(x)=0 \right\} \\
   & & \mbox{\quad} \,\cap\, \Fil^m \HdR{n}(X/K),
\end{eqnarray*}
where $K$ denotes the trivial filtered $(\varphi,N)$-module.
\qed\medskip

\begin{Remark}\label{remark on l vs p}
 If $\ell\neq p$, then we have
 $$
  \dim_{\QQ_\ell} \Het{n}\left(X_{\overline{K}},\QQ_\ell(m)\right)^{G_K} \,=\, 
  \dim_{\QQ_p} \Hlogcris{n}({\cal X}_0/K_0)^{\varphi=p^{m},N=0}
 $$
 by \cite{Katz Messing}.
 This explains why the dimensions of $\Het{n}(X_{\overline{K}},\QQ_\ell(m))^{G_K}$ behave differently 
 for $\ell=p$ and $\ell\neq p$ and it also shows that for $\ell\neq p$, these vector spaces
 capture information about the special fibre ${\cal X}_0$ only.
 We refer to Consani's article \cite{Consani} for background and some conjectures.
\end{Remark}

We end our discussion by presenting some probably well-known dimension estimates:
since $K$ is of characteristic zero, the Fr\"olicher spectral sequence 
\begin{equation*}
E_{1}^{r,s}=H^{s}(X,\Omega_{X/K}^{r})\Rightarrow H_{\rm dR}^{r+s}(X/K)
\end{equation*}
degenerates at $E_1$. 
In particular, by Proposition \ref{prop: logcrys} and Remark \ref{remark on l vs p} we obtain the dimension estimates
\begin{align*}
\dim_{\QQ_\ell}\Het{n}(X_{\overline{K}},\QQ_\ell(m))^{G_K}-\sum_{i=m+1}^{n}h^{i,n-i}(X)\,
& \leq\,\dim_{\QQ_p}\Het{n}(X_{\overline{K}},\QQ_p(m))^{G_K} \\
& \leq\,\dim_{\QQ_\ell}\Het{n} (X_{\overline{K}},\QQ_\ell(m))^{G_K}\,,
\end{align*}
where $h^{i,j}(X):=\dim_{K}H^{j}(X,\Omega_{X/K}^{i})$. 
In the case of interest to us, that is to say when $n=2$ and $m=1$, this gives
\begin{align*}
\dim_{\QQ_\ell}\Het{2}(X_{\overline{K}},\QQ_\ell(1))^{G_K}-h^{2,0}(X)\,
& \leq\,\dim_{\QQ_p}\Het{2}(X_{\overline{K}},\QQ_p(1))^{G_K} \\
& \leq\,\dim_{\QQ_\ell}\Het{2} (X_{\overline{K}},\QQ_\ell(1))^{G_K}.
\end{align*}

\subsection{Rational structures}
Next, we deal with the log-crystalline cohomology of the special fibre ${\cal X}_0$
of a proper and semi-stable model ${\cal X}\to\Spec\OO_K$ of $X$. Let
$$
  Y\,=\, \bigcup_{i\in I} Y_i
$$
be the decomposition of the special fibre $Y=\mathcal{X}_{0}$ into irreducible components. 
For a subset $J\subset I$, we denote by $Y_J$ the intersection of all $Y_j$ with $j\in J$.
Since $Y$ is strict normal crossing, each $Y_J$ is a smooth, proper, and geometrically integral 
scheme over $k$.
Moreover, we denote by $Y^{[i]}$ the disjoint union of all $Y_J$'s where $J$ 
has $(i+1)$ elements, that is, the subscript $i$ is equal to the codimension of $Y_J$ in $Y$.
By Mokrane \cite[\S3]{Mokrane} and Nakkajima \cite[\S4]{Nakkajima RZ},
there exists a $p$-adic Steenbrink-Rapoport-Zink spectral sequence
\begin{equation}
 \label{eq: SRZ spectral sequence}
   E_1^{-k,h+k}\,=\,\bigoplus_{j\geq\max\{-k,0\}}\Hcris{h-2j-k}(Y^{[2j+k]}/K_0)(-j-k) \,\Rightarrow\, \Hlogcris{h}(Y/K_0)\,,
\end{equation}
which degenerates at $E_2$.
This spectral sequence is compatible with the $F$-isocrystal structures on both sides
and induces a monodromy operator $N$ on the right hand side.

\begin{Definition}\label{def: cohomologically totally degenerate}
 A scheme $Y$ over a finite field $k$ is called \emph{cohomologically totally degenerate} if it is strictly normal crossing, say, equal to $\bigcup_{i\in I} Y_i$, where the $Y_i$ are the irreducible components of $Y$, 
 such that
 \begin{enumerate}
 \item for all $j$ and all odd integers $i$, the crystalline cohomology groups $\Hcris{i}(Y^{[j]}/K_0)$ are zero and
 \item for all $i$ and $j$, the cycle class maps
 $$
  {\rm CH}^{j}(Y^{[i]})\otimes K_0 \,\to\, \Hcris{2j}(Y^{[i]}/K_0)(j)
 $$
 are isomorphisms.
 \end{enumerate}
 \end{Definition}
 
 Put differently, all crystalline cohomology groups of all intersections of components of $Y$ are spanned by classes
 of algebraic cycles.
 In particular, since the Chow groups of a variety are $\QQ$-vector spaces and since the action of Frobenius
 on Chow groups is trivial up to Tate twist, this implies that the log-crystalline cohomology of ${\cal X}_0$
 is of a very simple form.
 For $\Hlogcris{2}({\cal X}_0/K_0)$, it leads to the following.
 
 \begin{Definition}\label{def: rational structure}
  Let $k$ be a finite field, let $K_0={\rm Frac}(W(k))$, let $\sigma$ be the Frobenius on $K_0$,
  and let $H$ be a $(\varphi,N)$-module over $K_0$.
   A {\em rational structure} on a  $(\varphi,N)$-module $H$ consists of a 
  finite-dimensional $\QQ$-vector $V$ space together with a direct sum decomposition
  $$
         V \,=\, A \,\oplus\, B_0\,\oplus\,B_1\,\oplus\, C
  $$
  and two $\QQ$-linear endomorphisms $\varphi_V$ and $N_V$ such that
  \begin{enumerate}
   \item $N_V$ is zero on $B_1$ and $A$, and $N_V$ induces isomorphisms
   $$
      C\,\stackrel{N}{\longrightarrow}\, N(C)=B_0\mbox{ \quad and \quad } 
      B_0\,\stackrel{N}{\longrightarrow} N(B_0)=A\,.
   $$ 
   \item $\varphi_V$ acts as identity on $A$, as multiplication by $p$ on $B_0\oplus B_1$,
    and as multiplication by $p^2$ on $C$.
   \item As $(\varphi,N)$-module, $H$ is isomorphic to $V\otimes_\QQ K_0$
    with $\varphi=\varphi_V\otimes\sigma$ and $N=N_V\otimes{\rm id}$.
  \end{enumerate}
 \end{Definition} 
 
If $V$ is a rational structure on a filtered $(\varphi,N)$-module $H$ over $K_0$,
then we have an isomorphism of $\QQ_p$-vector spaces
\begin{equation}
 \label{eq: phi and N}
   H^{\varphi=p, N=0} \,\cong\, B_1\otimes_\QQ\QQ_p.
\end{equation}
For example, if $m=1$ and $n=2$ in Proposition \ref{prop: logcrys} and 
if the $(\varphi,N)$-module $\Hlogcris{2}({\cal X}_0/K_0)$ there
comes with a rational structure, then \eqref{eq: phi and N} 
makes the right hand side of \eqref{eq: logcrys} much easier to compute.
Before establishing a natural rational structure on 
$\Hlogcris{2}({\cal X}_0/K_0)$, we need one more definition.

\begin{Definition}
 The \emph{dual graph of $Y$} is the simplicial complex $\Gamma$ that has
 one vertex $P_i$ for each component $Y_i$ of $Y$ and the simplex
 $\langle P_{i(0)},...,P_{i(k)}\rangle$ belongs to $\Gamma$ if and only if 
 $Y_J$ for $J=\{i(0),...,i(k)\}$ is non-empty.
 We define $H^*(\Gamma):=H_{\rm sing}^*(|\Gamma|,\QQ)$ to be 
 the singular cohomology
 of the topological realisation $|\Gamma|$ of $\Gamma$.
\end{Definition}

 \begin{Proposition}\label{prop: rational structure}
    Let $Y=\bigcup_{i\in I} Y_i$ be a cohomologically totally degenerate scheme
    over a finite field $k$
    and consider the $(\varphi,N$)-module $\Hlogcris{2}(Y/K_0)$.
    Then, 
    \begin{enumerate}
    \item The cycle class maps 
    $$
        {\rm CH}^{j}(Y^{[i]})\otimes\QQ \,\to\,  \Hcris{2j}(Y^{[i]}/K_0)(j)
     $$
     followed by the $p$-adic spectral sequence \eqref{eq: SRZ spectral sequence}
     induce a natural rational structure 
     $$
        \left( V=A\oplus B_0\oplus B_1\oplus C, \varphi_V, N_V \right)
     $$
     on $\Hlogcris{2}(Y/K_0)$.
     \item If $\Gamma$ denotes the dual complex of $Y$ and if
        $H^2(\Gamma)\neq0$, then $N$ has maximally nilpotent monodromy, that is, $N^2\neq0$.
     \item      If $Y^{\rm sm}$ denotes the smooth locus of $Y$ and if $Y$ 
     is equipped with its natural log-structure, then there exist homomorphisms
     $$
       \Pic(Y) \,\to\, \Pic(Y^{\rm sm}) \,\cong\, \Pic^{\log}(Y),
     $$
     where the first map is restriction and the second is an isomorphism.
     
     The $p$-adic spectral sequence \eqref{eq: SRZ spectral sequence}
      gives rise to an isomorphism
     $$
       \Pic(Y)\otimes\QQ \,\cong\, B_1.
     $$
     
     Moreover, the first Chern class maps give rise to a commutative diagram
     \begin{center}
\begin{tikzpicture}[descr/.style={fill=white,inner sep=1.5pt}]
        \matrix (m) [
            matrix of math nodes,
            row sep=2.5em,
            column sep=2.5em,
            text height=1.5ex, text depth=0.25ex
        ]
        {  \Pic(Y)\otimes F & \ \\
        \Pic^{\rm log}(Y)\otimes F & \Hlogcris{2}(Y/K_0)^{\varphi=p,N=0}, \\
        };

       \path[overlay,->, font=\scriptsize]
       (m-1-1) edge (m-2-1)
       (m-1-1) edge (m-2-2)
       (m-2-1) edge (m-2-2);
          
        \end{tikzpicture} 
\end{center}
where the images of both Chern class maps are equal 
     to $B_1$ if $F=\QQ$ and equal to $\Hlogcris{2}(Y/K_0)^{\varphi=p,N=0}$ 
     if $F=\QQ_p$. 
    \end{enumerate}
 \end{Proposition}
 
 \prf
 Using the cycle class maps, the assumption on cohomological degeneracy,
 and the spectral sequence \eqref{eq: SRZ spectral sequence},
 we obtain a $\QQ$-vector space $V$, such that $V\otimes K_0$ is naturally isomorphic 
 to $\Hlogcris{2}(Y/K_0)$.
 Moreover, since Frobenius acts on cohomology classes of cycles by multiplication
 by some power of $p$, we obtain a direct
 sum decomposition $V=V_0\oplus V_1\oplus V_2$ of $\QQ$-vector spaces
 together with a linear operator $\varphi_V$ that is multiplication by $p^i$ on $V_i$,
 such that the $F$-isocrystal structure on $\Hlogcris{2}(Y/K_0)$ is isomorphic to
$(V\otimes K_0,\varphi_V\otimes\sigma)$.
Also, the monodromy operator $N$ and the weight filtration arise from 
the spectral sequence and using the $\QQ$-vector space structures,
we obtain a monodromy operator on $V$.
For details, we refer to the discussions in \cite[\S 1]{BGS} or
\cite{Morrison}.

The fact that $N^2\neq0$ is equivalent to $H^2(\Gamma)\neq0$ is shown
in \cite[\S6]{Morrison} and although they are stated in the framework of
complex geometry in loc.cit., the arguments carry over literally to our situation.
This establishes claims (1) and (2).

Concerning claim (3):
first, we have a restriction homomorphism 
$\Pic(Y)\to\Pic(Y^{\rm sm})$ and an isomorphism 
$\Pic(Y^{\rm sm})\cong\Pic^{\log}(Y)$ by the 
discussion at the beginning of \cite[\S 2]{Yamashita}.

Next, let $\delta_{i}:Y^{[j]}\rightarrow Y^{[j-1]}$ be the morphisms induced by the obvious inclusions
$$
 Y_{\iota_{1}}\cap\ldots\cap Y_{\iota_{j+1}} \,\hookrightarrow\,  
 Y_{\iota_{1}}\cap\ldots\cap Y_{\iota_{i-1}}\cap Y_{\iota_{i+1}}\cap\ldots\cap Y_{\iota_{j+1}}\,.
$$
To give an invertible sheaf on $Y$ is equivalent to giving an invertible
sheaf on $Y^{[0]}$ plus compatibilities under the restriction maps 
$\delta_{1}^{\ast},\delta_{2}^{\ast}:Y^{[0]}\rightrightarrows Y^{[1]}$. 
In particular, we obtain an isomorphism
$$
\Pic(Y)\otimes\QQ\,\cong\,
\ker\left(\Pic(Y^{[0]})\otimes\QQ\xrightarrow{\delta_{2}^{\ast}-\delta_{1}^{\ast}}\Pic(Y^{[1]})\otimes\QQ\right)\,.
$$

Now, the boundary morphisms $\delta_{i}$ give an augmented simplicial scheme 
$Y^{[\bullet]}\rightarrow Y$, which is a proper smooth hypercovering of $Y$. 
Since each $Y^{[j]}$ is smooth, we have $\Hcris{\ast}(Y^{[j]}/K_{0})\cong \Hrig{\ast}(Y^{[j]}/K_{0})$ 
by \cite[Proposition 1.9]{Be97} and hence
$$
\Hcris{\ast}(Y^{[\bullet]}/K_{0}) \,\cong\, \Hrig{\ast}(Y^{[\bullet]}/K_{0})\cong H_{\mathrm{rig}}^{\ast}(Y/K_{0}),
$$
where the second isomorphism is because rigid cohomology satisfies cohomological descent 
for proper hypercoverings \cite[Corollary 2.2.3]{Ts03}. 
Now, consider the spectral sequence
$$
E_{1}^{s,t} \,=\, \Hrig{t}(Y^{[s]}/K_{0}) \,\Rightarrow\, \Hrig{s+t}(Y/K_{0})
$$
 of the (hyper)covering (see, for example \cite[Theorem 4.5.1]{Ts03}). 
 This degenerates at $E_{2}$ by a standard weight argument (our assumption that $Y$ is cohomologically 
 totally degenerate makes this argument very easy, but see also \cite[Corollary 5.2.4]{Ts03} for the 
 general statement). 
 In particular, we find isomorphisms
 \begin{align*}
\Hrig{2}(Y/K_{0})
 & \cong\ker\left(\Hrig{2}(Y^{[0]}/K_{0})\xrightarrow{\delta_{2}^{\ast}-\delta_{1}^{\ast}} \Hrig{2}(Y^{[1]}/K_{0})\right) \\
 & \cong\ker\left(\Hcris{2}(Y^{[0]}/K_{0})\xrightarrow{\delta_{2}^{\ast}-\delta_{1}^{\ast}} \Hcris{2}(Y^{[1]}/K_{0})\right) \\
& \cong\ker\left(\Pic(Y^{[0]})\otimes\mathbb{Q}_{p}\xrightarrow{\delta_{2}^{\ast}-\delta_{1}^{\ast}}\Pic(Y^{[1]})\otimes\mathbb{Q}_{p}\right) \\
& \cong\Pic(Y)\otimes\mathbb{Q}_{p},
 \end{align*}
 where the third isomorphism is induced by the crystalline Chern class map, and is an isomorphism 
 by the assumption that $Y$ is cohomologically totally degenerate. 
 The compatibility between rigid and crystalline Chern classes \cite[Th\'{e}or\`{e}me 5.2.3]{Pe03} implies 
 that the above isomorphism $\Pic(Y)\otimes\QQ_{p}\xrightarrow{\sim}H_{\mathrm{rig}}^{2}(Y/K_{0})$ 
 is induced by the rigid Chern class map $c_{1}^{\mathrm{rig}}:\Pic(Y)\rightarrow H_{\mathrm{rig}}^{2}(Y/K_{0})$.  
 
 Finally, one has a square
      \begin{center}
 \begin{tikzpicture}[descr/.style={fill=white,inner sep=1.5pt}]
        \matrix (m) [
            matrix of math nodes,
            row sep=2.5em,
            column sep=2.5em,
            text height=1.5ex, text depth=0.25ex
        ]
        {  \Pic(Y) & \Pic^{\log}(Y) \ \\
        H_{\mathrm{rig}}^{2}(Y/K_{0}) & \Hlogcris{2}(Y/K_0)^{\varphi=p,N=0}, \\
        };

       \path[overlay,->, font=\scriptsize]
       (m-1-1) edge (m-1-2)
       (m-1-1) edge node [right]{$c_{1}^{\mathrm{rig}}$} (m-2-1)
       (m-1-2) edge node [right]{$c_{1}$} (m-2-2);
       
       \path[overlay,->>, font=\scriptsize]
       (m-2-1) edge (m-2-2);
          
        \end{tikzpicture} 
\end{center}
where the top arrow is induced by the inclusion 
$\OO_{Y}^{\times}\hookrightarrow M^{\mathrm{gp}}$ (recall that $M$ is the log structure on $Y$). 
The surjection is the one given by the Clemens-Schmid exact sequence \cite{CT14}. 
Recall from loc. cit. that this map is the composition
$$
\begin{array}{lclcl}
 \Hrig{2}(Y/K_{0}) &\cong& H^{2}_{\mathrm{conv}}(Y/W(k)) 
  &\rightarrow& H_{\mathrm{log-conv}}^{2}((Y,M)/(W(k),\NN)) \\
  &&&\cong&\Hlogcris{2}(Y/K_0),
\end{array}
$$
where the first isomorphism is because $Y$ is proper and the second isomorphism 
is because $(Y,M)\rightarrow (\Spec k,\NN)$ is log-smooth. 
Using the compatibility of the rigid and crystalline (resp. log-crystalline) Chern classes, 
one checks that the square commutes. 
Tensoring the square with $\QQ$ (resp. $\QQ_p$) finishes the proof. 
\qed\medskip

\begin{Remark}
One can give an elementary argument for the surjectivity of 
$c_{1}:\Pic^{\log}(Y)\otimes\QQ_p\rightarrow\Hlogcris{2}(Y/K_0)^{\varphi=p,N=0}$ 
using the Hyodo-Kato complex, but we have chosen to present the above proof because 
it nicely demonstrates the relationship between $\Pic(Y)$ and $\Pic^{\log}(Y)$.
\end{Remark}
 
\subsection{Raskind-admissibility}
Next, we recall the following $p$-adic Lefschetz $(1,1$)-theorem, due to 
Berthelot and Ogus \cite[Theorem 3.8]{Berthelot Ogus} in the smooth case and to 
Yamashita \cite[Theorem 3.1]{Yamashita} in the semi-stable case.

\begin{Theorem}[Yamashita]\label{thm: yamashita}
 Let $X$ be a smooth and proper variety over $K$ and assume that there exists
 a proper and semi-stable model ${\cal X}\to\Spec\OO_K$ of $X$.
 Let ${\cal X}_0$ be the special fibre.
 \begin{enumerate}
 \item There exists a commutative diagram
 $$
 \begin{array}{ccccc}
    \Pic(X) &\leftarrow& \Pic({\cal X}) &\to& \Pic({\cal X}_0) \\
    \parallel && \downarrow && \downarrow \\
   \Pic(X) &=& \Pic^{\rm log}({\cal X}) &\to& \Pic^{\rm log}({\cal X}_0).
 \end{array}
 $$
 where the vertical maps are restrictions and the horizontal maps are specialisations.
 \item An invertible sheaf ${\cal L}\in\Pic({\cal X}_0)\otimes\QQ$ 
  (resp. ${\cal L}\in\Pic^{\rm log}({\cal X}_0)\otimes\QQ$)
  can be lifted to
  $\Pic({\cal X})\otimes\QQ$ (resp. $\Pic^{\rm log}({\cal X})\otimes\QQ$ ) 
  if and only if its first Chern class satisfies
 $$
    c_1({\cal L}) \,\in\, \Fil^1\,\cap\, \left(\Hlogcris{2}({\cal X}_0/K_0)\otimes_{K_0}K\right),
 $$
 where $\Fil^\bullet$ denotes the Hodge filtration on $\HdR{2}(X/K)$.
 \end{enumerate}
\end{Theorem}

\prf
Claim (1) is in the discussion at the beginning of \cite[\S 2]{Yamashita}.
Claim (2) is \cite[Theorem 3.1]{Yamashita}.
\qed\medskip

After these preparations, we make the following key observation.

\begin{Remark}\label{rem: key observation}
 If ${\cal X}_0$ has cohomologically totally degenerate reduction, it looks at first glance
as if the combination of Proposition \ref{prop: logcrys}, Proposition \ref{prop: rational structure}, 
and Theorem \ref{thm: yamashita} might prove Conjecture \ref{conj: raskind}.
However, it is crucial to note that Theorem \ref{thm: yamashita} deals with \emph{$\QQ$-classes} of invertible
sheaves, whereas the other results deal with \emph{$\QQ_p$-classes}.
\end{Remark}

\begin{Definition}\label{def: raskind admissible}
 Let $K$ be a $p$-adic field and let $K_0$ be the maximal unramified extension of $\QQ_p$ inside $K$.
 Let $(H\otimes_{K_0}K,\Fil^\bullet,\varphi,N)$ be a filtered $(\varphi,N)$-module over $K$ 
 and let $V$ be a rational structure on $H$. 
 Then, $H$ is called  \emph{Raskind-admissible} if the natural inclusion
 of $\QQ_p$-vector spaces
 $$
     \left(\Fil^1\,\cap\,V^{\varphi=p,N=0}\right)\otimes_\QQ\QQ_p   \,\subset\, \Fil^1 \,\cap\, H^{\varphi=p,N=0}
 $$
 is an equality.
\end{Definition}

We remark that $V^{\varphi=p,N=0}=B_1$ and $H^{\varphi=p,N=0}=B_1\otimes\QQ_p$ 
in the notation of Definition \ref{def: rational structure}, see also Equation \eqref{eq: phi and N}.
After these preparations, we now reformulate Raskind's conjecture for divisors
(Conjecture \ref{conj: raskind}) into semi-linear algebra.
In fact, the following equivalence holds under a slightly weaker
assumptions than Raskind's requirement of total degeneracy.

\begin{Theorem}\label{thm: raskind translation}
 Let $X$ be a smooth and proper variety over $K$ and assume that there exists
 a proper and semi-stable model ${\cal X}\to\Spec\OO_K$ of $X$.
 Assume that the special fibre ${\cal X}_0$ is cohomologically totally degenerate.
 Then, the following are equivalent:
 \begin{enumerate}
 \item The homomorphism \eqref{tatehom} is surjective for $\ell=p$.
 \item The filtered $(\varphi,N)$-module $\DD_{{\rm st}}(\Het{2}(X_{\overline{K}},\QQ_p))$ together
  with the rational structure arising from ${\cal X}_0$ is Raskind-admissible.
 \end{enumerate}
\end{Theorem}

\prf
By Proposition \ref{prop: rational structure}.(1), there exists a rational structure 
$(V=A\oplus B_0\oplus B_1\oplus C,\varphi_V,N_V)$
on $\DD_{{\rm st}}(\Het{2}(X_{\overline{K}},\QQ_p))$ associated to ${\cal X}_0$.
Next, by Proposition \ref{prop: rational structure}.(3), the first Chern class
induces an isomorphism
\begin{equation}
\label{eq: step one}
   \Pic({\cal X}_0)\otimes_\ZZ\QQ \,\cong\, B_1.
\end{equation}
Using the first Chern class on $X$ 
we obtain a map
\begin{equation}
\label{eq: step two}
\Pic(X)\otimes_\ZZ\QQ \,\to\, \Fil^1\,\cap\, \Hlogcris{2}({\cal X}_0/K_0)^{\varphi=p,N=0},
\end{equation}
whose image lies inside the subspace $B_1$ of $\Hlogcris{2}({\cal X}_0/K_0)^{\varphi=p,N=0}$.
Using Yamashita's theorem (Theorem \ref{thm: yamashita}),
it follows that the first map in
\begin{equation}
 \label{eq: step rational}
\Pic(X)\otimes_\ZZ\QQ \,\to\, \Fil^1\,\cap\, B_1\,\cong\,\Fil^1\,\cap\, \Pic({\cal X}_0)\otimes\QQ
\end{equation}
is surjective.

By Proposition \ref{prop: logcrys} and \eqref{eq: phi and N}, we have
\begin{eqnarray*}
  \Het{2}(X_{\overline{K}},\QQ_p(1))^{G_K} &\cong& \Fil^1\,\cap\, \Hlogcris{2}({\cal X}_0/K_0)^{\varphi=p,N=0} \\
  &=& \Fil^1\cap (B_1\otimes\QQ_p).
\end{eqnarray*}
Combining this with equation \eqref{eq: step two}, we see that the homomorphism \eqref{tatehom} 
is surjective for $\ell=p$ 
if and only if
$$
    \Pic(X)\otimes_\ZZ\QQ_p \,\to\, \Fil^1\,\cap\, (B_1\otimes_\QQ\QQ_p)
$$
is surjective.
In view of \eqref{eq: step rational}, this is equivalent to asking whether 
$$
   \left(\Fil^1\,\cap\, B_1\right) \otimes_\QQ\QQ_p \,\to\, \Fil^1\,\cap\, (B_1\otimes_\QQ\QQ_p)
$$
is surjective.
But this is equivalent to the rational structure $V$ on the filtered 
$(\varphi,N)$-module $\DD_{{\rm st}}(\Het{2}(X_{\overline{K}},\QQ_p))$ 
being Raskind-admissible.
\qed\medskip

 \subsection{Ordinary representations}
Using the Hyodo-Kato complex and work of 
Perrin-Riou \cite{Perrin-Riou} and
Illusie \cite{Illusie appendix}, we have the following description
of the interplay between the $p$-adic Galois represenation of
$G_K$ on $\Het{2}(X_{\overline{K}},\QQ_p)$
and the special fibre of a semi-stable model.

 \begin{Proposition}\label{prop: illusie}
  Let $X$ be a smooth and proper variety over $K$ and assume that there exists
   a proper and semi-stable model ${\cal X}\to\Spec\OO_K$ of $X$.
   Assume that the special fibre ${\cal X}_0$ is cohomologically 
   totally degenerate
   and that every component of ${\cal X}_0^{[i]}$ for every $i$
   is ordinary in the sense of 
   Bloch-Kato-Illusie-Raynaud.
   \begin{enumerate}
    \item The $G_K$-representation on $\Het{2}(X_{\overline{K}},\QQ_p)$
     is ordinary in the sense of Perrin-Riou \cite[1.2]{Perrin-Riou}.
     More precisely, if ${\cal F}^\bullet$ denotes the corresponding
     filtration, then there exist $G_K$-equivariant isomorphisms
     $$
        {\rm gr}^{2-i} \Het{2 }\left(X_{\overline{K}},\QQ_p\right) \,\cong\, 
        H^{2-i}\left( ({\cal X}_0)_{\overline{k}}, W\omega_{\rm log}^i\right)\otimes\QQ_p(-i)\,.
     $$
    \item Let $(V=A\oplus B_0\oplus B_1\oplus C,\varphi_V,N_V)$ 
     be the rational structure on 
     $\DD_{\rm st}(\Het{2}(X_{\overline{K}},\QQ_p))$ associated 
     to ${\cal X}_0$.
     Then,
     $$
      \begin{array}{llcl}
         A&\otimes_\QQ\QQ_p &\cong& H^{2}( {\cal X}_0, W\omega_{\rm log}^0)\otimes_{\ZZ_{p}}\QQ_{p} \\
        (B_0\oplus B_1)&\otimes_\QQ\QQ_p &\cong& H^{1}( {\cal X}_0, W\omega_{\rm log}^1)\otimes_{\ZZ_{p}}\QQ_{p}  \\
        C&\otimes_\QQ\QQ_p &\cong&  H^{0}( {\cal X}_0, W\omega_{\rm log}^2)\otimes_{\ZZ_{p}}\QQ_{p}\,.
      \end{array}
     $$
     \item We have the inequality
     $$ 
          \rho(X) \,\leq\, h^{1,1}(X) \,-\, h^{0,2}(X),
     $$
     where $\rho(X)$ denotes the Picard rank of $X$.
   \end{enumerate}
 \end{Proposition}
 
\prf
Since every component of ${\cal X}_0^{[i]}$ is ordinary for every $i$, 
so is ${\cal X}_0$ itself  \cite[Proposition 1.6]{Illusie appendix}.
Thus, claim (1) follows from \cite[Corollaire 2.7]{Illusie appendix}.
 
More precisely, we obtain a very explicit description of the filtered
$(\varphi,N)$-module $\DD_{\rm st}(\Het{2}(X_{\overline{K}},\QQ_p))$
via the rational structure $V$ and 
\cite[Corollaire 2.7]{Illusie appendix} provides us with an equally explicit
description of the $G_K$-action on $\Het{2}(X_{\overline{K}},\QQ_p)$
via $H^{2-i}( ({\cal X}_0)_{\overline{k}}, W\omega_{\rm log}^i)$.
Comparing these two descriptions, claim (2) follows.
 
Finally, the de Rham Chern class map 
$c_{1}:\Pic(X)\rightarrow \HdR{2}(X/K)$ is compatible with the log-crystalline Chern class 
and its image lies inside $\Fil^1$ and $\Hlogcris{2}({\cal X}_0/K_0)^{\varphi=p,N=0}$ 
(see, for example, \cite[\S2]{Yamashita}).
In our situation, the latter is isomorphic to $B_1\otimes\QQ_p$, from which claim (3) 
immediately follows (noting that $h^{i,j}(X)=\dim_{K_{0}}H^{j}(\mathcal{X}_{0},W\omega_{\mathcal{X}_{0}}^{i})\otimes_{W(k)}K_{0}=\dim_{\QQ_{p}}H^{j}( {\cal X}_0, W\omega_{\rm log}^i)\otimes_{\ZZ_{p}}\QQ_{p}$ by ordinarity).
 \qed\medskip
 
Thus, when tensored with $\QQ_p$, the rational structure on $\DD_{\rm st}(\Het{2}(X_{\overline{K}},\QQ_p))$
arising from ${\cal X}_0$ has an interpretation via the logarithmic Hodge-Witt cohomology groups of the special fibre.
Although this is not directly related to the conjectures discussed in
this article, it might be of independent interest.

\subsection{Concluding remarks}
The definitions and notions of this section are a little bit ad hoc, since we only 
deal with Conjecture \ref{conj: raskind}.
A more conceptual approach, which would be needed when studying Raskind's 
conjectures for cycles of higher codimension, could proceed along the following lines:
 \begin{enumerate}
  \item One can directly construct $\QQ$-structures on the groups
   $\Hlogcris{\ast}(Y/K_0)$ using the Chow complex of \cite{BGS}.
   Moreover, these would also come with $\QQ$-linearisations
   of the Frobenius $\varphi$ and the monodromy $N$.
  \item Concerning the definitions: one would have to define a \emph{weight} $w$ of a filtered $(\varphi,N)$-module
  (in Definition \ref{def: rational structure}, it would be $w=2$), one would have to define such a module
  to be of \emph{maximal nilpotent monodromy} if the monodromy operator $N$ satisfies $N^w\neq0$
  and then, a \emph{rational structure} would be a $\QQ$-vector space $V$ with a direct sum
  decomposition into subspaces upon which $\varphi_V$ acts as multiplication by $p^i$ for 
  $i=0,...,w$, etc.
 \item For an equivalence as in Theorem \ref{thm: raskind translation}, one would
   also need a version of Yamashita's theorem (Theorem \ref{thm: yamashita}) for deforming 
   cycles of higher codimension.
   This would be a semi-stable analogue of the $p$-adic variational Hodge conjecture (see for example \cite[Conjecture 1.2]{BEK14} for the good reduction case, where the conjecture is attributed to Fontaine-Messing). 
   This conjecture is open in codimension $\geq 1$, even in the case of good reduction, 
 but see \cite{BEK14} for the state of the art.    
\end{enumerate}
To keep the discussion in this section shorter, we have decided not to develop
the setup in general, but to stick to the case of divisors.

\section{Abeloid varieties}
\label{sec: abelian varieties}

From this section on, we study abelian varieties over $p$-adic fields with totally
degenerate reduction.
More precisely, we describe their morphisms, $\ell$-adic Tate modules,
and the filtered $(\varphi,N)$-modules associated to the latter via
$p$-adic uniformisation, that is, within the framework of abeloid varieties.
The results of this section might be of independent interest and some of them
 might already be known to the experts.
 
\subsection{Generalities}
We start with the behaviour of Question \ref{question: tate} under field extension
and under completion, similar to what we did for Raskind's conjecture for divisors 
(Conjecture \ref{conj: raskind}) 
in Section \ref{subsec: Raskind field extn} and Section \ref{sec: relation to Tate}.

Let $F$ be a field of characteristic $p\geq0$ and let $\ell$ be a 
prime, possibly equal to $\ell$.
If $A$ is an abelian variety of dimension $g$ over a field $F$, then the 
$\ell$-adic Tate module $T_\ell$
and the rational $\ell$-adic Tate module $V_\ell$ of $A$ 
are defined to be to be
$$
  T_\ell (A) \,:=\, \varprojlim_n A(F^{\rm sep})[\ell^n]
  \mbox{ \quad and \quad }
  V_\ell(A) \,:=\, T_\ell(A)\otimes_{\ZZ}\QQ
$$
together with their $G_F$-actions. 
If $\ell\neq p$, then $T_\ell(A)$ is a free $\ZZ_\ell$-module of rank $2g$.

\begin{Proposition}
 \label{fieldextension abelian}
  Let $A$ and $B$ be abelian varieties over a field $F$,
  let $F\subset F'$ be a finite Galois extension
  and let $\ell$ be a prime.
  If the map \eqref{tatehom2} (please see \S\ref{abelian and abeloid}) is surjective with r
  espect to $A_{F'}$, $B_{F'}$, $G_{F'}$, and $\ell$, then 
   \eqref{tatehom2} is surjective with respect to $A$, $B$, $G_F$, and $\ell$.
\end{Proposition}

\prf 
Suppose that \eqref{tatehom2} with respect to $A_{F'}$, $B_{F'}$, $G_{F'}$ is a surjection. 
Then, we get a surjection on the $\mathrm{Gal}(F'/F)$-invariants
\begin{equation*}\resizebox{1.0\hsize}{!}{$
(\Hom(A_{F'},B_{F'})\otimes_{\ZZ}\ZZ_{\ell})^{\mathrm{Gal}(F'/F)} \,\twoheadrightarrow\,
\Hom_{G_{F'}}(T_{\ell}(A),T_{\ell}(B))^{\mathrm{Gal}(F'/F)}=\Hom_{G_{F}}(T_{\ell}(A),T_{\ell}(B)) \,, $}
\end{equation*} 
and the left-hand side is $\Hom(A,B)\otimes_{\ZZ}\ZZ_{\ell}$ by Galois descent for morphisms.
\qed\medskip

\begin{Remark}
 For abeloid varieties over $p$-adic fields, we will see a second proof of this result in
 Corollary \ref{cor: Tate finite field extn} below.
\end{Remark}

\begin{Proposition}\label{prop: relation to Tate2}
  Let $A$ and $B$ be abelian varieties over a number field $F$
  and assume that there exists a  prime $\ell$ and
  a prime ideal $\idealp\subset\OO_F$ such that
  \eqref{tatehom2} is surjective for $A_{F_\idealp}$, $B_{F_\idealp}$,
  $G_{F_\idealp}$, and $\ell$.
  Then, \eqref{tatehom2} is surjective for $A$, $B$,
  $G_F$, and $\ell$.
\end{Proposition}

\prf
We view $G_{F_{\mathfrak{p}}}$ as a subgroup of $G_{F}$. 
Since we are allowed to make finite Galois extensions by Proposition \ref{fieldextension abelian}, we may assume that 
\begin{equation*}
\mathrm{Hom}(A,B)=\mathrm{Hom}(A_{\overline{F}},B_{\overline{F}}) \,.
\end{equation*}
In particular, $G_{F_{\mathfrak{p}}}$ and $G_{F}$ act trivially on $\mathrm{Hom}(A_{\overline{F}},B_{\overline{F}})$. 
Then we have the following commutative square
\begin{center}
\begin{tikzpicture}[descr/.style={fill=white,inner sep=1.5pt}]
        \matrix (m) [
            matrix of math nodes,
            row sep=2.5em,
            column sep=2.5em,
            text height=1.5ex, text depth=0.25ex
        ]
        { \mathrm{Hom}(A_{\overline{F}},B_{\overline{F}})\otimes_{\mathbb{Z}}\mathbb{Z}_{\ell} & \mathrm{Hom}(T_{\ell}(A),T_{\ell}(B))\\
          \mathrm{Hom}(A_{\overline{F_{\mathfrak{p}}}},B_{\overline{F_{\mathfrak{p}}}})\otimes_{\mathbb{Z}}\mathbb{Z}_{\ell} & \mathrm{Hom}(T_{\ell}(A_{F_{\mathfrak{p}}}),T_{\ell}(B_{F_{\mathfrak{p}}})) \\
        };

       \path[overlay,->, font=\scriptsize]
       (m-1-1) edge node [right]{$\cong$} (m-2-1)
       (m-1-2) edge (m-2-2)
        (m-1-1) edge (m-1-2)
        (m-2-1) edge (m-2-2);    
        \end{tikzpicture} 
\end{center}
Taking $G_{F_{\mathfrak{p}}}$-invariants gives 
\begin{center}
\begin{tikzpicture}[descr/.style={fill=white,inner sep=1.5pt}]
        \matrix (m) [
            matrix of math nodes,
            row sep=2.5em,
            column sep=2.5em,
            text height=1.5ex, text depth=0.25ex
        ]
        { \mathrm{Hom}(A,B)\otimes_{\mathbb{Z}}\mathbb{Z}_{\ell} & \mathrm{Hom}_{G_{F_{\mathfrak{p}}}}(T_{\ell}(A),T_{\ell}(B))\\
          \mathrm{Hom}(A_{F_{\mathfrak{p}}},B_{F_{\mathfrak{p}}})\otimes_{\mathbb{Z}}\mathbb{Z}_{\ell} & \mathrm{Hom}_{G_{F_{\mathfrak{p}}}}(T_{\ell}(A_{F_{\mathfrak{p}}}),T_{\ell}(B_{F_{\mathfrak{p}}})) \\
        };

       \path[overlay,->, font=\scriptsize]
       (m-1-1) edge node [right]{$\cong$} (m-2-1)
       (m-1-2) edge (m-2-2)
       (m-1-1) edge (m-1-2);    
        
       \path[overlay,->>, font=\scriptsize]
       (m-2-1) edge (m-2-2);
         
\end{tikzpicture} 
\end{center}
where the lower horizontal arrow is a surjection by assumption. 
We deduce therefore that the upper horizontal arrow is a surjection. 
This proves the proposition because of the inclusion
$\Hom_{G_{F}}(T_{\ell}(A),T_{\ell}(B))\subset\Hom_{G_{F_{\idealp}}}(T_{\ell}(A),T_{\ell}(B))$.
\qed\medskip

\subsection{Degenerations of abelian varieties}\label{subsec: Raskind}\label{subsec: degeneration AV}
In the simple cases treated in Proposition \ref{prop: simple case}, no
assumption on the degeneration was needed and in 
Section \ref{sec: translation}, we worked with a weak form of total
degeneration.
In \cite{Raskind, RX}, Raskind suggested a degeneration assumption which is rather involved.
For the purposes of this article, the following slight generalisation
of \cite[Example 1.(i)]{RX} suffices.

\begin{Lemma}[Raskind-Xarles + $\varepsilon$]\label{lem: raskind}
 Let $Y=\bigcup_{i\in I} Y_i$ be a strict normal crossing scheme 
 over a perfect field $F$.
 Assume that for every subset $J\subset I$ the intersection $Y_J$,
 if non-empty, is isomorphic to
 a successive blow-up of smooth, projective, and toric varieties along subvarieties
 that are also smooth and toric.
 Then, $Y$ is totally degenerate in the sense of Raskind.
\end{Lemma}

\prf
If every $Y_J$ is a smooth, projective, and toric variety, then this is
\cite[Example 1.(i)]{RX}.
Quite generally, if $\widetilde{X}$ is the blow-up of a smooth variety $X$ along 
a smooth subvariety $Y$, then one can express the $\ell$-adic cohomology groups, 
the crystalline cohomology groups, the Chow groups, and the cycle class maps 
for $\widetilde{X}$ in terms of $X$ and $Y$.
From these formulas, it follows that the requirements (a)-(c) of \cite[Definition 1]{RX} also
hold for our assumptions.
Moreover, if $X$ and $Y$ are ordinary in the sense of Bloch-Kato-Illusie-Raynaud, 
then so is $\widetilde{X}$ \cite[Proposition 1.6]{I90}, that is, also requirement (d) is satisfied.
\qed\medskip

\begin{Remark}\label{rem: CL}
 If $Y$ is a strict normal crossing scheme of dimension $N\leq2$, 
such that all the $Y_J$ are smooth and rational varieties, then the assumptions
of the lemma are satisfied.
For example, this applies to the combinatorial degenerations of type III 
of the surfaces from \cite[Definitions 5.4 to 5.7]{CL}.
\end{Remark}

Next, we discuss the notion of totally degenerate reduction 
for abelian varieties.
There are several obvious candidates, all of which are stable under finite field extensions
and all of which are equivalent up to finite field extensions.
The following is well-known, but maybe never explicitly stated in this way, which is why we 
include a short discussion with references.

\begin{Proposition}
  \label{prop: abelian totally degenerate}
  Let $\OO_K$ be a local and Henselian DVR with field of fractions $K$ and residue field $k$.
  Assume that $K$ is of characteristic zero and that $k$ is perfect of characteristic $p\geq0$.
  Let $A$ be an abelian variety of dimension $g\geq1$ over $K$.
  Consider the following properties:
  \begin{enumerate}
   \item $A$ admits uniformisation in the sense of Mumford \cite{Mumford}.
   \item The connected component of the special fibre of the N\'eron model of $A$ is a split torus.
   \item The special fibre of the projective regular 
    K\"unnemann-Mumford model \cite{Mumford, Kunnemann}
    of $A$ is a union of smooth and toric varieties.
   \item $A$ has totally degenerate reduction in the sense of Raskind \cite{Raskind}.
   \item For some (resp., all) $\ell\neq p$, the action of the inertia subgroup
     $I_K$ of $G_K$ on
     $\Het{1}(A_{\overline{K}},\QQ_\ell)$ is unipotent with maximal number
     of $2\times2$-Jordan blocks, that is, its Jordan normal form has 
     $g$ Jordan blocks of size $2\times 2$
     and generalised eigenvalue $1$.
    \item For some (resp., all) $\ell\neq p$, the $I_K$-action on
     $\Het{g}(A_{\overline{K}},\QQ_\ell)$ is maximally unipotent, that is, 
     its Jordan normal form has one Jordan block of size $g\times g$ and generalised
     eigenvalue $1$.
    \item Assume $g\geq2$. For some (resp., all) $\ell\neq p$, the $I_K$-action on
     $\Het{2}(A_{\overline{K}},\QQ_\ell)$ is unipotent and its Jordan normal form has
     no Jordan block of size $2\times2$ and $\frac{1}{2}g(g-1)$ Jordan blocks
     of size $3\times 3$.
   \end{enumerate}
   If $K$ is a $p$-adic field, consider also the following properties
   \begin{enumerate}
     \setcounter{enumi}{7}
      \item The $G_K$-representation $\Het{1}(A_{\overline{K}},\QQ_p)$ is semi-stable and the 
     monodromy operator $N$ of the associated filtered $(\varphi,N)$-module
     is nilpotent with minimal kernel, that is, $\ker(N)$ is $g$-dimensional. 
      \item The $G_K$-representation $\Het{g}(A_{\overline{K}},\QQ_p)$ is semi-stable and the 
     monodromy operator $N$ of the associated filtered $(\varphi,N)$-module
     is maximally nilpotent, that is, $N^g\neq0$.
     \item Assume $g\geq2$. The  $G_K$-representation
     $\Het{2}(A_{\overline{K}},\QQ_p)$ is semi-stable and the monodromy
     operator $N$ of the associated filtered $(\varphi,N)$-module has
     no Jordan block of size $2\times2$ and $\frac{1}{2}g(g-1)$ Jordan blocks
     of size $3\times 3$. 
  \end{enumerate}
  Then, 
   \renewcommand{\labelenumi}{(\roman{enumi})}
   \begin{enumerate}
  \item these properties are {\em stable under finite extension}, that is, if $A$ satisfies one of 
  these properties, then $A_L$ satisfies the same property for every finite field extension $K\subset L$.
  \item Moreover, these properties are {\em equivalent up to finite extension}, that is, if $A$ satisfies 
  one of the above properties,
  then there exists a finite field extension $K\subset L$ such that $A_L$ satisfies all of 
  these properties.
  \end{enumerate}
    \renewcommand{\labelenumi}{(\arabic{enumi})}
\end{Proposition}

\prf
The assertion that all these properties are stable under finite extension are well-known or easy and we leave
them to the reader.

The fact that (1) and (2) 
are equivalent up to finite extension is the main result of \cite{Mumford}, see also the
discussion in \cite[\S 5.6]{Lutkebohmert}.

The fact that (2) and (3) are equivalent up to finite extension is a special case of \cite[\S3]{Kunnemann}.

If an abelian variety satisfies (4), then the identity component of the N\'eron model is semi-abelian 
and the triviality of the conditions on the Chow group imply that it cannot have abelian parts, and thus, 
(4) is equivalent to (2) and (3) up to finite extension.
Conversely, if an abelian variety satisfies (3), then Lemma \ref{lem: raskind} applies and thus, 
(3) is equivalent to (4) up to finite extension, see also \cite[Example 1.(ii)]{RX}.

Assume that the N\'eron model of $A$ has semi-abelian reduction, let $t$ be the dimension 
of the toric part and $a$ be the dimension of the abelian part (and thus, $g=t+a$).
Then, the $I_K$-action on $\Het{1}(A_{\overline{K}},\QQ_\ell)$
is unipotent with index of unipotency at most $\leq2$. 
More precisely, there are $2a$ Jordan blocks of size $1\times 1$ and generalised
eigenvalue $1$ and there are $t$ Jordan blocks of size $2\times 2$ and generalised eigenvalue $1$.
We refer to \cite[Theorem 6 in Chapter 7.4]{BLR} or \cite[Expos\'e IX]{SGA7} for proofs and details.
Similarly, the $G_K$-representation on $\Het{1}(A_{\overline{K}},\QQ_p)$ is semi-stable
in this case.
The monodromy operator on the associated $(\varphi,N)$-module is nilpotent
with nilpotency at most $\leq2$.
More precisely, there are $2a$ Jordan blocks of size $1\times 1$ and
eigenvalue $0$ and there are $t$ Jordan blocks of size $2\times 2$ and generalised 
eigenvalue $0$.

From this discussion, it follows that (2) is equivalent to (5) up to finite extension and
that (2) is equivalent to (8) up to finite extension.

In any case and for all $i$, there exist $G_K$-equivariant isomorphisms between
$\Het{i}(A_{\overline{K}},\QQ_\ell)$ and $\wedge^i \Het{1}(A_{\overline{K}},\QQ_\ell)$.
Moreover, for all $\ell\neq p$, the $I_K$-actions on $\Het{i}(A_{\overline{K}},\QQ_\ell)$ 
are quasi-unipotent and the $G_K$-representations $\Het{i}(A_{\overline{K}},\QQ_p)$
are potentially semi-stable.
Replacing $K$ by a finite extension, we may assume that all these $I_K$-actions are
unipotent for $\ell\neq p$ and all the $G_K$-representations are semi-stable for $\ell=p$.

In particular, if $\ell\neq p$, then the $I_K$-action on $\wedge^g \Het{g}(A_{\overline{K}},\QQ_\ell)$
is unipotent with order of unipotency at most $g$.
More precisely, there exists at most one Jordan block of size $g\times g$ and generalised
eigenvalue $1$ and there exists such a block if and only if $g=t$.
We leave the straight forward exercise in linear algebra to the reader.
This implies that (5) is equivalent to (6) up to finite extension.
Similarly, (8) is equivalent to (9) up to finite extension.

If $g\geq2$ and $\ell\neq p$, then the $I_K$-inertia on $\wedge^2\Het{1}(A_{\overline{K}},\QQ_\ell)$ 
is unipotent with order of unipotency at most $3$.
More precisely, there are $r:=\frac{1}{2}t(t-1)$ Jordan blocks
of size $3\times 3$ with generalised eigenvalue $1$ and $s:=2at$ Jordan blocks
of size $2\times 2$ with generalised eigenvalue $1$.
Then, $t=\frac{r+4s}{2g-2}$ and thus, the $I_K$-action on 
$\wedge^2\Het{1}(A_{\overline{K}},\QQ_\ell)$ encodes $a$ and $t$.
Again, we leave the details to the reader.
These considerations show that (5) is equivalent to (7) up to finite extension.
Similarly, (9) is equivalent to (11) up to finite extension.
\qed\medskip

\subsection{Abeloid varieties}\label{subsec: abeloid varieties}
From now on, we will adopt the point of view of $p$-adic uniformisation.
In particular, this allows us to work with abeloid varieties, which are rigid analytic
varieties over $p$-adic fields that are not necessarily algebraisable.
We now establish and recall a couple of general facts about
abeloid varieties.

Let $q_1,...,q_g\in K^{\times g}$ be vectors, say $q_i=(q_{i,1},...q_{i,g})$,
such that $\nu_p(q_{i,j})>1$ for all $i,j$.
Let $Q=(q_{i,j})\in M_{g\times g}(K)$ be the $g\times g$-matrix, whose rows 
are the $q_i$.
Associated to $Q$, we have the matrix
\begin{align*}
  \ord_p(Q) &\,:=\, \left( \nu_p(q_{i,j}) \right) \,\in\,M_{g\times g}(\QQ).
\end{align*}
By definition, the abelian subgroup $\Lambda=q_1^\ZZ\cdots q_g^\ZZ$ of $K^{\times g}$ is called
a  \emph{lattice}, if the columns of $\ord_p(Q)$ span a lattice inside $\RR^g$ or, equivalently,
if the matrix $\ord_p(Q)$ is invertible.
(In the literature, this matrix is sometimes constructed with respect to a valuation
$\nu_K$ on $K^\times$ with $\nu_K(\pi_K)=1$ for some uniformiser $\pi_K\in\OO_K$ and then,
the matrix has integer entries rather then rational ones.
We have decided to work with the valuation $\nu_p$ instead, which rescales
the classical matrix by $\nu_p(\pi_K)$ and has the advantage of being stable 
under finite field extensions of $K$.)

Associated to a lattice $\Lambda\subset K^{\times g}$, there is a rigid analytic variety 
over $K$,
the \emph{abeloid variety} $\GG_{m,K}^g/\Lambda$, see \cite[Chapter 7]{Lutkebohmert}.
The $g\times g$-matrix $Q$ associated to a choice of basis for $\Lambda$ is called a
\emph{period matrix}.
The algebraisable abeloid varieties are precisely the totally degenerating abelian varieties
studied by Mumford \cite{Mumford}.
Moreover, if $g=1$, then a lattice $\Lambda\subset K^\times$ is generated by a single element
$q\in K^\times$ with $\nu_p(q)>1$, an abeloid variety of dimension one is always algebraisable, and 
these are precisely the \emph{Tate elliptic curves}.

We introduce the following notation:
 if $R$ is a commutative ring and if $A$ is an $R$-module, then
the set ${\rm Mat}_{m\times n}(A)$ of $m\times n$-matrices with 
values in $A$ is an Abelian group.
Let $X\in{\rm Mat}_{m\times n}(A)$.
Then, if $M\in{\rm Mat}_{s\times m}(R)$ and $N\in {\rm Mat}_{n\times t}(R)$ are matrices 
with entries in $R$ for some $s,t$, then the matrix products $M\odot X$ and $X\odot N$ are defined.
In particular, ${\rm Mat}_{n\times n}(A)$ is a ${\rm Mat}_{n\times n}(R)$-bimodule.
In the next theorem, we have $R=\ZZ$ and $A=K^\times$.

\begin{Theorem}[Gerritzen $+\varepsilon$]
 \label{thm: abeloid morphisms}
 Let $\Lambda_A\subset K^{\times g}$ and $\Lambda_B\subset K^{\times h}$
 be lattices and let $A:=\GG_{m,K}^g/\Lambda_A$ and
 $B:=\GG_{m,K}^h/\Lambda_B$ be the associated abeloid varieties.
 Let $Q_A$ and $Q_B$ be period matrices for $\Lambda_A$ and $\Lambda_B$.
 Then, there exist isomorphisms
  {\small
  \begin{eqnarray*} 
  & & \Hom(A,B)\\
   &\cong&
    \left\{ M\in{\rm Mat}_{g\times h}(\ZZ) \,|\, \Lambda_A\odot M \,\subseteq\, \Lambda_B
   \right\} \\
   &\cong&
   \left\{ M\in{\rm Mat}_{g\times h}(\ZZ) \,|\, \exists N\in{\rm Mat}_{h\times g}(\ZZ)\,:\,
    Q_A\odot M \,=\, N\odot Q_B
   \right\} \\
   &\cong&
   \left\{ M\in{\rm Mat}_{g\times h}(\ZZ) \,|\, 
    \left(\ord_p(Q_A)^{-1}\odot Q_A\right) \odot M \,=\, M\odot \left(\ord_p(Q_B)^{-1}\odot Q_B\right)
   \right\} \,.
 \end{eqnarray*}}
 In particular,
 \begin{enumerate}
   \item the natural map
   $$
    \Hom(A,B) \,\to\, \Hom(A_{\overline{K}},B_{\overline{K}})
   $$
   is an isomorphism of abelian groups.
   In particular, the $G_K$-action on the right hand side is trivial.
 \item $A$ and $B$ are $K$-isogenous if and only if they are $\overline{K}$-isogenous.
 \item $A$ is $K$-simple if and only if it is $\overline{K}$-simple.
 \end{enumerate}
\end{Theorem}

\prf
If $M=(m_{ij})\in M_{g\times h}(\ZZ)$ is a $g\times h$-matrix, then
it gives rise to a map $\psi_M:K^{\times g}\to K^{\times h}$ by
sending
$$
  (x_1,...,x_g) \,\mapsto\, \left( ..., \prod_{j=1}^g x_j^{m_{ji}},...\right).
$$
If $\psi_M(\Lambda_A)=\Lambda_A\odot M$ is contained in $\Lambda_B$, then $\psi_M$ descends
to a morphism $A\to B$ of abeloid varieties.
Conversely, every morphism of abeloid varieties is of this form
by the main theorem of \cite{Gerritzen}, see also the discussion in
 \cite[\S3]{Ka07}.
This establishes the first two isomorphisms describing of $\Hom(A,B)$.
Taking valuations in the equality $Q_A\odot M \,=\, N\odot Q_B$, we find
$$
   N\,=\, \ord_p(Q_A)\cdot M \cdot \ord_p(Q_B)^{-1},
$$
which implies the third isomorphism.

Now, since every morphism $A_{\overline{K}}\to B_{\overline{K}}$
can be defined over some finite field extension $L/K$ and since
the above description of homomorphisms is also valid for
homomorphisms over $L$, it follows from this description
that every homomorphism $\Hom(A_L,B_L)$ can be defined
over $K$.
This establishes claim (1).
In particular, $A$ and $B$ are isogenous over $K$ if and only if they
are isogenous over $\overline{K}$, which establishes claim (2).
Finally, $A$ is simple over $K$ if and only if there exists
no non-trivial idempotent in $\End(A)$ if and only if
there exists no non-trivial idempotent in $\End(A_{\overline{K}})$
(by the already established (1)) if and only if $A_{\overline{K}}$
is simple, which establishes claim (3).
\qed\medskip

We will give another description of $\Hom(A,B)\otimes\QQ$ in terms of $\MazurL$-invariants
in Proposition \ref{prop: abeloid morphisms via L} below.

\begin{Corollary}\label{cor: Tate finite field extn}
  Let $A$ and $B$ be abeloid varieties over a $p$-adic field $K$,
  let $L/K$ be a finite field extension, and let $\ell$ be a prime.
  If \eqref{tatehom2} is surjective for 
  $A_L$, $B_L$, and $\ell$, then \eqref{tatehom2} is
  surjective for $A$, $B$, and $\ell$.
\end{Corollary}

\prf
We have a commutative diagram
\begin{center}
\begin{tikzpicture}[descr/.style={fill=white,inner sep=1.5pt}]
        \matrix (m) [
            matrix of math nodes,
            row sep=2.5em,
            column sep=2.5em,
            text height=1.5ex, text depth=0.25ex
        ]
        { \Hom(A,B)\otimes\ZZ_\ell & \Hom_{G_K}\left(T_\ell(A),T_\ell(B)\right) \\
        \Hom(A_L,B_L)\otimes\ZZ_\ell & \Hom_{G_L}\left(T_\ell(A_L),T_\ell(B_L)\right) \\
        };

       \path[overlay,->, font=\scriptsize]
       (m-1-1) edge (m-1-2)
       (m-2-1) edge (m-2-2)
       (m-1-1) edge (m-2-1)  
       (m-1-2) edge (m-2-2);
              
        \end{tikzpicture} 
\end{center}
By Theorem  \ref{thm: abeloid morphisms}.(1), the left vertical arrow is an isomorphism, from which the 
statement immediately follows.
\qed\medskip

\subsection{The Tate module of an abeloid variety}
Let $A=\GG_{m,K}^g/\Lambda$ be an abeloid variety over a $p$-adic field $K$ 
and let $\ell$ be a prime, possibly equal to $p$.
It follows from the rigid analytic parametrisation 
$\overline{K}^{\times ^g}/\Lambda\cong A(\overline{K})$
that the Tate module $T_\ell(A)$ sits in a short exact sequence
\begin{equation}
 \label{eq: abeloid extension}
  0\,\to\, \ZZ_\ell(1)^g \,\to\, T_\ell(A) \,\to\, \ZZ_\ell^g\,\to\,0
\end{equation}
that is compatible with the $G_K$-actions.

To describe its extension class, we follow and generalise some results due to Serre 
\cite[Appendix A]{Serre}.
To state the result, we define $\mu_{\ell^\infty}(K)$ to be the group of those roots of unity of $K$ whose order is a power of $\ell$, 
we choose a uniformiser $\pi_K$, and we denote by 
$U^{(1)}:=1+\idealm_K=1+\pi_K\cdot\OO_K$ the group of $1$-units of $\OO_K$.
Then, there exists an isomorphism of abelian groups
$$
    K^\times \,\cong\, \mu(K)\,\times\,\pi_K^\ZZ\,\times U^{(1)}.
$$
In the sequel, we consider the $\ell$-adic completion
\begin{equation}
    \label{eq: gamma}
     \gamma_\ell\,:\,K^{\times} \,\to\,  \widehat{K^\times}^\ell \,:=\, \varprojlim_n K^{\times}/K^{\times \ell^n}.
\end{equation}
of $K^\times$.

\begin{Lemma}[Serre $+\varepsilon$]
  \label{lemma: Serre}
  Let $A:=\GG_{m,K}^g/\Lambda$ be an abeloid variety over $K$.
  \begin{enumerate}
  \item There exists an isomorphism
  $$
     \ker(\gamma_\ell) \,=\,\left\{ \begin{array}{ll}
      \mu_{\ell^\infty}(K)\times U^{(1)} &\mbox{if $\ell\neq p$,}\\
      \mu_{p^\infty}(K) &\mbox{if $\ell= p$.}\\
     \end{array}\right.
  $$
  In particular, $\ker(\gamma_\ell)$ is finite if and only if $\ell=p$.
  \item
  Taking Galois invariants in \eqref{eq: abeloid extension}, the
  boundary homomorphism in Galois cohomology gives rise
  to a homomorphism
  $$
   d_g \,:\, H^0(G_K,\ZZ_\ell^g) \,\to\, H^1(G_K,\ZZ_\ell(1)^g).
 $$
 Let $e_i=(0,...,1,0,...)$, $i=1,...g$ be the standard basis of 
 $\ZZ_\ell^g$.
 Then, the $\ZZ$-span $\Lambda'$ of $\{d_g(e_i)\}_{i=1,...,g}$
 determines the extension class of \eqref{eq: abeloid extension}.
 \item Kummer theory induces  an isomorphism
 $$ 
            (\widehat{K^\times}^\ell)^g \,\cong\, H^1(G_K,\ZZ_\ell(1)^g).
 $$
 Under this isomorphism, $\gamma_\ell(\Lambda)$ is equal to $\Lambda'$ from 
 assertion (2). 
 \item The image $\gamma_\ell(\Lambda)$ is a lattice, that is, a free $\ZZ$-module of rank $g$.
  In particular, the sequence \eqref{eq: abeloid extension} does not split.
  In fact, there does not even exist a non-trivial and $G_K$-equivariant 
  homomorphism $\ZZ_\ell\to T_\ell(A)$.
 \end{enumerate}
\end{Lemma}

\prf
The description of the kernel in claim (1) for $\ell=p$ is shown in the proof of the implication 
$(3)\Rightarrow (2)$ of the theorem of \cite[Appendix A.1.4]{Serre}.
If $\ell\neq p$, the valuation argument of loc.cit. still shows that $\pi_K^\ZZ$ has trivial intersection
with $\ker(\gamma_\ell)$.
This shows that $\ker(\gamma_\ell)$ is contained in $\tau(k^\times)\times U^{(1)}$,
where $\tau$ denotes the Teichm\"uller lift from $k^\times$ to $K_0^\times\subset K^\times$.
Hensel's lemma implies that $U^{(1)}\subset\ker(\gamma_\ell)$.
The intersection of $\tau(k^\times)$ with $\ker(\gamma_\ell)$ is $\mu_{\ell^\infty}(K)$.

Claims (2), (3), and (4) for $g=1$ are the proposition and the corollary of \cite[Appendix A.1.2]{Serre}.
The generalisations of claims (1) - (3) to arbitrary $g$ follow immediately by taking products and we leave them to
the reader.
For claim (4), we note that the valuation argument used in the proof of assertion (b) of the proposition in
\cite[Appendix A.1.2]{Serre}
still works, when being replaced by the valuation matrix $\ord_p(V)$ associated to a period matrix $V$ for
$\Lambda$,  which we introduced at the beginning of this section.
\qed\medskip

\begin{Remark}\label{rem: gamma}
 The $\ell$-adic completion is explicitly given by
 $$
  \widehat{K^\times}^\ell \,\cong\,
  \left(\mu(K)/\mu_{\ell^\infty}(K)\right) \,\times\, \pi_K^{\ZZ_\ell} \,\times\,
  \left\{
  \begin{array}{cl}
    \{1\} & \mbox{ if } \ell\neq p\\
    U^{(1)}& \mbox{ if } \ell=p\,.
  \end{array}
  \right.
 $$
 In this decomposition, the map $\gamma_\ell$ can be understood componentwise.
\end{Remark}

As a consequence, we now describe $G_K$-equivariant homomorphisms between
$\ell$-adic Tate modules of abeloid varieties.
In the one-dimensional case, this result is implicit in \cite[Appendix A.1.4]{Serre}.

\begin{Proposition}
 \label{prop: tate module}
 Let $\Lambda_A\subset K^{\times g}$ and $\Lambda_B\subset K^{\times h}$
 be lattices and let $A:=\GG_{m,K}^g/\Lambda_A$ and
 $B:=\GG_{m,K}^h/\Lambda_B$ be the associated abeloid varieties.
 Let $Q_A$ and $Q_B$ be period matrices for $\Lambda_A$ and $\Lambda_B$.
 Then, there exist isomorphisms of $\ZZ_\ell$-modules
 
 \begin{equation*}
 \begin{adjustbox}{width=12cm}{$
\begin{array}{l}
\Hom_{G_K}\left(T_\ell(A),T_\ell(B)\right)\\
\cong\left\{ M\in{\rm Mat}_{g\times h}(\ZZ_\ell) \,|\, \gamma_\ell(\Lambda_A)\odot M \,\subseteq\, \gamma_\ell(\Lambda_B)   \right\} \\
\cong\left\{ M\in{\rm Mat}_{g\times h}(\ZZ_\ell) \,|\, \exists N\in{\rm Mat}_{h\times g}(\ZZ_\ell)\,:\,    \gamma_\ell(Q_A)\odot M \,=\, N\odot \gamma_\ell(Q_B)   \right\}  \\  
\cong\left\{ M\in{\rm Mat}_{g\times h}(\ZZ_\ell) \,|\,      \left(\ord_p(Q_A)^{-1}\odot \gamma_\ell(Q_A)\right) \odot M      \,=\, M\odot \left(\ord_p(Q_B)^{-1}\odot \gamma_\ell(Q_B)\right)
   \right\} \,.
\end{array} $}
\end{adjustbox}
\end{equation*}
Here, $\gamma_\ell:K^\times\to \widehat{K^\times}^\ell$ denotes the $\ell$-adic completion from Lemma \ref{lemma: Serre}.
\end{Proposition}

\prf
Given a $G_K$-equivariant morphism $\varphi:T_\ell(A)\to T_\ell(B)$, we obtain a commutative
diagram
\begin{center}
\begin{tikzpicture}[descr/.style={fill=white,inner sep=1.5pt}]
        \matrix (m) [
            matrix of math nodes,
            row sep=2.5em,
            column sep=2.5em,
            text height=1.5ex, text depth=0.25ex
        ]
        { 0 & \ZZ_\ell(1)^g & T_\ell(A) & \ZZ_\ell^g & 0 \\
        0 & \ZZ_\ell(1)^h & T_\ell(B) & \ZZ_\ell^h & 0 \\
        };

       \path[overlay,->, font=\scriptsize]
       (m-1-1) edge (m-1-2)
       (m-1-2) edge (m-1-3)
       (m-1-3) edge (m-1-4)
       (m-1-4) edge (m-1-5)
       (m-2-1) edge (m-2-2)
       (m-2-2) edge (m-2-3)
       (m-2-3) edge (m-2-4)
       (m-2-4) edge (m-2-5)
       (m-1-2) edge node [right]{$\rho$} (m-2-2)
       (m-1-3) edge node [right]{$\varphi$} (m-2-3)
       (m-1-4) edge node [right]{$\sigma$} (m-2-4);
              
        \end{tikzpicture} 
\end{center}
for some matrices $\rho,\sigma\in{\rm Mat}_{g\times h}(\ZZ_\ell)$.
Taking $G_K$-invariants and passing to cohomology,
it follows that the diagram
\begin{equation}\label{eq: extension diagram}
\begin{tikzpicture}[descr/.style={fill=white,inner sep=1.5pt}]
        \matrix (m) [
            matrix of math nodes,
            row sep=2.5em,
            column sep=2.5em,
            text height=1.5ex, text depth=0.25ex
        ]
        { \ZZ_\ell^g &  H^1\left(G_K,\ZZ_\ell(1)^g\right) \\
        \ZZ_\ell^h & H^1\left(G_K,\ZZ_\ell(1)^h\right) \\
        };

       \path[overlay,->, font=\scriptsize]
       (m-1-1) edge node [above]{$d_{g}$} (m-1-2)
       (m-2-1) edge node [above]{$d_{h}$} (m-2-2)
       (m-1-1) edge node [right]{$\sigma$} (m-2-1)
       (m-1-2) edge node [right]{$\rho_{\ast}$} (m-2-2);
              
        \end{tikzpicture} 
\end{equation}
commutes.

By Lemma \ref{lemma: Serre}, the images of $\Lambda'_A$ and $\Lambda_B'$ under
$d_g$ and $d_h$ are lattices, and thus, the homomorphisms $d_g$ and $d_h$
are injective.
In particular, $\rho$ determines $\rho_\ast$, which determines $\sigma$ uniquely 
and vice versa.

Using the results and identifications of Lemma \ref{lemma: Serre}, the commutativity of
the above diagram implies
$$
 \gamma_\ell(Q_A)\odot\rho \,=\, \sigma\odot \gamma_\ell(Q_B)
$$
with respect to the notation introduced above.
Thus, $\gamma_\ell(\Lambda_A)\odot \rho$ lies in the $\ZZ_\ell$-span of 
$\gamma_\ell(\Lambda_B)$.

Conversely, let $M\in{\rm Mat}_{g\times h}(\ZZ_\ell)$ be such that
$\gamma_\ell(\Lambda_A)\odot M$ is contained in the $\ZZ_\ell$-span of 
$\gamma_\ell(\Lambda_B)$.
Then, $M$ gives rise to map $\ZZ_\ell(1)^g\to \ZZ_\ell(1)^h$ and
we can find a unique matrix $N\in{\rm Mat}_{g\times h}(\ZZ_\ell)$
defining a map $\ZZ_\ell^g\to\ZZ_\ell^h$
such that the diagram \eqref{eq: extension diagram} commutes.
This commutativity implies that $M$ and $N$ determine a unique $G_K$-equivariant 
map $\varphi:T_\ell(A)\to T_\ell(B)$.

The last isomorphism follows from taking valuations as in the proof
of Theorem \ref{thm: abeloid morphisms}.
\qed\medskip

\subsection{The $p$-adic Galois representations}\label{The p-adic Galois representations}
Let $A=\GG_{m,K}^g/\Lambda_A$ be an abeloid variety over a $p$-adic field $K$ and 
let $Q_A=(q_{i,j})$ be a period matrix for $\Lambda_A$.
As seen in \eqref{eq: abeloid extension}, the $p$-adic Galois representation of $G_K$ on the rational 
Tate module $V_p(A)$ is an extension of 
$\QQ_p(1)^g$ by $\QQ_p^g$.
We denote by $\log_p$ Iwasawa's $p$-adic logarithm, normalised such that
$\log_p(p)=0$. 
Associated to this data, we construct a filtered $(\varphi,N)$-module over $K$
as follows:
\begin{enumerate}
 \item Let $V$ be the $2g$-dimensional vector space 
   over $\QQ$ with basis $x_1,...,x_g,y_1,...,y_g$
   together with two linear operators $\varphi, N$:
   $$
   \begin{array}{lclclcl}
    \varphi(x_i) &=& p^{-1}\cdot  x_i & \mbox{ \qquad } & \varphi(y_i) &=& y_i \\
    N(x_i) &=& 0 & & N(y_i) &=& \sum_{j=1}^g \nu_p(q_{i,j})\cdot x_j
   \end{array}
   $$
   that is, these operators are given by matrices
   $$
    \begin{pmatrix}
      p^{-1}\cdot \mathrm{Id}_{g\times g} & 0 \\
      0 & \mathrm{Id}_{g\times g}
    \end{pmatrix}
    \mbox{ \quad and \quad }
    \begin{pmatrix}
      0 & {\rm ord}_p(Q) \\
      0 & 0
    \end{pmatrix}\,.
   $$
   We equip $V_{K_0}:=V\otimes_\QQ K_0$ with the $K_0$-linear extension
   $N\otimes{\rm id}_{K_0}$ of $N$ and with the $K_0$-semi-linear extension 
   $\varphi\otimes\sigma$ of $\varphi$.
   Here, $\sigma$ denotes the lift of Frobenius on $K_0$ and by abuse of notation,
   we will denote these extensions again by $\varphi$ and $N$.
   This turns $(V_{K_0},\varphi,N)$ into a $(\varphi,N)$-module.
   \item A filtration on $V_K:=V\otimes_\QQ K$ defined by 
   ${\rm Fil}^i=0$ for $i\geq1$, by ${\rm Fil}^i=V_K$ for $i<0$, and ${\rm Fil}^0$
   is the $g$-dimensional $K$-vector space spanned by
   $$
      y_i\,+\, \sum_{j=1}^g \log_p(q_{i,j})\cdot x_j
   $$
   for $i=1,...,g$.
\end{enumerate}

After these preparations, we obtain the following result, which was already known
for Tate elliptic curves, that is, in the case where $g=1$, see also 
Remark \ref{rem: LeStum} below.

\begin{Theorem}\label{thm: abeloid phi N}
  Let $A=\GG_{m,K}^g/\Lambda_A$ be an abeloid variety over a $p$-adic field $K$.
  Then, the filtered $(\varphi,N)$-module $\DD_{\rm st}(V_p(A))$ associated to the
  rational Tate module of $A$ is isomorphic to $(V_{K_0},\varphi,N,{\rm Fil}^\bullet)$
  constructed above.
\end{Theorem}

\prf
We use the notations of \cite[\S II.4]{Berger} and generalise the computations there
from $g=1$ to arbitrary $g$.
To obtain an explicit description of the rational Tate module $V_p(A)$,
we fix a compatible system $\{\varepsilon^{(n)}\}_n$ of $p^n$-th roots of unity,
as well as a compatible system $\{q_{i,j}^{(n)}\}_n$ of $p^n$-th roots of $q_{i,j}$.
Via $p$-adic uniformisation, we obtain a $G_K$-equivariant parametrisation
$\overline{K}^{\times g} / \Lambda_A \to A(\overline{K})$, which 
identifies the $p^n$-torsion subgroup of $A(\overline{K})$ with 
$$
  \{ x\in \overline{K}^{\times g} \,|\, x^{p^n}\in\Lambda_A \},
$$
which shows that the $2g$ vectors
$$
 \begin{array}{lcl}
  e_i  &:=& \varprojlim_n (1,...,\varepsilon^{(n)},1,...) \\
  f_i  &:=& \varprojlim_n (q_{i,1}^{(n)},...,q_{i,g}^{(n)})
 \end{array}
$$ 
with $i=1,...,g$ form a $\ZZ_p$-basis of the Tate module $T_p(A)$.
Thus, if $g\in G_K$, we compute
$$
 g(e_i) \,=\, \varprojlim_n (1,...,g(\varepsilon^{(n)}),1,...)
  \,=\, \chi(g)\cdot \varprojlim_n (1,...,\varepsilon^{(n)},1,...)
  \,=\, \chi(g)\cdot e_i,
$$
where $\chi:G_K\to\ZZ_p^\times$ denotes the $p$-adic cyclotomic character.
Moreover, we define a $g\times g$ matrix $C(g)=\left(c_{i,j}(g)\right)$
with entries in $\ZZ_p$ via 
\begin{eqnarray*}
 g(f_i)  &=& \varprojlim_n \left(g(q_{i,1}^{(n)}),...,g(q_{i,g}^{(n)})\right)  \\
  &=& \varprojlim_n \left(q_{i,1}^{(n)} (\varepsilon^{(n)})^{c_{i,j}(g)} ,...,q_{i,g}^{(n)} (\varepsilon^{(n)})^{c_{i,j}(g)} \right) \\
  &=& f_i + \sum_{j=1}^g c_{i,j}(g)\cdot e_j\,.
\end{eqnarray*}
Thus, the action of $g\in G_K$ on $T_p(A)$ is given
by the matrix
$$
    \begin{pmatrix}
      \chi(g)\cdot \mathrm{Id}_{g\times g} & C(g) \\
      0 & \mathrm{Id}_{g\times g}
    \end{pmatrix}\,.  
$$
To determine the $p$-adic periods, we have
$t=\log_p([\varepsilon])\in B_{\rm dR}^+\subset B_{\rm dR}$
and set
$$
   u_{i,j} \,:=\, \log_p(q_{i,j}) \,-\,\sum_{n=1}^\infty \frac{ (1-[\widetilde{q}_{i,j}])^n }{n}.
$$
were 
$$
  \widetilde{q}_{i,j} \,=\, (q_{i,j}^{(0)},q_{i,j}^{(1)},\ldots) \,\in\, \widetilde{\mathbb{E}}^{+} \,=\, 
  \displaystyle\varprojlim_{x\mapsto x^{p}}\mathcal{O}_{\widehat{\overline{K}}}
$$ 
and $[\widetilde{q}_{i,j}]\in W(\widetilde{\mathbb{E}}^{+})$ denotes its Teichm\"{u}ller lift. 
This series converges in $B_{\rm dR}^+$ and as explained in \cite[\S II.4.3]{Berger},
one should think of $u_{i,j}$ as being equal to $\log_p([\widetilde{q}_{i,j}])$ and one has
$g(u_{i,j})=u_{i,j}+c_{i,j}(g)t$.
In this explicit description, it is easy to see that the $2g$ vectors
$$
 \begin{array}{lcl}
  x_i  &:=& \frac{1}{t}\otimes e_i \\
  y_i  &:=& -\sum_{j=1}^g \frac{u_{i,j}}{t}\otimes e_j \,+\, 1\otimes f_i
 \end{array}
$$ 
with $i=1,...,g$ lie in $\DD_{dR}(V_p(A))=(B_{\rm dR}\otimes_{\QQ_p} V_p(A) )^{G_K}$.
These elements form a $2g$-dimensional vector space, which shows explicitly that
the $G_K$-representation on $V_p(A)$ is de~Rham.
Now, $t$ and the $u_{i,j}$ lie in the subring $B^+_{\rm cris}\subset B^+_{\rm dR}$ and we have
$\varphi(t)=p\cdot t$, $\varphi(u_{i,j})=p\cdot u_{i,j}$, $\varphi(e_i)=e_i$, and $\varphi(f_i)=f_i$,
see also \cite[\S II.4.3]{Berger}.
This implies 
$$
 \varphi(x_i) \,=\, \frac{1}{pt}\otimes e_i \,=\, p^{-1}\cdot x_i
 \mbox{ \quad and \quad }
  \varphi(y_i) \,=\, -\sum_{j=1}^g \frac{p\cdot u_{i,j}}{p\cdot t}\otimes e_j \,+\, 1\otimes f_i\,=\,y_i\,.
$$
Then, we have $B_{\rm st}=B_{\rm cris}[Y]$ and the normalisation $\log_p(p)=0$ determines
an embedding of $B_{\rm st}$ into $B_{\rm dR}$, see \cite[\S II.3.3]{Berger}.
Then,
$$
  u_{i,j} \,=\, \log_p [\widetilde{q}_{i,j}] \,=\, \nu_p(q_{i,j})\cdot Y \,+\, \log_p \left[ \frac{\widetilde{q}_{i,j}}{\nu_p(q_{i,j})} \right].
$$
The monodromy operator is given by $N:=-\frac{d}{dY}$ and we compute $N(x_i)=0$.
Rewriting $y_i$ as
$$
  y_i \,=\, -\sum_{j=1}^g \frac{\nu_p(q_{i,j})}{t}\cdot Y\otimes e_j 
   \,-\, \sum_{j=1}^g \frac{1}{t}\cdot \log_p\left[\frac{\widetilde{q}_{i,j}}{\nu_p(q_{i,j})}\right]\otimes e_j
   \,+\, 1\otimes f_i\,,
$$
we compute 
$$
 N(y_i) \,=\, -\frac{d}{dY}\, y_i 
  \,=\, \sum_{j=1}^g \frac{\nu_p(q_{i,j})}{t}\otimes e_j  
  \,=\, \sum_{j=1}^g \nu_p(q_{i,j}) e_j.
$$
By definition, the filtration on $B_{\rm dR}$ is defined by ${\rm Fil}^i(B_{\rm dR})=t^i\cdot B^+_{\rm dR}$, from 
which it is easy to see that the induced filtration ${\rm Fil}^i$ on the $K$-span of the $x_i,y_i$, that is, 
the intersection with $t^i\cdot B_{\rm dR}^+$, is zero for $i\geq1$ and it is equal to the whole space for $i<0$.
Moreover, the elements 
$$
  y_i\,+\, \sum_{j=1}^g \log_p(q_{i,j})\cdot x_j,\mbox{ \quad }i=1,...,g
$$
lie in ${\rm Fil}^0$, see also \cite[\S II.4.3]{Berger} 
and we leave it to the reader to show that these $g$ vectors actually span ${\rm Fil}^0$.
This establishes the claim and we see from these explicit computations that
the $G_K$-representation on $V_p(A)$ is semi-stable. 
We remark that the semi-stability of $V_p(A)$ is a special case 
of \cite[Corollary 5.26]{CN17}.
\qed\medskip

\begin{Remark}\label{rem: LeStum}
 For Tate elliptic curves, this result was established by Le~Stum \cite[\S 9]{LeStum} and
 our computations extended the exposition in \cite[\S II.4]{Berger}.
 For the description of $\varphi$ and $N$, we also refer to \cite{Co00, CI99}.
 Since we followed the exposition for Tate elliptic curves in \cite[\S II.4]{Berger}, we have chosen to 
 use the notation found therein. 
 We note that $\widetilde{\mathbb{E}}^{+}\cong\displaystyle\varprojlim_{x\mapsto x^{p}}\mathcal{O}_{\widehat{\overline{K}}}/p=\mathcal{O}_{\widehat{\overline{K}}}^{\flat}$
 and also that $W(\widetilde{\mathbb{E}}^{+})$ is commonly referred to as $\mathbb{A}_{\mathrm{inf}}$. 
\end{Remark}

We end this discusssion by an analogue of Proposition \ref{prop: tate module} in the context of filtered
$(\varphi,N)$-modules.
Before stating the result, we extend the definition of the $\MazurL$-invariant of a Tate elliptic curve 
\cite[Ch.~II \S1]{MTT86} to abeloid varieties of aribitrary dimension:
if $Q_A$ is a period matrix for the abeloid variety $A:=\GG_{m,K}^g/\Lambda_A$, we set
$$
 \MazurL(Q_A) \,:=\, \ord_p(Q_A)^{-1}\cdot\log_p(Q_A)\in\mathrm{Mat}_{g\times g}(K^{\times}).
$$
Here, $\log_p(Q_A)$ denotes the matrix obtained by applying $\log_p$
to every entry of $Q_A$.
Note that by a definition of period matrices, $\ord_p(Q_A)\in{\rm Mat}_{g\times g}(\QQ)$ 
is an invertible matrix. 

\begin{Proposition}\label{prop: abeloid morphisms via L}
 Let $A=\mathbb{G}_{m,K}^{g}/\Lambda_{A}$ and $B=\mathbb{G}_{m,K}^{h}/\Lambda_{B}$ be abeloid varieties 
 and let $Q_A$ and $Q_B$ be period matrices for $A$ and $B$, respectively.
\begin{enumerate}
\item If $A$ is isomorphic to $B$, then there exists a $M\in{\rm GL}_{g\times g}(\ZZ)$ such that
  \begin{equation*}
  \MazurL(Q_B) \,=\, M^{-1}\cdot \MazurL(Q_{A}) \cdot M \,.
  \end{equation*} 
  In particular, this describes how $\MazurL$
  transforms under a change of period matrix of one abeloid variety.
\item 
  There exists an isomorphism of $\QQ$-vector spaces 
   \begin{eqnarray*}
   & & \Hom(A,B)\otimes_\ZZ\QQ\\
   &\cong&
   \left\{ M\in{\rm Mat}_{g\times h}(\QQ) \,|\, 
   \MazurL(Q_A)\cdot M \,=\, M\cdot\MazurL(Q_B)
   \right\} \,.
 \end{eqnarray*}
 In particular, $A$ and $B$ are isogenous if and only if $g=h$ and there exists a $M\in{\rm GL}_g(\QQ)$
 such that $\MazurL(Q_A)\cdot M=M\cdot\MazurL(Q_B)$.
\end{enumerate}
\end{Proposition} 

\prf
We start with the following computation:
let $R$ be a subring of $\QQ_p$, 
let $A=\mathbb{G}_{m,K}^{g}/\Lambda_{A}$ and $B=\mathbb{G}_{m,K}^{h}/\Lambda_{B}$
be abeloid varieties, and let $Q_A$ and $Q_B$ be period matrices.  
Moreover, let $M\in{\rm Mat}_{g\times h}(R)$ be a matrix such that there exists
a $N\in{\rm Mat}_{g\times h}(R)$ with
\begin{equation}
 \label{eq: period matrix}
   Q_A\,\odot\, M  \,=\, N\,\odot\,Q_B.
\end{equation}
We have seen in the proof of Theorem \ref{thm: abeloid morphisms} that
we have $N=\ord_p(Q_A)\cdot M\cdot \ord_p(Q_B)^{-1}$ in this case, which then yields
$$
   (\ord_p(Q_A)^{-1}\odot Q_A) \odot M \,=\, M \odot (\ord_p(Q_B)^{-1}\odot Q_B).
$$
Taking the Iwasawa logarithm on both sides, we obtain
$$
  \MazurL(Q_A)\cdot M \,=\, M\cdot\MazurL(Q_B).
$$
This already implies claim (1): if $A\cong B$, then 
we can find $M,N\in{\rm GL}_g(\ZZ)$ satisfying 
equation \eqref{eq: period matrix} and the assertion follows.

We now establish claim (2): given 
a homomorphism $\Hom(A,B)\otimes\QQ$, there
exist by Theorem \ref{thm: abeloid morphisms} two 
matrices $M,N\in{\rm Mat}_{g\times h}(\QQ)$ that satisfy \eqref{eq: period matrix}.
By the above computations, we find that 
$\MazurL(Q_A)\cdot M=M\cdot\MazurL(Q_B)$.

Conversely, assume that we are given a matrix 
$M\in{\rm Mat}_{g\times h}(\QQ)$ such 
that $\MazurL(Q_A)\cdot M=M\cdot \MazurL(Q_B)$.
We set $N:=\ord_p(Q_A)\cdot M\cdot \ord_p(Q_B)^{-1}$ and find
$$
  \log_p(Q_A)\cdot M \,=\, N\cdot \log_p(Q_B) \mbox{ \qquad and \qquad }
  \ord_p(Q_A)\cdot M \,=\, N\cdot \ord_p(Q_B),
$$
which shows that 
$$
  \log_p(Q_A\odot M) \,=\, \log_p(N\odot Q_B) \mbox{ \qquad and \qquad }
  \ord_p(Q_A\odot M) \,=\, \ord_p(N\odot Q_B).
$$
Using properties of the Iwasawa logarithm,
this shows that there exist roots of unity $\varepsilon_{i,j}\in K$ 
such that the $(i,j)$ entries of the matrices
$Q_A\odot M$ and $N\odot Q_B$ differ by the factor $\varepsilon_{i,j}$,
see also the proof of \cite[Proposition 6]{LeStum}.
Thus, if $R$ is a positive integer such that $\varepsilon_{i,j}^R=1$ for all $i,j$, then
$Q_A\odot (RM)=(RN)\odot Q_B$.
In particular, $RM$ and $RN$ define an element of $\Hom(A,B)$ and
claim (2) follows.
\qed\medskip
 
Our definition of the $\MazurL$-invariant of an abeloid variety generalises the 
$\MazurL$-invariant of a Tate elliptic curve \cite[Ch.~II \S1]{MTT86}. 
We refer to \cite{Dasgupta} for a survey of various $\MazurL$-invariants for 
varieties that are uniformised by Drinfeld's upper half plane $\widehat{\Omega}^1_K$.
Before continuing with our discussion of Proposition \ref{prop: tate module} in the context of 
filtered $(\varphi,N)$-modules, we first show that the $\MazurL$-invariant $\MazurL^{\mathrm{Col}}(C)$ 
of a Mumford curve $C$, as defined by Besser-de~Shalit \cite{BdS16}, coincides with our 
$\MazurL$-invariant of the Jacobian $J$ of $C$.
 
\begin{Proposition}\label{prop: besser}
 Let $\Gamma\subset\mathrm{PGL}_{2}(K)$ be a Schottky group
 and let $C=\Gamma\backslash\widehat{\Omega}^1_{K}$ be the associated Mumford 
 curve of genus $g$. 
 Let $J=\GG_{m,K}^{g}/\Lambda_{J}$ be the Jacobian of $C$. 
 Then, there exists a choice of period matrix $Q_J$ for $J$ such that 
 \begin{equation*}
  \MazurL^{\mathrm{Col}}(C) \,=\,\MazurL(Q_J)\,.
 \end{equation*} 
\end{Proposition}  

\prf
(We refer the reader to \cite[\S2.1 and \S2.2]{BdS16} for a summary of the de Rham cohomology of varieties that are uniformised 
by Drinfeld upper half spaces, and also for the appropriate references to the literature.)

Since we only consider de Rham cohomology in cohomological degree $1$, we are in the situation $d=1$ and $i=0$ 
in the notation of \cite{BdS16}. 
We have already fixed a uniformiser of $K$, so we shall drop this label from the notation of \cite{BdS16}. 
In summary, we simplify the notation by setting $\MazurL^{\mathrm{Col}}(C):=\MazurL_{\pi,1}^{\mathrm{Col}}(C)$. 

By definition, we have $\MazurL^{\mathrm{Col}}(C)=\nu^{-1}\circ\lambda^{\mathrm{Col}}$ where 
\begin{equation*}
  \nu \,:=\, \mathrm{gr}N \,:\, \mathrm{gr}_{\Gamma}^{0}\HdR{1}(C/K) \,\to\, \mathrm{gr}_{\Gamma}^{1}\HdR{1}(C/K)
\end{equation*}
is the map induced by the monodromy operator on the graded pieces of the covering filtration, 
using the fact that the covering filtration coincides with the weight filtration up to a shift in index. 
It is an isomorphism by the monodromy-weight conjecture, which is a theorem in our situation. 
The covering filtration is opposite to the Hodge filtration, which gives the identifications
\begin{equation*}
\mathrm{gr}_{\Gamma}^{0}H_{\mathrm{dR}}^{1}(C/K)\cong H^{0}(C,\Omega_{C/K}^{1})
\end{equation*} 
and 
\begin{equation*}
\mathrm{gr}_{\Gamma}^{1}H_{\mathrm{dR}}^{1}(C/K)\cong H^{1}(C,\mathcal{O}_{C})\,.
\end{equation*} 
The monodromy operators of $C$ and $J$ coincide under the identification 
$\HdR{1}(C/K)=\HdR{1}(J/K)$, by \cite[\S3]{CI99}. 
Choose a period matrix $Q_{J}$ of $J$. 
Then we have computed in Theorem \ref{thm: abeloid phi N} that $\nu$ is given by the matrix $\ord_p(Q_{J})$. 

The map
\begin{equation*}
\lambda^{\mathrm{Col}}:\mathrm{gr}_{\Gamma}^{0}\HdR{1}(C/K)\rightarrow\mathrm{gr}_{\Gamma}^{1}\HdR{1}(C/K)
\end{equation*}
is defined using harmonic cochains on the Bruhat-Tits building $\mathcal{T}$ of $\mathrm{PGL}_{2}(K)$ 
and the identifications
\begin{equation*}
\mathrm{gr}_{\Gamma}^{0}\HdR{1}(C/K) \,\cong\, H^{0}(\Gamma,C_{\mathrm{har}}^{1}(\mathcal{T}))
\end{equation*} 
and 
\begin{equation*}
\mathrm{gr}_{\Gamma}^{1}\HdR{1}(C/K) \,\cong\, H^{1}(\Gamma,C_{\mathrm{har}}^{0}(\mathcal{T}))\,.
\end{equation*} 
Let $\CC_{p}:=\widehat{\overline{K}}$. 
It is shown in \cite[\S3.1.2]{Gro00} that for Mumford curves, the map $\lambda^{\mathrm{Col}}\otimes_{K}\CC_{p}$ 
coincides with Coleman integration. 

Fix $z'\in\widehat{\Omega}^1_{K}(\mathbb{C}_{p})$. 
For $\alpha\in\Gamma$, define 
\begin{equation*}
 u_{\alpha}(z) \,:=\, \displaystyle\prod_{\gamma\in\Gamma}\,\frac{z-\gamma(z')}{z-\gamma\alpha(z')}
\end{equation*}
(this is independent of the choice of $z'$), 
and for $\alpha,\beta\in\Gamma$, define
\begin{equation*}
 Q_{\alpha,\beta} \,:=\,\frac{u_{\alpha}(z)}{u_{\alpha}(\beta z)}\in K^{\times}\,.
\end{equation*}
Then, $\Lambda_{J}$ is generated by elements of the form $Q_{\alpha,\beta}$, that is the $Q_{\alpha,\beta}$ form a period matrix $Q_{J}$ for $J$ \cite[Ch.~VI \S2]{GvdP80}. 
Recall as well that $H^{0}(C,\Omega^{1}_{C/K})$ is generated over $K$ by differential forms of the form
\begin{equation*}
\omega_{\alpha} \,:=\,\frac{du_{\alpha}}{u_{\alpha}},
\end{equation*}
where $\alpha\in\Gamma':=\Gamma/[\Gamma,\Gamma]\simeq\Lambda_{J}$. 
Using the notation of \cite{Gro00}, all paths $\gamma$ are linear combinations of the paths of the 
form $[z,\beta\cdot z]$ with $\beta\in\Gamma$, and hence to compute $\lambda^{\mathrm{Col}}$ 
we see that it suffices to calculate
\begin{equation*}
 \int_{z}^{\beta\cdot z}\frac{du_{\alpha}}{u_{\alpha}} \,=\, \log_{p}u_{\alpha}(\beta\cdot z)-\log_{p}u_{\alpha}(z) \,=\, \log_{p}Q_{\alpha,\beta}\,.
\end{equation*}
Altogether, we see that $\MazurL^{\mathrm{Col}}=\nu^{-1}\circ\lambda^{\mathrm{Col}}=\mathrm{ord}_{p}(Q_{J})^{-1}\cdot\log_{p}(Q_{J})=\MazurL(Q_J)$.
\qed\medskip

We now return to our discussion of \ref{prop: tate module} in the context of filtered $(\varphi,N)$-modules. 
\begin{Proposition}
 \label{prop: tate module2}
 Let $A:=\GG_{m,K}^g/\Lambda_A$ and
 $B:=\GG_{m,K}^h/\Lambda_B$ be abeloid varieties and let 
 $Q_A$ and $Q_B$ be period matrices for $A$ and $B$, respectively.
 Then, there exists an isomorphism of $\QQ_p$-vector spaces
 \begin{eqnarray*}
   &&  
   \Hom_{\mathrm{MF}_{K}^{\mathrm{wa},\varphi,N}} \left(\DD_{\mathrm{st}}(V_{p}(A)),\DD_{\mathrm{st}}(V_{p}(B))\right) \\
   &\cong&
   \left\{ M\in{\rm Mat}_{g\times h}(\QQ_p) \,|\, \MazurL(Q_A)\cdot M \,=\, M\cdot \MazurL(Q_B)
   \right\} \,.
 \end{eqnarray*}
   
  In particular, $\DD_{\mathrm{st}}(V_{p}(A))$ and $\DD_{\mathrm{st}}(V_{p}(B))$ are isomorphic
  if and only if $g=h$ and there exists a $M\in{\rm GL}_g(\QQ_p)$ such that
  $\MazurL(Q_A)\cdot M = M\cdot \MazurL(Q_B)$.
\end{Proposition}

\prf
Quite generally, let $D,D'$ be objects of $\mathrm{MF}_{K}^{\mathrm{wa},\varphi,N}$. 
Then, an element of $\Hom_{\mathrm{MF}_{K}^{\mathrm{wa},\varphi,N}}(D,D')$ is a $K_{0}$-linear map from the underlying 
$K_{0}$-vector space of $D$ to the underlying $K_{0}$-vector space of $D'$ that commutes 
with the Frobenius and monodromy operators, and such that the $K$-linear extension of this map sends 
the Hodge filtration of $D_{K}$ into the Hodge filtration of $D'_{K}$. 
After choosing bases, we represent this $K_{0}$-linear map by a matrix $M$ and let $F,F'$ be matrices representing 
the Frobenius operators on $D$ and $D'$, respectively. 
Then, compatibility with Frobenius is the condition that $F M=M F'$.

We claim that the matrix $M$ has coefficients in $\QQ_p$:
first, by considering 
$$
 \left(\begin{smallmatrix}
\mathrm{Id} & 0 \\
M & \mathrm{Id}
\end{smallmatrix}\right) \,:\,D\oplus D' \,\to\, D\oplus D' ,
$$ 
we see that it suffices to assume $D=D'$, $\varphi=\varphi'$, and $F=F'$.
Recall that the Frobenius operator is semi-linear with respect to $\sigma$ by definition, that is,
$F M=\sigma(M) F$. 
This together with condition $F M=M F$ shows that $\sigma(M)=M$, which proves the claim.

\sloppy Now let us return to the explicit description of $\Hom_{\mathrm{MF}_{K}^{\mathrm{wa},\varphi,N}}(\DD_{\mathrm{st}}(V_{p}(A)),\DD_{\mathrm{st}}(V_{p}(B)))$. 
Using the bases from Section \ref{The p-adic Galois representations}, we see that giving an element of 
this space is equivalent to giving a matrix $M\in {\rm Mat}_{2g\times 2h}(\QQ_{p})$ that satisfies the following
three conditions
\begin{equation*}
\tag{a}\label{a}
\begin{pmatrix}
p^{-1}\cdot\mathrm{Id}_{g\times g} & 0 \\
0 & \mathrm{Id}_{g\times g}
\end{pmatrix}
M\,=\,M
\begin{pmatrix}
p^{-1}\cdot\mathrm{Id}_{h\times h} & 0 \\
0 & \mathrm{Id}_{h\times h}
\end{pmatrix}
\end{equation*}
\begin{equation*}
\tag{b}\label{b}
\begin{pmatrix}
0 & \mathrm{ord}_{p}(Q_{A}) \\
0 & 0
\end{pmatrix}
M\,=\,M
\begin{pmatrix}
0 & \mathrm{ord}_{p}(Q_{B}) \\
0 & 0
\end{pmatrix}
\end{equation*}
\begin{equation*}
\tag{c}\label{c}
\begin{pmatrix}
\log_{p}(Q_{A}) & \mathrm{Id}_{g\times g}
\end{pmatrix}
M \,=\, 
\begin{pmatrix}
Z\cdot\log_{p}(Q_{B}) & Z
\end{pmatrix}
\end{equation*}
for some $Z\in {\rm Mat}_{g\times h}(K)$.
We see that \eqref{a} holds if and only if $M$ is of the form 
$\left(\begin{smallmatrix}
X & 0 \\
0 & Y
\end{smallmatrix}\right)$
for some matrices $X,Y\in {\rm Mat}_{g\times h}(\QQ_{p})$. 
It then follows from \eqref{b} and \eqref{c} that to give an element of 
$\Hom_{\mathrm{MF}_{K}^{\mathrm{wa},\varphi,N}}(\mathbb{D}_{\mathrm{st}}(V_{p}(A)),\mathbb{D}_{\mathrm{st}}(V_{p}(B)))$ 
is the same as giving $X,Y\in {\rm Mat}_{g\times h}(\mathbb{Q}_{p})$ that satisfy
\begin{equation*}\tag{d}\label{d}
\mathrm{ord}_{p}(Q_{A})\cdot X\,=\,Y\cdot\mathrm{ord}_{p}(Q_{B})
\end{equation*}
and
\begin{equation*}\tag{e}\label{e}
\log_{p}(Q_{A})\cdot X\,=\, Y\cdot\log_{p}(Q_{B}) \,.
\end{equation*}
Since $Q_B$ is a period matrix, $\ord_p(Q_B)$ is invertible, and we find
$Y=\ord_p(Q_A)\cdot X\cdot \ord_p(Q_B)^{-1}$.
Plugging this into \eqref{e}, the proposition follows.
\qed\medskip

\subsection{A translation of Question \ref{question: tate} into linear algebra}\label{subsec: reformulation question}
Let $A=\mathbb{G}_{m,K}^{g}/\Lambda_{A}$ and $B=\mathbb{G}_{m,K}^{h}/\Lambda_{B}$ 
be abeloid varieties over a $p$-adic field $K$ and let $\ell$ be a prime, possibly equal to $p$. 
We choose period matrices $Q_A$ and $Q_B$.

\begin{enumerate}
\item Under the identifications of Theorem \ref{thm: abeloid morphisms} and Proposition \ref{prop: tate module}
the homomorphism \eqref{tatehom2}
\begin{equation*}
  \mathrm{Hom}(A,B)\otimes_{\ZZ}\ZZ_{\ell} \,\to\, \Hom_{G_K}\left(T_{\ell}(A),\,T_{\ell}(B)\right)
\end{equation*}
is given by the homomorphism
\begin{eqnarray*}
  & & \{M\in\mathrm{Mat}_{g\times h}(\ZZ)\,|\,\Lambda_{A}\odot M\subset\Lambda_{B}\}\otimes_{\ZZ}\ZZ_{\ell} \\
  &\to& \{M\in\mathrm{Mat}_{g\times h}(\ZZ_{\ell})\,|\,\gamma_{\ell}(\Lambda_{A})\odot M\subset\gamma_{\ell}(\Lambda_{B}) \},
\end{eqnarray*}
or, equivalently,
\begin{equation*}
 \begin{adjustbox}{width=11.1cm}{$
\begin{array}{l}
\left\{ M\in{\rm Mat}_{g\times h}(\ZZ) \,|\, 
    \left(\ord_p(Q_A)^{-1}\odot Q_A\right) \odot M \,=\, M\odot \left(\ord_p(Q_B)^{-1}\odot Q_B\right)
   \right\} \otimes_\ZZ\ZZ_\ell\\
  \to \left\{ M\in{\rm Mat}_{g\times h}(\ZZ_\ell) \,|\, 
    \left(\ord_p(Q_A)^{-1}\odot \gamma_\ell(Q_A)\right) \odot M \,=\, M\odot \left(\ord_p(Q_B)^{-1}\odot \gamma_\ell(Q_B)\right)
   \right\}\,.
\end{array} $}
\end{adjustbox}
\end{equation*}
\item Moreover, if $\ell=p$, then under the identifications of Theorem \ref{thm: abeloid morphisms} 
and Proposition \ref{prop: tate module2}, the analog of homomorphism \eqref{tatehom2}
\begin{equation*}
\Hom(A,B)\otimes_{\ZZ}\QQ_{p} \,\to\,
\Hom_{\mathrm{MF}_{K}^{\mathrm{wa},\varphi,N}}\left(\DD_{\mathrm{st}}(V_{p}(A)),\DD_{\mathrm{st}}(V_{p}(B))\right)
\end{equation*}
is given by 
\begin{equation*}
 \begin{adjustbox}{width=11.1cm}{$
\begin{array}{l}
\left\{ M\in{\rm Mat}_{g\times h}(\QQ) \,|\, 
    \left(\ord_p(Q_A)^{-1}\odot Q_A\right) \odot M \,=\, M\odot \left(\ord_p(Q_B)^{-1}\odot Q_B\right)
   \right\} \otimes_\QQ\QQ_p\\
  \to \left\{ M\in{\rm Mat}_{g\times h}(\QQ_p) \,|\, 
    \MazurL(Q_A) \cdot M \,=\, M\cdot \MazurL(Q_B)
   \right\}\,,
\end{array} $}
\end{adjustbox}
\end{equation*}
or, equivalently, 
\begin{eqnarray*}
  & &   \left\{ M\in{\rm Mat}_{g\times h}(\QQ) \,|\, 
    \MazurL(Q_A) \cdot M \,=\, M\cdot \MazurL(Q_B)
   \right\} \otimes_\QQ\QQ_p \\
  &\to& \left\{ M\in{\rm Mat}_{g\times h}(\QQ_p) \,|\, 
    \MazurL(Q_A) \cdot M \,=\, M\cdot \MazurL(Q_B)
   \right\}\,.
\end{eqnarray*}
\end{enumerate}

This gives reformulations of Question \ref{question: tate} in terms of linear algebra.
Note that it is particularly easy to see that \eqref{tatehom2} is injective in some of 
these reformulations.
We also see that the surjectivity of \eqref{tatehom2} is equivalent to the
interplay of $\QQ$-structures versus $\QQ_\ell$-structures, which is similar
to Raskind-admissibility and the results of Section \ref{sec: translation}.
In Lemma \ref{lemma: Serre} and Remark \ref{rem: gamma}, 
we have seen that $\gamma_\ell$ behaves very differently
depending on whether $\ell\neq p$ or $\ell=p$:
\begin{enumerate}
 \item  In Proposition \ref{prop: totally deg ell}, we will see that surjectivity of \eqref{tatehom2} 
   may fail if $A$ and $B$ are Tate elliptic curves and $\ell\neq p$.
 \item Surjectivity of \eqref{tatehom2} may look plausible if $\ell=p$: 
  we will give a positive result in Proposition \ref{prop: isogeny product tate} below and disprove it in general in
  Theorem \ref{thm: counterexample}.
\end{enumerate}

\section{Products of Tate elliptic curves}
\label{sec: product tate}

In this section, we use the results of the previous section to study
Raskind's conjecture for divisors (Conjecture \ref{conj: raskind}) 
and Question \ref{question: tate} for abelian varieties that are isogenous to
products of Tate elliptic curves.
For the product $X$ of two Tate elliptic curves, we determine 
the rational structure on the filtered $(\varphi,N)$-module 
$\DD_{\rm st}(\Het{2}(X_{\overline{K}},\QQ_p))$,
which leads to a direct verification of Conjecture \ref{conj: raskind}.
We also classify all filtered $(\varphi,N)$-modules 
that are ordinary in the sense of Perrin-Riou that have 
a fixed rational structure that is modeled on a surface 
with $p_g=1$.

\subsection{Products of Tate elliptic curves}

We recall that a Tate elliptic curve over a $p$-adic field $K$ is the same as 
an abeloid variety of dimension one over $K$.
In particular, they are of the form $E(q):=\GG_{m,K}/q^\ZZ$ for some $q\in K^\times$ with $\nu_p(q)>0$.
As in the previous section, we set $\MazurL(x):=\log_p(x)/\nu_p(x)$.

\begin{Theorem}[Le~Stum, Serre]\label{thm: LeStum Serre}
  Let $E(q_i)$, $i=1,2$ be two Tate elliptic curves over a $p$-adic field $K$
  associated to $q_i\in K^\times$ with $\nu_p(q_i)>0$.
  Then, the following are equivalent:
  \begin{enumerate}
  \item $E(q_1)$ and $E(q_2)$ are isogenous.
  \item There exist positive integers $A_i$, $i=1,2$ such that $q_1^{A_1}=q_2^{A_2}$.
  \item The rational Tate modules $V_p(E(q_i))$, $i=1,2$ are isomorphic as $p$-adic $G_K$-representations.
  \item The $\DD_{\rm st}(V(E(q_i)))$, $i=1,2$ are isomorphic as filtered $(\varphi,N)$-modules over $K$.
  \item $\MazurL(q_1)=\MazurL(q_2)$.
  \end{enumerate}
\end{Theorem}

\prf
The equivalences $(1)\Leftrightarrow(2)\Leftrightarrow(3)$ are shown in \cite[\S A.1.4]{Serre}
and the equivalences $(1)\Leftrightarrow(4)\Leftrightarrow(5)$ are shown in \cite[Proposition 6]{LeStum}.
\qed\medskip

\begin{Remark}
 In fact, using the results from Section \ref{sec: abelian varieties}, it is easy to deduce
 the equivalence $(1)\Leftrightarrow(2)$  from Theorem \ref{thm: abeloid morphisms}, 
 to deduce the equivalence $(2)\Leftrightarrow(3)$ from Proposition \ref{prop: tate module}, 
 to prove the equivalence $(1)\Leftrightarrow(4)$ using Proposition \ref{prop: abeloid morphisms via L},
 and to prove the equivalences
 $(3)\Leftrightarrow(4)\Leftrightarrow(5)$ using Theorem \ref{thm: abeloid phi N}.
\end{Remark}

Concerning Question \ref{question: tate}, we have a positive answer in the following special case.

\begin{Proposition}\label{prop: isogeny product tate}
 Let $K$ be a $p$-adic field and let $A$ and $B$ be abelian varieties over $K$, both of which are
 isogenous to products of Tate elliptic curves. 
 Then, \eqref{tatehom2} (please see \S\ref{abelian and abeloid}) is surjective for $A$, $B$, and $\ell=p$, that is, Question \ref{question: tate} 
 has a positive answer for  $A$ and $B$.
\end{Proposition}

\prf
First, the map \eqref{tatehom2} is injective \cite[IV.19.3]{MumfordBook} and since the cokernel is 
torsion-free \cite[Lemma 1]{TateAbelian}, we are reduced to showing the surjectivity 
of the map
\begin{equation}\label{rational hom}
 \Hom(A,B) \otimes_{\ZZ} \QQ_{p} \,\longrightarrow\,  
 \Hom_{G_{K}}\left(V_p(A), V_p(B)\right) \,.
\end{equation}

If $A$ is isogenous to a product $\prod_{i=1}E(q_i)$ of Tate elliptic curves, then an isogeny gives rise to an 
isomorphism of $p$-adic $G_K$-representations
$$
   V_p(A) \,\cong\, \bigoplus_{i=1}^g\, V_p(E(q_i)).
$$

We are thus reduced to the case where both $A$ and $B$ are Tate elliptic curves.
If they are not isogenous, then ${\rm Hom}(A,B)=0$ and ${\rm Hom}(V_p(A),V_p(B))=0$ by
Theorem \ref{thm: LeStum Serre}.
On the other hand, if they are isogenous, then we have
${\rm Hom}(A,B)=\ZZ$ and ${\rm Hom}(V_p(A),V_p(B))=\QQ_p$.
This verifies the claim in both cases and the proposition follows.
\qed\medskip

As a corollary, we establish Conjecture \ref{conj: raskind} in a special case.

\begin{Corollary}\label{cor: product tate}
  Let $K$ be a $p$-adic field and let $A$ be an abelian variety over $K$ that is
  isogenous to a product of Tate elliptic curves. 
  Then, \eqref{tatehom} is surjective for $A$ and $\ell=p$, that is, 
  Raskind's conjecture for divisors (Conjecture \ref{conj: raskind}) is true for
  $A$.
\end{Corollary}

\prf
We may assume that $A$ is in fact isomorphic to a product of Tate elliptic curves.
The argument to deduce surjectivity of \eqref{tatehom} 
from the surjectivity of \eqref{tatehom2} is \cite[Theorem 3]{TateAbelian}.
\qed\medskip

In Appendix \ref{appendix}, we will see that both results fail to be true if $\ell\neq p$.

\subsection{The product of two Tate curves}\label{subsec: product tate}
Let $E(q):=\GG_{m,K}/q^\ZZ$ be the Tate elliptic curve associated to an element
$q\in K^\times$ with $\nu_p(q)>0$.
Associated to $q$, we construct a filtered $(\varphi,N)$-module as follows
\begin{enumerate}
 \item Let $V$ be a $2$-dimensional vector space 
   over $\QQ$ with basis $e_1,e_2$
   together with two linear operators $\varphi, N$ defined by
   $$
   \begin{array}{lclclclc}
    \varphi(e_1) &=& e_1 & \mbox{ \qquad } & \varphi(e_2) &=& p\cdot e_2 \\
    N(e_1) &=& 0 & & N(e_2) &=& e_1&.
   \end{array}
   $$
   We equip $V_{K_0}:=V\otimes_\QQ K_0$ with the $K_0$-linear extension
   of $N$ and the $K_0$-semi-linear extension of $\varphi$,
   which defines a $(\varphi,N)$-module over $K_0$.
   \item Define a filtration on $V_K:=V\otimes_\QQ K$ defined by 
   ${\rm Fil}^i=V_K$ for $i\geq0$, by ${\rm Fil}^i=0$ for $i\geq2$, and ${\rm Fil}^1$
   is the one-dimensional $K$-vector space spanned by
   $$
      \MazurL(q)\cdot e_1\,+\,e_2.
   $$
\end{enumerate}
Now, the $p$-adic $G_K$-representations $V_p(E(q))$ and
$\Het{1}(E(q)_{\overline{K}},\QQ_p)$ are dual.
Using the explicit description of $\DD_{\rm st}(V_p(E(q)))$ 
provided by Theorem \ref{thm: abeloid phi N}
(although in the case of Tate elliptic curves, this was classically known, see \cite[II.4.2]{Berger}), 
the formulae from \cite[4.2.4 and 4.3.4]{Fontaine} to compute the dual, 
and after suitably rescaling,
it is not difficult to see that the just constructed filtered $(\varphi,N)$-module $V$
is isomorphic to $\DD_{\rm st}(\Het{1}(E(q)_{\overline{K}},\QQ_p))$.
Alternatively, one can also use Le~Stum's computations  \cite[\S 9]{LeStum}.

Next, let $q_i\in K^\times$ with $\nu_p(q_i)>0$ for $i=1,2$, let $E(q_i)$ be the
assocated Tate elliptic curves over $K$, and set $X:=E(q_1)\times E(q_2)$.
Then, we have the following description of the filtered $(\varphi,N)$-module
$\DD_{\rm st}(\Het{2}(X_{\overline{K}},\QQ_p))$
\begin{enumerate}
\item Let $V^{(i)}$, $i=1,2$ be two $2$-dimensional vector spaces with bases
$\{e_1^{(i)},e_2^{(i)}\}$, Frobenius, monodromy, and filtration on $V^{(i)}\otimes K$
associated to $q_i$ as above.
Set $V:=V^{(1)}\oplus V^{(2)}$.
\item  Then, $W:=\wedge^2(V)=\wedge^2(V^{(1)}\oplus V^{(2)})$ is a $6$-dimensional 
$\QQ$-vector  space with basis
$$
\begin{array}{lll}
a\,=\, e_1^{(1)}\wedge e_1^{(2)},\\
b_0\,=\,e_1^{(1)}\wedge e_2^{(2)}+e_2^{(1)}\wedge e_1^{(2)},\\
b_1\,=\,e_1^{(1)}\wedge e_2^{(2)}-e_2^{(1)}\wedge e_1^{(2)},\quad&
b_2\,=\, e_1^{(1)}\wedge e_2^{(1)},\quad&
b_3\,=\, e_1^{(2)}\wedge e_2^{(2)}\\
c\,=\, e_2^{(1)}\wedge e_2^{(2)}.
\end{array}
$$
We set 
$$
 A:=\langle a\rangle,\, B_0:=\langle b_0\rangle,\, B_1:=\langle b_1,b_2,b_3\rangle,\, B:=B_0\oplus B_1,
 \mbox{ and } 
 C:=\langle c\rangle.
$$
\item The $\varphi$'s on $V^{(1)}$ and $V^{(2)}$ induce a linear endomorphism $\varphi_W$
with
$$
   \varphi_W|_A \,=\, {\rm id}_A, \mbox{ \quad }
   \varphi_W|_B \,=\, p\cdot{\rm id}_A, \mbox{ \quad and \quad }
   \varphi_W|_C \,=\, p^2\cdot {\rm id}_A\,.
$$
Similarly, we obtain a linear endomorphism $N_W$ with
$$
  N_W(c)\,=\,b_0,\mbox{ \quad }N_W(b_0)\,=\,2a,\mbox{ \quad and \quad }N_W|_{A\oplus B_1}\,=\,0.
$$
We extend $\varphi_W$ semi-linearly and $N_W$ linearly to $W\otimes K_0$
and thus obtain a $(\varphi,N)$-module.
\item  Moreover, the decomposition $(W=A\oplus B_0\oplus B_1\oplus C,\varphi_W,N_W)$ 
defines a rational structure in the sense of Definition \ref{def: rational structure}.
\item
The filtrations on $V^{(i)}\otimes K$ give rise to a filtration
on $W\otimes K$
with $\Fil^0=W$, $\Fil^3=0$, and $\Fil^2$ is the one-dimensional
$K$-span of
\begin{equation} 
 \label{eq: Mazur vector}
\underbrace{\MazurL(q_1)\MazurL(q_2)\cdot a}_{\in A} \,+\,
\underbrace{\MazurL(q_1)\cdot e_1^{(1)}\wedge e_2^{(2)} \,+\, \MazurL(q_2)\cdot e_2^{(1)}\wedge e_1^{(2)}}_{\in B}
\,+\,\underbrace{e_2^{(1)}\wedge e_2^{(2)}}_{=c\,\in\, C}.
\end{equation}
To explain $\Fil^1$, we note that
since $W=\wedge^2V$, the wedge product induces a non-degenerate
symmetric bilinear form $W\times W\to\wedge^4V\cong\QQ$.
Then, $\Fil^2$ is isotropic with respect to this pairing and it is easy to see that
$$
   \Fil^1 \,=\, (\Fil^2)^\perp.
$$
\end{enumerate}
Thus, we obtain a filtered $(\varphi,N)$-module 
$(W\otimes K_0,\Fil^\bullet,\varphi_W\otimes\sigma,N_W\otimes{\rm id}_{K_0})$
and a rational structure $(W=A\oplus B_0\oplus B_1\oplus C,\varphi_W,N_W)$.
We will see in Proposition \ref{prop: raskind for product tate} that it is isomorphic
to $\DD_{\rm st}(\Het{2}(X_{\overline{K}},\QQ_p))$.

\begin{Lemma}\label{lem: is raskind admissible}
 This just-constructed filtered $(\varphi,N)$-module together with its
 rational structure 
 is Raskind-admissible.
 More precisely, we have 
 $$
  \dim_\QQ \left(\Fil^1\cap B_1\right)\,=\, 
  \left\{ \begin{array}{ll}
    2 & \mbox{ if }\MazurL(q_1)\neq\MazurL(q_2) \mbox{ and }\\
    3 & \mbox{ if }\MazurL(q_1)=\MazurL(q_2).
  \end{array}\right.
 $$
\end{Lemma}

\prf
First, we note that an element lies in $\Fil^1$ if and only if it has zero intersection
with the vector \eqref{eq: Mazur vector}.
This makes computations very easy.

If $\MazurL(q_1)=\MazurL(q_2)$, then
$B_1\subset \Fil^1$ from which the claim on the dimension and 
Raskind-admissibility easily follows.

If $\MazurL(q_1)\neq\MazurL(q_2)$, then 
$b_2,b_3\in\Fil^1$ and thus, $\dim_\QQ(\Fil^1\cap B_1)\geq2$.
Moreover, we have that $\Fil^1\cap(B_1\otimes K)$ is strictly
contained in $B_1\otimes K$, which yields the chain of inequalities
$\dim_\QQ(\Fil^1\cap B_1)\leq\dim_{\QQ_p}(\Fil^1\cap(B_1\otimes\QQ_p))\leq
\dim_K(\Fil^1\cap(B_1\otimes K))\leq 2$.
Together with the previous inequality, this implies
that we have equality everywhere and
establishes the claimed dimension, as well as
Raskind-admissibility.
\qed\medskip

\begin{Proposition}\label{prop: raskind for product tate}
  Let $K$ be a $p$-adic field, let $q_i\in K^\times$ with $\nu_p(q_i)>0$ for $i=1,2$, let
   $E(q_1), E(q_2)$ be the associated Tate elliptic curves over $K$, and set
  $X:=E(q_1)\times E(q_2)$.
  \begin{enumerate}
  \item $X$ admits a proper and semi-stable model ${\cal X}\to\Spec\OO_K$,
   whose special fibre ${\cal X}_0$ is cohomogically totally degenerate.
  \item The filtered $(\varphi,N)$-module $\DD_{\rm st}(\Het{2}(X_{\overline{K}},\QQ_p))$
   together with the rational structure associated to ${\cal X}_0$ are 
   isomorphic to the one constructed at the beginning of this subsection.
  \end{enumerate}
  In particular, this filtered $(\varphi,N)$-module is Raskind-admissible and
  Conjecture \ref{conj: raskind} is true for $X$.
  More precisely, we have
$$
  \rho(X)\,=\, 
  \left\{ \begin{array}{ll}
    2 & \mbox{ if }\MazurL(q_1)\neq\MazurL(q_2) \mbox{ and }\\
    3 & \mbox{ if }\MazurL(q_1)=\MazurL(q_2).
  \end{array}\right.
 $$
 for the Picard rank of $X$.
\end{Proposition} 

\prf
Using the decomposition
$$
  \Het{2}(X_{\overline{K}},\QQ_p)\,\cong\, \bigwedge^2 \Het{2}(X_{\overline{K}},\QQ_p) \,\cong\,
  \bigwedge^2\left( \Het{1}(E(q_1))_{\overline{K}},\QQ_p)\oplus\Het{1}(E(q_2)_{\overline{K}},\QQ_p) \right)
$$
and the explicit description of $\DD_{\rm st}(\Het{1}(E(q_i)_{\overline{K}},\QQ_p))$ given
at the beginning of Section \ref{subsec: product tate},
it is easy to see that the filtered $(\varphi,N)$-module constructed above
is isomorphic to $\DD_{\rm st}(\Het{2}(X_{\overline{K}},\QQ_p))$.

Let ${\cal E}_i\to\Spec\OO_K$ the ``standard'' proper and semi-stable model
of $E_i:=E(q_i)$, whose special fibre ${\cal E}_{i,0}$ is a cycle of $\PP^1$'s.
It follows that $\HdR{\ast}(E_i/K)$ and
$\Hlogcris{\ast}({\cal E}_{i,0}/K_0)$
carry structures of $\QQ$-vector spaces arising from 
classes of algebraic cycles of the special fibre.
Tensoring with $K$, we obtain the filtered 
$(\varphi,N)$-module $V^{(i)}\otimes K$ constructed at the beginning of
Section \ref{subsec: product tate}.
Conversely, the vector $e_2^{(i)}\in V^{(i)}$ arises as 
an eigenvector of $\varphi$, which makes it canonical up to 
a factor $\lambda_i\in\QQ_p^\times$.
Since $e_1^{(i)}=N(e_2^{(i)})$, this also determines $e_1^{(i)}$.
Thus, given $V^{(i)}\otimes K$, the just discussed rational structure must be the
$\QQ$-span $\langle\lambda_i e_1^{(i)},\lambda_i e_2^{(i)}\rangle$
for some $\lambda_i\in\QQ_p^\times$.

We obtain a proper and semi-stable model ${\cal X}\to\Spec\OO_K$ of $X$
via a partial resolution of singularities of ${\cal E}_1\times{\cal E}_2\to\Spec\OO_K$.
Thus, the  $\QQ$-vector space structures on $\HdR{\ast}(X/K)$ and
$\Hlogcris{\ast}({\cal X}_{0}/K_0)$ arise (via K\"unneth) from the 
$\QQ$-vector space structures on $\HdR{\ast}(E_i/K)$ and
$\Hlogcris{\ast}({\cal E}_{i,0}/K_0)$.
Therefore, the rational structure on $\Hlogcris{\ast}({\cal X}_{0}/K_0)$
arising from ${\cal X}_0$ is of the form 
as discussed at the beginning of Section \ref{subsec: product tate}
and isomorphic to it, where the isomorphism is induced
by multiplication by scalars $\lambda_i\in\QQ_p^\times$ 
as just explained.
We note that such a rescaling  multiplies the vector
\eqref{eq: Mazur vector} by $\lambda_1\cdot\lambda_2$, that is,
it still spans the same $K$-vector space $\Fil^2$.
Since $\Fil^1=(\Fil^2)^\perp$, we see that rescaling also leaves $\Fil^1$
invariant.
Thus, the filtered $(\varphi,N)$-module together with its rational structure
discussed at the beginning of Section \ref{subsec: product tate}
is isomorphic to $\DD_{\rm st}(\Het{2}(X_{\overline{K}},\QQ_p))$ with its
rational structure arising from ${\cal X}_0$.

By Lemma \ref{lem: is raskind admissible}, it is Raskind-admissible
and thus, Conjecture \ref{conj: raskind} for $X$ follows from
Theorem \ref{thm: raskind translation}.
The claim on the Picard ranks follows from the dimensions computed
in Lemma \ref{lem: is raskind admissible}.
\qed\medskip

\begin{Remark}
 We note that Conjecture \ref{conj: raskind} for the product of two Tate curves
 was already established Tate, as is explained in \cite[Appendix A.1.4]{Serre}  (see also \cite[\S 4.1, Corollary 19]{RXJacobians}), and it also follows from the more general Corollary \ref{cor: product tate}.
\end{Remark} 

The previous proposition raises the question, whether 
every admissible $(\varphi,N)$-module with rational structure is 
Raskind-admissible.
This is \emph{not} the case, as the following example shows.

\begin{Example}
 \label{ex: non-Raskind-admissible}
 We keep the notations and assumptions of Proposition \ref{prop: raskind for product tate}.
 Choose $\gamma\in\QQ_p\backslash\QQ$, choose
 $\lambda\in K\backslash\{0\}$, and define
 $$
    v \,:=\, 2(\lambda^2-\gamma) a \,+\, \lambda b_0 \,-\, \gamma b_2 \,+\, b_3 \,+\,c.
 $$
 Then, $v$ is isotropic with respect to the pairing introduced at the beginning
 of Section \ref{subsec: product tate}, which allows us to define a filtration
 on $V\otimes K$ via
 $$
     \Fil^2 \,:=\, K\cdot v \,\mbox{ \quad and \quad } \Fil^1 \,:=\, (\Fil^2)^\perp\,.
 $$
 This filtration is  admissible by 
 Proposition \ref{prop: classification admissible} below.
 We leave it to the reader to check that 
 $$
  \begin{array}{lcll}
     \Fil^1\cap (B_1\otimes\QQ_p) &=& \QQ_p\langle b_1, \gamma b_2+b_3 \rangle \\ 
     \Fil^1\cap B_1 &=& \QQ\,\langle b_1 \rangle, 
  \end{array}
 $$
 which implies that $\Fil^\bullet$ is \emph{not} Raskind-admissible.
 Varying $\gamma$ and $\lambda$, we obtain a whole family of such modules.
\end{Example}

\subsection{Admissibility}
Given a $(\varphi,N)$-module and a rational structure
in the sense of Definition \ref{def: rational structure} 
with $\dim A=\dim C=1$, we now address the question
when a filtration over $K$ is ordinary in the sense of 
Perrin-Riou \cite[1.2]{Perrin-Riou}.
We have the following result, which should be 
the framework for rational structures
on $\DD_{\rm st}(\Het{2}(X_{\overline{K}},\QQ_p))$
where $X$ is a smooth and proper surface over $K$ with cohomological
total degeneration and 
$h^{0,2}=h^{2,0}=1$ (that is, $p_g=1$ in classical terminology).

\begin{Proposition}\label{prop: classification admissible} 
 Let $(V=A\oplus B_0\oplus B_1\oplus C,\varphi_V,N_V)$
 be a rational structure in the sense of Definition \ref{def: rational structure}.
 Assume moreover that
 \begin{enumerate}
  \item $\dim A=\dim C=1$.  
    We also fix a non-zero element $c\in C$, that is, a basis of
    this vector space.
  \item There exists a non-degenerate, symmetric, and bilinear
pairing $Q:V\times V\to\QQ$.
 \end{enumerate}
 Let $K$ be a $p$-adic field and let $\Fil^\bullet$ be a filtration
 on $V\otimes K$ with $\Fil^0=V\otimes K$ and $\Fil^3=0$.
 \begin{enumerate}
  \item If $\Fil^\bullet$ is ordinary, then $\Fil^2$ is $1$-dimensional and
  spanned by a unique vector of the form
  \begin{equation} 
  \label{eq: parameter vector}
     v\,=\,  v'\,+\, c
  \end{equation}
  with $v'\in (A\oplus B_0\oplus B_1)\otimes K$.
  \item If $v$ is of the form \eqref{eq: parameter vector} with $Q(v,v)=0$
   and $Q(v, N^2(c))\neq0$, then 
  $$
     \Fil^2 \,:=\, K\cdot v\mbox{ \quad and \quad } \Fil^1:=(\Fil^2)^\perp
  $$
  defines an ordinary filtration on $V\otimes K$.
 \end{enumerate}
In particular, there exists a bijection
  $$  \small
 \left\{ v\in V\otimes K \,:\, \begin{array}{l}
  \mbox{$Q(v,v)=0$, $Q(v,N^2(c))\neq0$, }\\
  \mbox{and $v$ is of the form \eqref{eq: parameter vector}}
  \end{array}
  \right\}
   \,\to\,
   \left\{
   \begin{array}{c} 
  \mbox{ordinary filtrations}\\
  \mbox{of $V\otimes K$ with}\\
  \Fil^1=(\Fil^2)^\perp
 \end{array}
  \right\}
  $$
that is defined by sending $v$ to the filtration with 
$\Fil^1=(Kv)^\perp$ and $\Fil^2=Kv$.
\end{Proposition} 

\prf
By \cite[2.6]{Perrin-Riou} it follows that $\Fil^\bullet$ is ordinary 
if and only if we have direct sum decompositions
$$
  V\otimes K \,=\, \Fil^1 \,\oplus\, \left(A\otimes K\right) \,=\, \Fil^2 \,\oplus\, \left((A\oplus B)\otimes K\right),
$$
where $B=B_0\oplus B_1$,

In particular, $\Fil^2$ is one-dimensional and thus, generated by one element $v\in V\otimes K$.
Being ordinary, we have $v\not\in (A\oplus B)\otimes K$.
Thus, after possibly rescaling, we may assume that $v$ is of the form
$v'+c$ for some $v'\in (A\oplus B)\otimes K$ and this vector is unique.
This establishes claim (1).

If $v\in V\otimes K$ satisfies $Q(v,v)=0$, 
then $\Fil^2:=Kv$ is contained in $\Fil^1:=(Kv)^\perp$, and we obtain a filtration.
Since $Q$ is non-degenerate, it follows that $(Kv)^\perp$ is of codimension $1$ in $V\otimes K$.
If $v$ is moreover of the form \eqref{eq: parameter vector}, that is, $v=v'+c$ with 
$v'\in (A\oplus B_0\oplus B_1)\otimes K$, then $V\otimes K = \Fil^2\oplus (A\oplus B)\otimes K$.
If $Q(v,N^2(c))\neq0$, then $A$, which is spanned by $N^2(c)$, is not contained in
$(Kv)^\perp=\Fil^1$, that is, we have $V\otimes K = \Fil^1\oplus A\otimes K$.
This shows that $\Fil^\bullet$ is ordinary and establishes claim (2).

We leave the remaining assertion to the reader.
\qed\medskip

Thus, ordinary filtrations on $V\otimes K$ with
$\Fil^1=(\Fil^2)^\perp$ are parameterised by the $K$-rational points of a 
quasi-affine scheme over $\QQ$.
Explicitly: first, $b_0:=N(c)$ is a basis of $B_0$, then, $c:=N(b_0)=N^2(c)$ is a basis
of $A$ and we choose a basis $b_1,...,b_s$ of $B_1$.
Then, we define an affine quadric  $\cal Q$ by the equation
$$
  Q\left(\lambda a +\sum_{i=0}^s\mu_i b_i+c,\, \lambda a +\sum_{i=0}^s\mu_i b_i+c\right) \,=\,0
$$
in the $(s+2)$-dimensional affine space with coordinates
$(\lambda,\mu_0,...,\mu_s)$.
We define the Zariski open subset ${\cal U}\subset{\cal Q}$ 
by the condition 
$$
  Q\left(\lambda a +\sum_{i=0}^s\mu_i b_i+c,\, N^2(c)\right) \,\neq\,0.
$$
Then, $\cal U$ is a quasi-affine scheme of dimension $(s+1)$ over $\QQ$,
whose $K$-rational points are in bijection to ordinary
filtrations $\Fil^\bullet$ on $V\otimes K$ 
with $\Fil^1=(\Fil^2)^\perp$.

\section{A counter-example}\label{sec: counterexample}

For abelian varieties that are isogenous to products of Tate elliptic curves, 
we established Raskind's conjecture for divisors (Conjecture \ref{conj: raskind})
and showed that Question \ref{question: tate} has a positive answer in the previous section.
In this section, we will show that in general, Question \ref{question: tate} has a negative answer 
and that in general Raskind's conjecture for divisors is false.

\subsection{Totally degenerate reduction and $\ell= p$}
In view of Proposition \ref{prop: isogeny product tate}, the first place to look for counter-examples
are abeloid surfaces over $p$-adic fields.

\begin{Theorem}\label{thm: counterexample}
 Let $p$ be a prime with $p\geq5$ and $p\equiv 1\mod3$.
 Then, there exists a Tate elliptic curve $A$ and an algebraisable abeloid surface $B$
 over $\QQ_p$, such that 
 \begin{enumerate}
 \item
 the natural maps
 \begin{eqnarray*}
    \Hom(A,B)\otimes_\ZZ\ZZ_p &\to& \Hom_{G_{\QQ_p}}\left( T_p(A),\, T_p(B) \right) \\
     \End(B)\otimes_\ZZ\ZZ_p &\to& \End_{G_{\QQ_p}}\left( T_p(B) \right)
 \end{eqnarray*}
 are not surjective.
 In particular, \eqref{tatehom2} (please see \S\ref{abelian and abeloid}) is not surjective and Question \ref{question: tate}
 has a negative answer for $\ell=p$.
 \item The natural map induced by the first Chern class map
 \begin{eqnarray*}
     \Pic(B)\otimes_\ZZ\QQ_p &\to& \Het{2}( B_{\overline{\QQ}_p},\QQ_p(1))^{G_{\QQ_p}}
 \end{eqnarray*}
 is not surjective.
 In particular, \eqref{tatehom} is not surjective for $\ell=p$ 
 and Raskind's conjecture for divisors (Conjecture \ref{conj: raskind}) is false.
 \end{enumerate}
\end{Theorem}

\begin{Remark}
As seen in Proposition \ref{prop: abelian totally degenerate}, an algebraisable abeloid variety over $K$ 
is the same as an abelian variety over $K$ with totally degenerate reduction.
\end{Remark}

\prf
Let $\varepsilon\in 1+p\cdot\ZZ_p$ be a non-trivial $p$-adic unit, set $q_1:=p$, and set $q_2:=\varepsilon\cdot p$.
Let $A:=E(q_1)$ be the Tate elliptic curve over $\QQ_p$ with respect to the lattice $q_1^\ZZ\subset\QQ_p^\times$.
First, we set
$$
V_B'\,:=\,
\begin{pmatrix}
q_1 & 1 \\
1 & q_2 \\
\end{pmatrix}
\,\in\, {\rm Mat}_{2\times 2}(\QQ_p^\times)
$$
and note that abeloid surface $B'$ over $\QQ_p$ associated to the matrix $V_B'$ 
is isomorphic to the product of the two Tate elliptic curves $E(q_1)\times E(q_2)$.

For $v_1,v_2\in\ZZ_p$, we define
$$
S\,:=\,{\rm Id}_{2\times2}\,-\,2\cdot
\begin{pmatrix}
v_1\\ v_2
\end{pmatrix}
\cdot
\begin{pmatrix}
v_1 & v_2 
\end{pmatrix}
\,=\,
\begin{pmatrix}
1-2v_1^2 & -2v_1v_2 \\
-2v_1v_2 & 1-2v_2^2 \\
\end{pmatrix}
\,\in\, {\rm Mat}_{2\times 2}(\ZZ_p)\,.
$$
Clearly, $S$ is symmetric and if $v_1^2+v_2^2=1$, which we will assume from now on,  
then $S^{-1}=S^t$ and thus, we find
$$
   S \,=\, \begin{pmatrix} 
    a & b \\ 
    b & -a 
   \end{pmatrix}
   \mbox{ \qquad with \qquad }
   a:=1-2v_1^2, \mbox{ \qquad } b:=-2v_1v_2,
$$
and note that we have $a^2+b^2=1$
(one should think of this matrix as the analog of an
orthogonal matrix over $\RR$ that describes the reflexion along  the axis spanned by $(v_1\, v_2)$). 
We set
$$
  V_B \,:=\, S^{-1}\odot V'_B\odot S \,\in\, {\rm Mat}_{2\times 2}(\widehat{\QQ_p^\times}^p),
$$
(here $\widehat{\QQ_p^\times}^p$ denotes the $p$-adic completion $\gamma_p$
from Lemma \ref{lemma: Serre}), which is equal to 
$$
V_B\,=\, 
 \begin{pmatrix}
  q_1^{a^2}\cdot q_2^{b^2} & q_1^{ab}\cdot q_2^{-ab} \\
  q_1^{ab}\cdot q_2^{-ab} & q_1^{b^2}\cdot q_2^{a^2}
\end{pmatrix}
\,=\,
 \begin{pmatrix}
  \varepsilon^{b^2}\cdot p & \varepsilon^{-ab} \\
  \varepsilon^{-ab} & \varepsilon^{a^2}\cdot p
\end{pmatrix}\,.
$$
In particular, this matrix has actually coefficients in $\QQ_p^\times$ rather
than merely in $\widehat{\QQ_p^\times}$, see also Remark \ref{rem: gamma}.
Moreover, the valuation
$\nu_p:\QQ_p^\times\to\QQ$ sends 
$V_B$ to the identity matrix.
Since $V_B$ is a symmetric matrix and since it is definite with respect to the
$\QQ$-linear functional $\nu_p$, it is a Riemann matrix in the 
sense of Gerritzen \cite{GerritzenRiemann}.
We let $B$ be the abeloid surface over $\QQ_p$ associated
to $V_B$.
By \cite[Theorem 11]{GerritzenRiemann}, this surface is algebraisable,
that is, $B$ is an abelian surface over $\QQ_p$.

In order to determine $\Hom(A,B)$ and $\Hom_{G_{\QQ_p}}(T_p(A),T_p(B))$,
we are looking at the equation
$$
 \resizebox{1.0\hsize}{!}{$
 \begin{pmatrix}
 q_1 
 \end{pmatrix}
 \,\odot\,
 \begin{pmatrix}
 x & y
 \end{pmatrix}
 \,=\,
  \begin{pmatrix}
 p^x & p^y
 \end{pmatrix}
 \,\stackrel{!}{=}\,
 \begin{pmatrix}
 x' & y'
 \end{pmatrix}
 \,\odot\,V_B
 \,=\,
 \begin{pmatrix}
 \varepsilon^{x'b^2-y'ab}\cdot p^{x'}
 & \varepsilon^{y'a^2-x'ab}\cdot p^{y'}
 \end{pmatrix}\,.$}
$$
Taking valuations, we find $x=x'$ and $y=y'$.
To avoid trivialities, we will also assume that $a\neq 0$ and $b\neq0$,
which leads to the general solution
$$
  y\,=\, \frac{b}{a}\cdot x\,.
$$
Since there is always the $p$-adic solution $x=a$, $y=b$,
Proposition \ref{prop: tate module} implies that
$\Hom(V_p(A),V_p(B))$ is non-zero and in fact, isomorphic to $\QQ_p$.
On the other hand, Theorem \ref{thm: abeloid morphisms} implies that
$\Hom(A,B)\otimes\QQ$ is non-zero 
if and only we can find a solution with $x\in\QQ$ and $y\in\QQ$.
Thus, $\Hom(A,B)\otimes\QQ$ is non-zero if and only if
$a/b\in\QQ$.

In order to establish the first claim, we have to show that we can 
find $v_1,v_2\in\ZZ_p$ that satisfy all the restrictions we made during the previous
discussion and such that $a/b=(2v_1^2-1)/(2v_1v_2)$ does not lie in $\QQ$.
If $v_1=2$, then $v_2=\sqrt{-3}$ does not lie in $\QQ$, but it is an element of $\ZZ_p$
if $p\geq5$ and $p\equiv1\mod3$ (the last statement easily follows from the law of
quadratic reciprocity).

Similarly, to compute the endomorphism algebras $\End(B)\otimes\QQ$ and $\End_{G_{\QQ_p}}(V_p(B))$,
we have to solve the equation
$$
  V_B\,\odot\, C \,=\, C'\,\odot\, V_B
$$
for $2\times2$-matrices $C$ and $C'$ with entries in $\QQ$ and $\QQ_p$, respectively.
Taking valuations, we find $C=C'$.
Moreover, if $C=(c_{ij})_{1\leq i,j\leq 2}$, then we leave it to the reader to show that 
the general solution is given by
$$
   c_{12}\,=\,c_{21}\mbox{ \quad and \quad }
   c_{22}-c_{11}\,=\, c_{12}\cdot \frac{b^2-a^2}{ab}.
$$
We thus always have the solution $c_{12}=c_{21}=0$ and $c_{11}=c_{22}$, that is,
multiplication by a scalar.
If $(b^2-a^2)/(ab)$ does not lie in $\QQ$ (which is the case if $v_1=2$ and $v_2=\sqrt{-3}$), 
then these are the only solutions in $\QQ$,
that is, $\End(B)\otimes\QQ\cong\QQ$.
On the other hand, the above equation always has more solutions in $\QQ_p$,
that is, $\End(V_p(B))$ is strictly larger than $\QQ_p$.

These computations establish the first claim.
The second claim follows from the first claim by the same arguments 
as in the proof of Proposition \ref{cor: product tate}.
\qed\medskip

\begin{Remark}
The restrictions on the prime $p$ in Theorem \ref{thm: counterexample} are artificial in that 
they are only used to state a clean counterexample. 
The method of proof should give counterexamples to Conjecture \ref{conj: raskind} and Question \ref{question: tate} 
for any prime $p$. 
\end{Remark}

\appendix

\section{Further (counter-)examples}
\label{appendix}

So far, we have studied the surjectivity of the maps \eqref{tatehom} and
\eqref{tatehom2} (defined in \S\ref{intro}  )  in the case where $\ell=p$ and where 
the varieties in question are smooth and proper
over a $p$-adic field with totally degenerate reduction.
If $\ell\neq p$ or if the variety has good reduction,
then it was already more or less well-known to the experts 
that one cannot hope for such surjectivity results, but for the sake
of completing the picture, 
we have decided to collect some examples.
As a byproduct, we see that also ``independence of $\ell$''
fails.
In this section, we claim only little originality.

\subsection{Good reduction and $\ell\neq p$}
Let $X$ be a smooth and proper variety over a $p$-adic field $K$ that
has good reduction, say, via a smooth and proper model
${\cal X}\to\Spec\OO_K$ with special fibre ${\cal X}_0$ over the residue field $k$.
By base-change in \'etale cohomology, the 
$G_K$-action on $\Het{*}(X_{\overline{K}},\QQ_\ell)$ is unramified
and factors through the $G_k$-action 
on $\Het{*}({\cal X}_{0,\overline{k}},\QQ_\ell)$.
\begin{enumerate}
 \item In particular, we have
   $$
   \Het{2}({X}_{\overline{K}},\QQ_\ell(1))^{G_K} \,=\, 
   \Het{2}({\cal X}_{0,\overline{k}},\QQ_\ell(1))^{G_k}
 $$
 and if ${\cal X}_0$ satisfies the Tate conjecture for divisors over the finite field $k$,
 then these spaces are isomorphic to $\Pic({\cal X}_0)\otimes_\ZZ\QQ_\ell$.
\item Let $A$ and $B$ be abelian varieties with good reduction over $K$ and let
 let ${\cal A}_0$ and ${\cal B}_0$ be the special fibres of their N\'eron models.
 In this case, base-change implies that 
 $$
    {\rm Hom}_{G_K}\left(T_\ell(A),T_\ell(B)\right) \,=\,
    {\rm Hom}_{G_k}\left(T_\ell({\cal A}_0),T_\ell({\cal B}_0)\right),
 $$
 which is isomorphic to ${\rm Hom}({\cal A}_0,{\cal B}_0)\otimes_\ZZ\ZZ_\ell$
 by Tate's theorem \cite{TateAbelian}. 
\end{enumerate}
In particular, the right hand sides of \eqref{tatehom} and \eqref{tatehom2}
compute invariants of the special fibre, see also Remark \ref{remark on l vs p}.
After these preparations, it is easy to give the desired counter-examples - in fact,
``almost all'' elliptic curves provide counter-examples.

\begin{Proposition}
 Let $K$ be a $p$-adic field, let $E$ be an elliptic curve over $K$ with good reduction, 
 and let ${\cal E}_0$ be the special fibre of its N\'eron model.
 If ${\cal E}_0$ is supersingular or if ${\cal E}_0$ is ordinary
 and $E$ does not have CM, then
 \eqref{tatehom} is not surjective for $X=E\times E$ and $\ell\neq p$
 and \eqref{tatehom2} is not surjective for $A=B=E$ and $\ell\neq p$.
\end{Proposition}

\prf
The first claim follows from the second claim by the same arguments 
as in the proof of Proposition \ref{cor: product tate}.
Therefore, it suffices to show that the natural map 
$\End(E)\otimes_\ZZ\QQ_\ell\to \End_{G_K}(V_\ell(E))$
is not surjective.

First, assume that ${\cal E}_0$ is supersingular.
Then, $\End({\cal E}_0)$ is an order in a quaternion algebra
and thus, $\End({\cal E}_0)\otimes\QQ_\ell$ is $4$-dimensional.
Hence, $\End_{G_K}(V_\ell(E))$, which is isomorphic to
$\End_{G_k}(V_\ell({\cal E}_0))$, is $4$-dimensional 
by Tate's theorem \cite[Main Theorem]{TateAbelian}.
On the other hand, $\End(E)$ is isomorphic to 
$\ZZ$ or to an order in a quadratic imaginary field, which implies
that $\End(E)\otimes\QQ_\ell\to\End_{G_K}(V_\ell(E))$ cannot
be surjective.

Similarly, if $E$ does not have CM, then $\End(E)\otimes\QQ_\ell=\QQ_\ell$.
Moreover, if ${\cal E}_0$ is ordinary, then 
$\End_{G_k}({\cal E}_0)$ is an order in a quadratic imaginary
field and $\End_{G_k}({\cal E}_0)\otimes\QQ_\ell$
is $2$-dimensional.
\qed\medskip

\subsection{Good reduction and $\ell= p$}
Next, we show that Conjecture \ref{conj: raskind} and Question \ref{question: tate}
have a negative answer if $\ell=p$ and in the case of good reduction.

\begin{Proposition}
 Let $K$ be a $p$-adic field
 and let ${\cal E}_0$ be an ordinary elliptic curve over $k$.
 \begin{enumerate}
  \item Let $A$ be a lift of ${\cal E}_{0}$ over $K$ with CM, for example, the canonical lift, and
  \item let $B$ be a lift of ${\cal E}_{0}$ over $K$ without CM.
 \end{enumerate} 
 Then, \eqref{tatehom} is not surjective for $X=A\times B$ and $\ell=p$
 and \eqref{tatehom2} is not surjective for $A$, $B$, and $\ell=p$.
\end{Proposition}

\prf
The first claim follows from the second claim by the same arguments 
as in the proof of Proposition \ref{cor: product tate}.
Therefore, it suffices to show that the natural map 
${\rm Hom}(A,B)\otimes_\ZZ\QQ_p\to \Hom_{G_K}(V_p(A),V_p(B))$
is not surjective.

Since $A$ cannot be isogenous to 
$B$, the source $\Hom(A,B)\otimes_{\mathbb{Z}}\mathbb{Q}_{p}$ is trivial.
Next, for $E\in\{A,B\}$ there exists a short exact sequence of $p$-adic 
$G_K$-representations
\begin{equation}
 \label{eq: ordinary}
     0\,\to\, X \,\to\, V_p(E) \,\to\, Y \,\to\, 0.
\end{equation}
More precisely, $X$ corresponds to the Tate module associated to
the connected component of the $p$-divisible group ${\cal E}_0[p^\infty]$ 
and $Y$ corresponds to the Tate module associated to the \'etale
quotient.
In particular, the $G_K$-representations $X$ and $Y$ only depend on
${\cal E}_0$ and not on the choice of lift $E$.
Moreover, the sequence of $G_K$-representations \eqref{eq: ordinary}
splits if and only if the lift of ${\cal E}_0$ has CM, that is, 
it splits for $A$ but not for $B$.
We refer to \cite[Appendix A.2.4]{Serre} for
details and proof. 
But this implies that the target 
$\Hom_{G_{K}}(V_p(A),V_p(B))$  is non-trivial:
taking the monomorphism $X\to V_p(B)$ from \eqref{eq: ordinary}
and the zero map $Y\to V_p(B)$, we obtain a non-trivial
and $G_{K}$-equivariant map
$$
V_p(A) \,=\, X\oplus Y \,\to\, V_p(B).
$$
This establishes the second claim.
\qed\medskip

\begin{Remark}
In \cite[\S3.5]{LT66}, Lubin and Tate constructed elliptic curves $E$ over $p$-adic fields having good 
and supersingular reduction such that the the monomorphism
$$
\End(E)\otimes_{\ZZ}\ZZ_p \,\hookrightarrow\, \End_{G_{K}}\left(T_p(E)\right)
$$
is not surjective. 
Therefore, \eqref{tatehom} is not surjective for $X=E\times E$ and $\ell=p$ and
\eqref{tatehom2} is not surjective for $A=E$, $B=A$, and $\ell=p$.
\end{Remark}

\subsection{Totally degenerate reduction and $\ell\neq p$}

Now, we show that Proposition \ref{prop: isogeny product tate}
and Corollary \ref{cor: product tate} fail if $\ell\neq p$.

\begin{Proposition}\label{prop: totally deg ell}
 For every prime $p$, there exist Tate elliptic curves $A$ and $B$ over $\QQ_p$, 
 such that
 \eqref{tatehom} is not surjective for $X=A\times B$ and $\ell\neq p$
 and \eqref{tatehom2} is not surjective for $A$, $B$, and $\ell\neq p$.
\end{Proposition}

\prf
Let $\varepsilon\in 1+p\cdot \ZZ_p$ be a non-trivial $p$-adic unit
and let $A:=E(p)$ and $B:=E(\varepsilon\cdot p)$ be the Tate elliptic curves
associated to $p\in\QQ_p^\times$ and $\varepsilon\cdot p\in\QQ_p^\times$. 
Then, $A$ and $B$ are not isogenous by Serre's criterion 
(Theorem \ref{thm: LeStum Serre}.(2)) and thus, $\Hom(A,B)=0$.

We let $\gamma_\ell$ be the $\ell$-adic
completion from Lemma \ref{lemma: Serre}.
By Proposition \ref{prop: tate module}, we have 
$$
\Hom_{G_K}(T_\ell(A),T_\ell(B)) \,\cong\, \left\{ m\in\ZZ_\ell \,|\, \exists n\in\ZZ_\ell\,:\,
\gamma_\ell(p)^m=\gamma_\ell(\varepsilon\cdot p)^n \right\} \,.
$$
By Lemma \ref{lemma: Serre} and Remark \ref{rem: gamma}, we have 
$\gamma_\ell(\varepsilon)=1$, that is, we 
have $\gamma_\ell(p)=\gamma_\ell(\varepsilon\cdot p)$, which implies
that $\Hom_{G_K}(T_\ell(A),T_\ell(B))$ is non-zero
(in fact, isomorphic to $\ZZ_\ell$).
\qed\medskip

\subsection{\'Etale fundamental groups and all primes at once}
If $A$ is an abelian variety of dimension $g$ over a field characteristic zero field $F$, 
then there exists an isomorphism of $G_F$-representations
$$
   \piet(A_{\overline{F}}) \,\cong\, \prod_\ell\, T_\ell(A),
$$ 
where the product is taken over all primes.
As an abelian group, this is a free $\widehat{\ZZ}$-module of rank $2g$.
Now, if $A$ and $B$ are abelian varieties over $F$, then one may ask whether the 
natural map
\begin{equation*}\tag{$\star\star\star$}
\label{tatehom3}
   \Hom(A,B)\otimes_\ZZ\widehat{\ZZ}\,\to\,{\rm Hom}_{G_F}\left( \piet(A_{\overline{F}}),\, \piet(B_{\overline{F}})\right)
\end{equation*}
is an isomorphism.

\begin{Proposition}
 For every prime $p$, there exist Tate elliptic curves $A$ and $B$ over $\QQ_p$, 
 such that \eqref{tatehom3} is not surjective.
\end{Proposition}

\prf
Fix a prime $\ell_0\neq p$ and let $A$ and $B$ be counterexamples 
as provided by Proposition \ref{prop: totally deg ell}, that is,
\eqref{tatehom2} is not surjective for $A$, $B$, and $\ell_0$.
Since \eqref{tatehom3} factors through
$$
\prod_{\ell} \Hom(A,B)\otimes_{\ZZ}\ZZ_\ell \,\to\, \prod_{\ell}\Hom_{G_{F}}(T_{\ell}(A),T_{\ell}(B))
$$
and since this map is not surjective at the factor corresponding to $\ell_0$,
the claim follows.
\qed\medskip

\subsection{Independence of $\ell$}
From the previous computations, we conclude that
also ``independence of $\ell$'' fails in the $p$-adic 
world, see also Remark \ref{remark on l vs p} and the subsequent discussion.

\begin{Proposition}
  For every prime $p$, there exist Tate elliptic curves $A$ and $B$ over $\QQ_p$, 
 such that
 $$
    \dim_{\QQ_\ell} \Het{2}((A\times B)_{\overline{\QQ}_p},\QQ_\ell(1))^{G_{\QQ_p}}
    \,=\,\left\{
    \begin{array}{cl}
    2 & \mbox{ \quad if }\ell=p\\
    3 & \mbox{ \quad if }\ell\neq p
    \end{array}
    \right. \,.
 $$
 In particular, this dimension depends on the prime $\ell$.
\end{Proposition}

\prf
Let $A$ and $B$ the Tate elliptic curves from the proof of Proposition \ref{prop: totally deg ell}.
There, we have seen that  $A$ and $B$ are not isogenous and that
$\Hom_{G_{\QQ_p}}(V_\ell(A),V_\ell(B))$ is one-dimensional if $\ell\neq p$.
On the other hand, since $A$ and $B$ are not isogenous,
$\Hom_{G_{\QQ_p}}(V_p(A),V_p(B))$ is zero by Proposition \ref{prop: isogeny product tate}.

The arguments from the proof of Proposition \ref{cor: product tate}
show that the sought $\QQ_\ell$-dimensions are equal to 
$2+\dim_{\QQ_\ell} \Hom_{G_K}(V_\ell(A),V_\ell(B))$
and the claim follows.
\qed

\end{document}